\documentclass[12pt,centertags]{amsr}  % Specially designed style file amsr.sty
\usepackage{amsthm}
\usepackage{amssymb}        % Useless symbols
\usepackage{calc}           % Basic arithmetic with lengths...
\usepackage{xspace}         % Smart space for shortcuts
\usepackage[all]{xy}        % Commutative Diagrams
\CompileMatrices            % Faster Tex-ing
\UseTips                    % Use Computer Modern Arrowheads

\usepackage[bookmarks=true]{hyperref}   % Hyperref in DVI and PDF (like HTML Links between sections)

\makeindex
\newcommand{\idx}[1]{#1\index{#1}}

\newcommand{\enm}[1]{\ensuremath{#1}}          % Shortcuts
\newcommand{\op}[1]{\operatorname{#1}}
\newcommand{\cal}[1]{\mathcal{#1}}
\newcommand{\wt}[1]{\widetilde{#1}}
\newcommand{\wh}[1]{\widehat{#1}}
\renewcommand{\bar}[1]{\overline{#1}}

\newcommand{\CC}{\enm{\mathbb{C}}}             % All Number domains easily accssable
\newcommand{\NN}{\enm{\mathbb{N}}}

\newcommand{\QQ}{\enm{\mathbb{Q}}}
\newcommand{\ZZ}{\enm{\mathbb{Z}}}
\newcommand{\FF}{\enm{\mathbb{F}}}
\renewcommand{\AA}{\enm{\mathbb{A}}}

\newcommand{\Aa}{\enm{\cal{A}}}           % All caligraphy letters easily accessable
\newcommand{\Bb}{\enm{\cal{B}}}

\newcommand{\Dd}{\enm{\cal{D}}}

\newcommand{\Hh}{\enm{\cal{H}}}

\newcommand{\Kk}{\enm{\cal{K}}}
\newcommand{\Ll}{\enm{\cal{L}}}
\newcommand{\Mm}{\enm{\cal{M}}}
\newcommand{\Nn}{\enm{\cal{N}}}
\newcommand{\Oo}{\enm{\cal{O}}}

\renewcommand{\phi}{\varphi}        % Dont know how to not loose the original ones???
\renewcommand{\theta}{\vartheta}
\renewcommand{\epsilon}{\varepsilon}

\newcommand{\Ann}{\op{Ann}}         % Standard Operators

\newcommand{\Spec}{\op{Spec}}

\newcommand{\Hom}{\op{Hom}}
\newcommand{\Ext}{\op{Ext}}

\newcommand{\End}{\op{End}}
\newcommand{\Aut}{\op{Aut}}

\newcommand{\height}{\op{ht}}

\newcommand{\id}{\op{id}}
\newcommand{\dirlim}{\varinjlim}
\newcommand{\invlim}{\varprojlim}
\newcommand{\length}{\op{length}}
\newcommand{\Image}{\op{Im}}

\newcommand{\tensor}{\otimes}         % Symbols with meaning
\newcommand{\dirsum}{\oplus}

\newcommand{\intsec}{\cap}
\newcommand{\Union}{\bigcup}

\newcommand{\set}[1]{\left\{#1\right\}}
\newcommand{\xn}[1][x]{\enm{#1_1,\ldots,#1_n}}
\newcommand{\xd}[1][x]{\enm{#1_1,\ldots,#1_d}}
\newcommand{\xc}[1][x]{\enm{#1_1,\ldots,#1_c}}

\renewcommand{\to}[1][]{\xrightarrow{\ #1\ }}
\newcommand{\To}{\xrightarrow}
\newcommand{\onto}{\twoheadrightarrow}
\newcommand{\into}{\hookrightarrow}

\newcommand{\diff}[1][x]{\textstyle{\frac{\partial}{\partial #1}}}

\newcommand{\usc}[1][m]{\underline{\phantom{#1}}}

\newcommand{\note}[1]{\marginpar{\small #1}}

\newcommand{\defeq}{\stackrel{\scriptscriptstyle \op{def}}{=}}
\renewcommand{\note}[1]{}    %ignores marginal notes (for debugging to avoid underfull h-boxes

\newcommand{\Czech}{\v{C}ech\xspace}               % The Czech guy
\newcommand{\Cz}[1]{\text{\v{C}}(#1)}              % and his complex
\newcommand{\modR}{{\text{mod--}R}}                % Right R--modules
\renewcommand{\mod}{\enm{\op{mod}}\xspace}         % To make sure R-mod is not typed cursively in Theorems
\newcommand{\Rmod}[1][R]{\enm{#1\text{--}\op{mod}}} % R--mod inside math environment (good for Functors T: R--mod --> A--mod

\newcommand{\F}[1][]{F^{#1*}}                  % Frobenius operator
\newcommand{\D}[1][]{D^{(#1)}}                 % Differential operators of level e
\newcommand{\Hhom}{\Hh om}                     % Curly Hom
\newcommand{\Fful}[1]{F^{\infty}{#1}}
\newcommand{\Fnil}[1]{{#1}_{\op{nil}}}         % F-nilpotent submodule
\newcommand{\Fred}[1]{{#1}_{\op{red}}}         % F-reduced quotient
\newcommand{\Ffred}[1]{{#1}_{\op{fred}}}       % F-full reduced subquotient
\newcommand{\deR}{\op{deR}}                  % DeRham Functor of RH correspondence

\newcommand{\ie}{\textit{i.e.}\ }           % i.e. in italics and with proper spacing afterwards
\newcommand{\eg}{\textit{e.g.}\ }           % e.g.   ....

\newcommand{\ltag}{\enm{\star}}                  % Local tags in 3 depths
\newcommand{\lltag}{\ltag\,\ltag}

\newcommand{\eqnref}[1]{(\ref{#1})}              % makes sure that equation numbers get referenced with brackets around

\theoremstyle{plain}
\newtheorem{theorem}{Theorem}[chapter]
\newtheorem*{theorem*}{Theorem}
\newtheorem{proposition}[theorem]{Proposition}
\newtheorem*{proposition*}{Proposition}
\newtheorem{corollary}[theorem]{Corollary}
\newtheorem*{corollary*}{Corollary}
\newtheorem{lemma}[theorem]{Lemma}

\newtheorem*{maintheorem1}{Main Theorem 1}
\newtheorem*{maintheorem2}{Main Theorem 2}
\newtheorem{claim}[theorem]{Claim}
\newtheorem{problem}[theorem]{Problem}

\theoremstyle{definition}
\newtheorem{definition}[theorem]{Definition}

\theoremstyle{remark}
\newtheorem{remark}[theorem]{Remark}
\newtheorem{example}[theorem]{Example}
\newtheorem{examples}[theorem]{Examples}
\newtheorem{question}[theorem]{Question}

\title{The Intersection Homology $D$--Module in Finite Characteristic}
\candidate{Manuel Blickle}
\degree{Doctor of Philosophy}
\field{Mathematics}
\chair{Associate Professor Karen E Smith}
\restofcommittee{ Professor Juha Heinonen \\ Professor Melvin Hochster \\ Professor Robert Lazarsfeld \\ Professor Tobias Stafford \\ Assistant Professor Jamie Tappenden}
\copyrighted

\begin{document}

% ----------  Frontmatter -------------------------------------------------------------------

\frontmatter   % this takes care of all front crap such as titlepage, copyrightpage...

\oneandahalfspacing
\begin{acknowledgements}
%%%%%%%%%%%%%%%%%%%%%%%%%%%%%%%%%%%%%%%%%%%%%%%%%%%%%%%%%%%%%%
%%                                                          %%
%%   This is file: aknowledgements.tex                      %%
%%   It contains the Aknowledgements of my                  %%
%%   dissertation: dissertation.tex                         %%
%%                                                          %%
%%%%%%%%%%%%%%%%%%%%%%%%%%%%%%%%%%%%%%%%%%%%%%%%%%%%%%%%%%%%%%

I thank my supervisor, Karen Smith, for her excellent guidance during the
years of my doctoral studies. Without her support, expertise and, most
importantly, encouragement this dissertation would not exist the way it
does now. Also, I thank Matthew Emerton for some helpful conversations and
for entrusting me with early manuscripts of \cite{Em.Kis,Em.Kis2}, both of
which enhanced my understanding of the subject of this dissertation
greatly.

I thank my parents for their support and trust, that, whatever it is their
son is doing (mathematics?), he must know what is best for him. I don't
know where they get this idea from, but it sure is comforting. Many other
people had a subtle impact on my life, thereby influencing my mathematical
career too. Thanks to all of you for being there for me when I needed you
most: Alexander Schwarzhaupt, Andreas Hoelzl, Barbara M\"uller, David
Voormann, Diane Miller, Elizabeth Rhoades, Elliot Lawes, Enrico Sbarra,
Eric Rupley, Friderike Tirpitz, Ingo Bertsche, J\"org Przybilla, Joel
Pitkin, Len Puch, Mike Reiskind, Nick Buchler, Peter Stahn, Rebecca
Kennedy, Sebastian Gengnagel, Stephen Keith, Thanh Tran,  Thomas
Leonhardi, Ursula Brandt, Valeri\'e Vivancos, Valentin Blickle.

I thank the Forschergruppe ``Arithmetik und Geometrie'' at the University
of Essen for their hospitality during the first half of the Year 2001.
Especially, Eckart Viehweg and H\'el\`ene Esnault who made my very
fruitful stay in Essen possible. During the month of March 2001 I stayed
at the University of Jyv\"askyl\"a. This visit was partly sponsored by
Pekka Koskela through his grant with the Academy of Finland, and by the
Deutsche Forschungsgemeinschaft.

During my graduate studies I received financial support from the
Studienstiftung des Deutschen Volkes, and from my supervisor's NSF grants
DMS 0070722 and DMS 96-25308. The dissertation is typeset with \LaTeX.

\end{acknowledgements}

%\begin{preface}
%\input{parts/preface}
%\end{preface}

\oneandahalfspacing
\tableofcontents
%\listoftables
%\listoffigures

%---------   Mainmatter  ----------------------------------------------------------------

\mainmatter   % this makes further formatting adjustments.

%%%%%%%%%%%%%%%%%%%%%%%%%%%%%%%%%%%%%%%%%%%%%%%%%%%%%%%%%%%%%%
%%                                                          %%
%%   This is file: chapter1.tex                             %%
%%   It contains the First Chapter of my                    %%
%%   dissertation: dissertation.tex                         %%
%%                                                          %%
%%%%%%%%%%%%%%%%%%%%%%%%%%%%%%%%%%%%%%%%%%%%%%%%%%%%%%%%%%%%%%

\chapter{Introduction} Let $Y \subseteq \AA^n_k$ be a
$d$--dimensional affine subvariety of affine $n$--space over the
field $k$. In the study of the singularities of $Y$, there are
traditionally (at least) two different approaches one can take.

If $y \in Y$ is a point of $Y$, one can consider the local ring of $Y$,
call it $A$, at this point. This local ring contains information about the
nature of the singularity of the variety $Y$ at $y$, such as normality,
Cohen--Macaulayness, being a complete intersection or rationality. Highly
developed tools from commutative algebra can be used to extract this
information. For example, whether the point $y$ is a rational singularity
of $Y$ can be detected by the vanishing of a certain submodule of the
local cohomology module $H^d_{m}(A)$ where $m$ is the maximal ideal of $A$
corresponding to the point $y \in Y$.

Alternatively, one can take a more geometric approach and study $Y$ via
its embedding in the smooth space $X=\AA^n_k$. Then one studies the
singularities of $Y$ through certain objects (modules, generally) on $X$.
By working completely on $X$, the smoothness of $X$ makes a host of
techniques, most notably, the theory of $D$--modules, available. As an
example, consider the cohomology module $H^{n-d}_{[Y]}(\Oo_X)$ of global
sections of $\Oo_X$ supported on $Y$. This is an $\Oo_X$--module which
contains information about the singularities of $Y$. Better still, its
elements can be ``differentiated'' and this additional structure as a
$D$--module is a powerful tool for understanding the singularities of $Y$.

Our goal is to show a certain connection between these two viewpoints,
while simultaneously showing the existence of the object featured in the
title of this dissertation. Concretely, if $k$ is of finite characteristic
$p$, then the local cohomology module $H^d_{m}(A)$ contains the \idx{tight
closure} of zero, $0^*_{H^d_{m}(A)}$, as its unique maximal proper
submodule stable under the natural action of the Frobenius. The vanishing
of this submodule is equivalent to $A$ being $F$--rational, a property
essentially equivalent to rational singularities in characteristic zero.
We want to relate this Frobenius stable $A$--submodule $0^*_{H^d_{m}(A)}$
to an equally significant $D_R$--submodule of $H^{n-d}_I(\Oo_X)$, namely
the Brylinski--Kashiwara\index{Brylinski}\index{Kashiwara}
\idx{intersection homology} $D$--module. If $k$ is of characteristic zero,
this $D$--module is characterized as the unique simple $D$--submodule of
$H^{n-d}_{[Y]}(\Oo_X)$. In Section \ref{sec.inthom} we explain the
significance of this object in the context of the \idx{Riemann--Hilbert
correspondence}.

The first difficulty we encounter is that the two objects we would like to
relate live in different characteristics. For the maximal Frobenius stable
submodule $0^*_{H^d_m(A)}$ we assumed finite characteristic. The
Brylinski--Kashiwara $D$--module, denoted $\Ll(Y,X)$, however, is only
defined in characteristic zero. Thus we are led to the natural question of
the existence of an analog of \index{$\Ll(Y,X)$}$\Ll(Y,X)$ in finite
characteristic:

\begin{question}\label{ques.LnFiniteChar}
Let\/ $k$ be a field of positive characteristic and\/ $Y$ a closed
subscheme of a smooth\/ $k$--scheme\/ $X$ of constant codimension\/ $c$.
Does\/ $\Hh^c_{[Y]}(\Oo_X)$ have a unique simple\/ $\Dd_X$--submodule?
\end{question}

To answer this question, the methods mentioned below in Section
\ref{sec.inthom} from characteristic zero are not directly applicable.
There is no obvious analog of the Riemann--Hilbert correspondence in the
finite characteristic setting, so this line of reasoning does not apply.
Although there are proofs in characteristic zero of the existence of
$\Ll(Y,X)$ that do not explicitly refer to the Riemann--Hilbert
correspondence \cite{Bry.Kash,Bjork.IntCohom}, these make substantial use
of properties of $\Dd_X$ which again fail in finite characteristic. Most
notably it is the fact that $\Dd_X$ is coherent which leads to the
important notion of a \index{holonomic}holonomic $\Dd_X$--module. Thus, in
the quest for the unique simple $\Dd_X$--submodule of $\Hh^c_{[Y]}(\Oo_X)$
in finite characteristic, a completely different approach is necessary.

In this dissertation we give a positive answer to a local version of
Question \ref{ques.LnFiniteChar}. Our main result (Main Theorem 2) on the
analog of the intersection homology $\Dd_X$--module $\Ll(Y,X)$ in finite
characteristic implies the following result.

\begin{theorem*}
Let\/ $y$ be a point in a normal algebraic variety\/ $Y$ of codimension\/
$c$ in a smooth algebraic variety\/ $X$ over perfect field of finite
characteristic. Then the local cohomology module\/ $H^c_I(R)$ has a unique
simple\/ $D_R$--submodule, where\/ $R$ is the local ring of\/ $X$ at\/ $y$
and\/ $I$ is the ideal defining\/ $Y \subseteq X$ at $y$.
\end{theorem*}

Roughly speaking, we construct this unique simple $D_R$--submodule of
$H^c_I(R)$ by applying a generalization of the Matlis duality functor that
incorporates Frobenius actions to the tight closure of zero in the local
cohomology module $H^d_m(R/I)$. This amounts to an explicit description of
this submodule, arising as a direct limit of the Frobenius powers of a
certain submodule $\tau_{\omega_A}$ (the parameter test
module\index{parameter test module}) of the \idx{canonical module}
$\omega_A$ of $A$. This explicit description yields a precise criterion
(Theorem \ref{thm.HcIisSimpleCrit}) for the simplicity of $H^c_I(R)$ as a
$D_R$--module. One of its implications is the following.
\begin{corollary*}
Let\/ $A = R/I$ be a normal domain, where\/ $R$ is a regular local\/
$F$--finite ring, and\/ $I$ is an ideal of height\/ $c$. If\/ $A$ is\/
$F$--rational, then\/ $H^c_I(R)$ is\/ $D_R$--simple.
\end{corollary*}
The connection to $F$--rationality stems from our construction of
$\Ll(A,R)$ as the dual of the tight closure of zero in the top local
cohomology module $H^d_m(A)$ of the local ring $A$. The aforementioned
parameter test module $\tau_{\omega_A}$ which gives rise to $\Ll(A,R)$ is
exactly the annihilator of $0^*_{H^d_m(A)}$, which vanishes if and only if
$A$ is $F$--rational \cite{Smith.rat}. Therefore, built into the
construction of the unique simple $D_R$--submodule $\Ll(A,R)$ is its
relation to the unique maximal $R[F]$--submodule $0^*_{H^d_m(R)}$.

One main idea of the proof of this Theorem is the following: Instead of
dealing with the $D_R$--module structure of $H^c_I(R)$ directly we
consider its structure as a unit $R[F^e]$--module; \ie we replace the
action of the differential operators by the, in some sense, complementary
(see \cite{Em.Kis2}) action of the Frobenius on $H^c_I(R)$. The
$D_R$--module structure of $H^c_I(R)$ is in fact determined by this unit
$R[F^e]$--module structure. We then go on to solve the problem in this
context of Frobenius actions by explicitly constructing the unique simple
unit $R[F^e]$--module $L(A,R)$ of $H^c_I(R)$ from $0^*_{H^d_m(A)}$. Of
course, after this was achieved, we must get back to the $D_R$--structure
to obtain the simple $D_R$--module $\Ll(A,R)$ we are looking for; at the
end we will show that in fact $\Ll(A,R) = L(A,R)$, \ie the simple unit
$R[F^e]$--submodule is $D_R$--simple.

The step of proving that $L(A,R)$ is $D_R$--simple turned out to be quite
subtle and a large part of the dissertation consists of dealing with the
interplay of unit $R[F^e]$--structures and their underlying
$D_R$--structure. This study eventually leads to the following result.
\begin{theorem*}
    Let\/ $k$ be an uncountable algebraically closed field of
    characteristic\/ $p>0$ and\/ $R$ be a regular ring, essentially of
    finite type over\/ $k$. Then a simple unit\/ $R[F^\infty]$--module is
    simple as a\/ $D_R$--module.
\end{theorem*}

We show that one can reduce the problem of the existence of a unique
simple $D_R$--submodule to the case that the coefficient field is huge. In
this case, this theorem shows that $L(A,R)$ is $D_R$--simple.

Summarizing, the key techniques used to establish the analog of the
intersection homology $\Dd_X$--module in finite characteristic are the
theory of tight closure on one side, and a theory of modules with a
Frobenius action on the other. Tight closure is used to establish the
unique maximal $R[F^e]$--submodule of $H^d_m(A)$. Then our extended Matlis
duality and constructions from the theory of $R[F^e]$--modules allow us to
construct the unique simple unit $R[F^e]$--submodule $L(A,R)$ of
$H^c_I(R)$. After that, the intimate relationship between unit
$R[F^e]$--modules and their underlying $D_R$--module structure allows one
to conclude that $L(A,R)$ is the unique $D_R$--submodule of $H^c_I(R)$.

Before we go on to discuss the contents of this dissertation in greater
detail, we include a summary of the Riemann--Hilbert correspondence and
its relation to the Brylinski--Kashiwara $\Dd_X$--module in characteristic
zero.

\section{The intersection homology $\Dd$--module $\Ll(Y,X)$}\label{sec.inthom}
This subsection is intended to explain some of the background on
\index{$\Ll(Y,X)$}$\Ll(Y,X)$ in characteristic zero and is not essential
for the understanding of this dissertation. Thereby, the reader should not
be discouraged by the appearance of many unexplained terms. What is
important to us is the phenomenology and not the details of the theory.

Intersection homology\index{intersection homology} was developed by
\idx{Goreski} and \idx{Mac\-Pherson} \cite{Gor.MacPh} as a theory of
homology groups for possibly singular varieties which satisfies many of
the nice properties that singular homology does for smooth varieties. For
example, there are versions of Poincar\'e duality, the K\"unneth formula
and the Lefschetz theorems. It has many applications to algebraic geometry
and beyond. For an excellent survey on the subject, see Kleiman
\cite{Klei.IntHom}. Its connection to $\Dd$--module theory is obtained
through the \idx{Riemann--Hilbert correspondence}, which we shall review
now.

Let $k$ be an algebraically closed field of characteristic zero (say, $k =
\CC$, the complex numbers), and $X$ a smooth algebraic variety over $k$
and let $Y$ be a closed subvariety of $X$ of codimension $c$. The
Riemann--Hilbert correspondence grew out of \index{Hilbert}Hilbert's 21st
problem \cite{Hilb.prob} on the existence of differential equations with
given singular points and monodromy. Its solution, the Riemann--Hilbert
correspondence, roughly asserts an equivalence between local systems of
$\CC$--vectorspaces and coherent $\Oo_X$--modules together with an
integrable connection. The local system is given by a representation of
the fundamental group of $X$ (via the monodromy action) and an integrable
connection for a module is nothing but an action of the ring of
differential operators.

This correspondence is vastly generalized if one requires that the objects
are preserved under pushforward by algebraic maps. This requirement forces
one to enlarge the categories since neither of them is closed under
pushforward. In this generality the Riemann--Hilbert correspondence gives
an equivalence between the derived category of bounded complexes of
sheaves of $\CC$--vectorspaces with constructible cohomology on one side
with the derived category of bounded complexes of coherent (left)
$\Dd_X$--modules with regular holonomic cohomology on the other. Here,
constructible means with respect to some stratification of $X$ and a
$\Dd_X$--module is a $\Oo_X$--module together with an action of the sheaf
of differential operators $\Dd_X$ on $X$. In this algebraic formulation
the Riemann--Hilbert correspondence was obtained by \idx{Beilinson} and
\idx{Bernstein} (see \cite{Borel.Dmod} and
\cite{Mebk.Phd,Kashiwara.FaisConst} for analytic version).

Given a $\Dd_X$--module $\Mm$ (or in general a complex of
$\Dd_X$--modules), the corresponding complex of constructible sheaves is
$\deR(\Mm) = \op{RHom}_{\Dd_X}(\Mm, \Oo_X)$ which is isomorphic (in the
derived category) to
\[
    0 \to \Mm \to \Omega^1_X \tensor \Mm \to \ldots \to \Omega^n_X \tensor
    \Mm \to 0,
\]
the \idx{deRham} complex associated to the $D_X$--module
$\Mm$.\footnote{Here we are being imprecise in the sense that one should
take the Verdier--Borel--Moore dual of $\deR(\Mm)$ to recover the
correspondence.} If $\Mm$ is a single $\Dd_X$--module, then the first map
$\Mm \to[\nabla] \Omega^1_X \tensor \Mm$ is the \idx{connection}
associated to the $\Dd_X$--module $\Mm$. If $\xn$ are local coordinates it
is given by $\nabla(m)=\sum dx_i \tensor \diff[x_i](m)$ where $d$ denotes
the universal derivation $d: \Oo_X \to \Omega^1_X$.

An important object on the side of bounded complexes with constructible
cohomology of this correspondence is the intersection homology complex
$IC_\cdot(Y)$ with respect to the middle perversity. If we consider the
pushforward of $IC_\cdot(Y)$ to $X$ under the inclusion $j: Y \into X$ we
get a simple \idx{perverse} sheaf on $X$. Under the Riemann--Hilbert
correspondence this arises as the deRham complex of a \idx{holonomic}
$\Dd_X$--module $\Ll(Y,X)$, \ie
\[
    \deR(\Ll(Y,X)) = j_*IC_\cdot(Y)[-c]
\]
where $c$ denotes the codimension of $Y$ in $X$. The fact that
$IC_\cdot(Y)$ is perverse implies that $\Ll(Y,X)$ is represented by a
\emph{single} holonomic $\Dd_X$--module. The simplicity of $IC_\cdot(Y)$
implies that $\Ll(Y,X)$ is a simple $\Dd_X$--module.

If $Z$ denotes the singular locus of $Y$, then there is a natural map
$IC_\cdot(Y) \to {}^pj_*\CC_{Y-Z}$ where ${}^pj_*$ denotes the
extraordinary direct image as in \cite{BeiBerDel.pervers}. Under the
Riemann--Hilbert correspondence the latter belongs to the local cohomology
sheaf $\Hh^c_{[Y]}(\Oo_X)$ of local sections supported on $Y$
\cite[Proposition 2]{Vil}. This establishes $\Ll(Y,X)$ as a
$\Dd_X$--submodule of $\Hh^c_{[Y]}(\Oo_X)$. In \cite{Bry.Kash}, Brylinski
and Kashiwara show that
\[
    \Ll(Y,X) \subseteq \Hh^c_{[Y]}(X,\Oo_X)
\]
is the \emph{unique simple} $\Dd_X$--submodule. This proves the existence
of the unique simple $\Dd_X$--submodule of the local cohomology sheaf
$\Hh^c_{[Y]}(\Oo_X)$ but hardly gives any concrete information about it.
The best result toward a concrete description is due to \idx{Vilonen}
\cite{Vil} who characterizes $\Ll(Y,X)$ for a complete intersection with
isolated singularities via the vanishing of local residues.

\section{Outline of Dissertation}
The main technique we use is the theory of modules with a Frobenius action
and its connection to $D_R$--module theory over a regular ring $R$. This
was first systematically developed in \cite{HaSp} of Hartshorne and
Speiser and is extensively treated in \index{Lyubeznik}Lyubeznik's
\cite{Lyub}. The preprints of \cite{Em.Kis,Em.Kis2} of Emerton and Kisin
also give a detailed study of the subject with a more arithmetic flavor; a
great part of their study of $R[F]$--modules originates with ideas from
crystalline cohomology. Many of the ideas in this dissertation are from,
or inspired by, one these three sources.

\subsection*{Chapter \ref{chap.Fmod}} We set up the necessary notation
and formalism for working with modules over rings of finite
characteristic. Most importantly, the properties of Peskine and Szpiro's
Frobenius functor are summarized. After the basic definitions are given,
we begin with a self contained, fairly detailed, introduction to the
theory of $R[F^e]$--modules, loosely following \cite{Lyub} in content, and
Emerton and Kisin \cite{Em.Kis} in notation. All the basics of the theory
are summarized, and complete arguments are given for most of the results.

\subsection*{Chapter \ref{chap.DRsub}} We focus on interaction of Frobenius
with $D_R$--module theory. The main observation is that unit
$R[F^e]$--modules carry a natural \index{DR-module@$D_R$--module}
$D_R$--module structure. This is further exploited via a version of
Frobenius descent, a powerful tool in the study of differential operators
of finite characteristic. It asserts that the Frobenius functor $F^*$
defines an equivalence of the category of $D_R$--modules with itself. The
result finds its origin in the so-called \idx{Cartier descent} which
establishes an equivalence between $R$--modules and $R$--modules equipped
with an integrable connection of $p$--curvature zero. In its most general
form, this appears as \idx{Berthelot}'s \index{Frobenius!descent}Frobenius
descent (this is where we took the name), a part of his powerful theory of
arithmetic $D_R$--modules
\cite{Ber.IntroDmodArith,Ber.OpDiff,Ber.FrobDesc}. The treatment in
Lyubeznik \cite{Lyub} of the connection between $R[F^e]$--modules and
$D_R$--modules only implicitly makes use of this theory. In our exposition
we show how the explicit use of Frobenius descent illuminates some of the
ideas of \cite{Lyub}.

After the foundational material is laid down, Chapter \ref{chap.DRsub}
contains the proof of the following result.

\newcommand{\mainthmonecontent}{
\begin{maintheorem1}
    Let\/ $k$ be an uncountable algebraically closed field of
    characteristic\/ $p>0$ and\/ $R$ be a regular ring, essentially of
    finite type over\/ $k$. Then a simple finitely generated unit\/
    $R[F^\infty]$--module is simple as a\/ $D_R$--module.
\end{maintheorem1}}
\mainthmonecontent

Examples showing that the result is sharp are given; in particular the
assumption on the algebraically closedness of $k$ can not be
weakened\footnote{We speculate that an equivalent of Main Theorem 1 can be
obtained for non-algebraically closed fields if one permits finite
algebraic extensions of $k$.}. These examples also provide a
counterexample to \cite[Remark 5.6a]{Lyub}, where he asserts a stronger
version of Main Theorem 1 without assuming that $k$ is algebraically
closed.

The somewhat awkward assumption of uncountability of the ground field in
Main Theorem 1 is due to the fact that its proof heavily uses that
$\End_{D_R}(N)=k$ for a simple $D_R$--module when $k$ is algebraically
closed. Unfortunately, this is only known in the case that $k$ has
strictly bigger cardinality than a $k$--basis of $R$, due to an argument
of Dixmier \cite{Dix}. We expect that this might be true without the
assumption that $k$ be uncountable. In characteristic zero, this is known
due to \idx{Quillen}'s lemma \cite{Quill}, but again the characteristic
zero proof does not transfer to finite characteristic since it crucially
uses the noetherianess of the ring of differential operators in
characteristic zero, a property not shared by the differential operators
in finite characteristic.

The proof of Main Theorem 1 can roughly be summarized as follows: From
\cite[Theorem 5.6]{Lyub} one can easily deduce that a finitely generated
simple unit $R[F^e]$--module $M$ must be semisimple as a $D_R$--module and
of finite length. Using this we can identify $D_R$--submodules of $M$ with
subspaces of an appropriate $k$--vectorspace $V$, also equipped with a
Frobenius action. The point is then to observe that the $F$--stable
$D_R$--submodules of $M$ correspond to the $F$--stable subspaces of $k$.
This allows us to reduce the problem to the case of a $k$--vectorspace,
where it can be worked out explicitly.

\subsection*{Chapter \ref{chap.Functors}} We collect some general
observations surrounding the following question: If $K: R\text{--mod} \to
A\text{--mod}$ is a functor, can we naturally extend $K$ to a functor from
$R[F^e]$--modules to $A[F^e]$--modules? If $K$ is covariant, it is easily
seen that this only depends upon the existence of a natural transformation
of functors $\F[e] \circ K \to R \circ \F[e]$. For a contravariant $K$,
though, it is \emph{a priori} not clear how this can be achieved. The
concept of a generator of an $R[F^e]$-module of Lyubeznik enables one to
find a natural extension of $K$ which respects Frobenius actions. Given a
natural transformation $\F[e] \circ K \to K \circ \F[e]$ of functors one
can extend $K$ to a functor $\Kk: R[F^e]\text{--mod} \to
A[F^e]\text{--mod}$. After some basic properties of $\Kk$ are derived we
show that under appropriate assumptions an equivalent definition using
Hartshorne and Speisers's \emph{leveling functor} $G$ can be given.

The purpose for developing this abstract machinery is to apply it to the
Matlis duality functor of a complete local ring $R$. It follows that the
Matlis duality functor $D = \Hom(\usc,E_R)$ has a natural extension $\Dd$
to a functor on $R[F^e]$--modules (here $E_R$ denotes the injective hull
of the residue field of $R$). For $R[F^e]$--modules which are cofinite as
$R$--modules and supported on $\Spec A= \Spec R/I$ this extension of
Matlis duality $\Dd$ specializes to the functor $\Hh_{A,R}$ introduced in
\cite{Lyub}. This places $\Hh_{A,R}$ in the natural, more general,
framework we introduced.

The extended Matlis dual $\Dd$ is our main tool to translate the result
that $H^d_m(A)$ has a unique maximal $F$--stable submodule into the
statement that $H^c_I(R)$ has a unique simple unit $R[F^e]$--submodule,
\ie it is an integral part of the proof of the existence and construction
of the unique simple unit $D_R$--submodule of $H^c_I(M)$. The following
properties of $\Dd$ are most important in this context (Proposition
\ref{prop.DdGraded}).
\begin{itemize}
\item $\Dd$ is an exact functor.
\item $\Dd(M)=0$ if and only if $M$ is $F$--nilpotent.
\item If $M$ is a simple $R[F^e]$--module, then $\Dd(M)$ is a simple
      unit $R[F^e]$--module.
\end{itemize}
In particular, the last fact, which is implicit in Lyubeznik's treatment
of $\Hh_{A,R}$, shows that $\Dd(H^d_m(A)/0^*_{H^d_m(A)})$ must be a simple
unit $R[F^e]$--module. Indeed, by maximality of the $R[F^e]$--submodule
$0^*_{H^d_m(A)}$ the quotient $H^d_m(A)/0^*_{H^d_m(A)}$ is a simple
$R[F^e]$--module.

The critical example is then to observe that the extended Matlis dual
$\Dd$ applied to the top local cohomology module $H^d_m(A)$ gives the
local cohomology module $H^c_I(R)$, \ie $\Dd(H^d_m(A))=H^c_I(R)$. Thus
$\Dd(H^d_m(A)/0^*_{H^d_m(A)})$ is the simple unit $R[F^e]$--submodule of
$\Dd(H^d_m(A))=H^c_I(R)$.

\subsection*{Chapter \ref{chap.IntHomD_R--mod}} We give a brief
introduction to tight closure theory and summarize the constructions and
results relevant to us. Most important is the already mentioned result
that (under the hypothesis that the completion of $A$ is a domain) the top
local cohomology module of $A$ has a unique maximal nontrivial
$A[F^e]$--submodule; this submodule is the tight closure of zero
$0^*_{H^d_m(A)}$ and its annihilator in the canonical module $\omega_A$
under the Matlis duality pairing $H^d_m(A) \times \omega_A \to E_A$ is
called the parameter test module $\tau_{\omega_A}$. Ultimately, we show
that the parameter test module is the root, in the sense of \cite{Lyub},
of the unique simple $R[F^e]$--submodule $L(A,R)$; \ie $L(A,R)$ arises as
the increasing union (interpreted suitably) of the Frobenius powers
$\F[e]\tau_{\omega_R}$ of the parameter test module. Combining these
results from tight closure theory with all other results discussed so far
we are able to prove the following result.
\begin{theorem*}
    Let\/ $(R,m)$ be regular, local and\/ $F$--finite. Let\/ $I$ be an
    ideal with\/ $\height I = c$ such that\/ $A=R/I$ is analytically
    irreducible. Then the local cohomology module\/ $H^c_I(R)$ has a
    unique simple\/ $D_R$--submodule\/ $\Ll(A,R)$.
\end{theorem*}

Due to the necessity of two reductions, the first proof of this theorem
gives only the existence and not a concrete description of $\Ll(A,R)$. One
reduction involves completing in order to apply our extension of the
Matlis dual Functor $\Dd$, which is only defined for complete rings. The
second reduction involves extending the coefficient field of $R$ so it is
big enough (uncountable) and algebraically closed, so that we are able to
apply Main Theorem 1.

In a more careful analysis following the proof of Main Theorem 2 we are
able to show that $\Ll(A,R)$ is in fact well behaved under completion as
well as extension of the ground field. This leads to the following
concrete  description of $\Ll(A,R)$ (Theorem \ref{thm.Llconcrete}).
\begin{theorem*}
    Let\/ $A = R/I$ be as in the last theorem. Then the unique simple\/
    $D_R$--submodule\/ $\Ll(A,R)$ is the unique simple\/
    $R[F^\infty]$--module of\/ $H^c_I(R)$. It arises naturally as the
    direct limit
    \[
        \tau_{\omega_A} \to R^1 \tensor \tau_{\omega_A} \to R^2 \tensor
        \tau_{\omega_A} \to \ldots
    \]
    where the map is the restriction of the natural map\/ $\omega_A \to
    R^1 \tensor \omega_A$ dual to the\/ $R[F^e]$--module structure map\/
    $R^1 \tensor H^d_m(R/I) \to[\theta] H^d_m(R/I)$ via local duality.
\end{theorem*}

It is this concrete description of $\Ll(A,R)$ which leads to a precise
criterion for the $D_R$--simplicity of $H^c_I(R)$ (Theorem
\ref{thm.HcIisSimpleCrit}).
\begin{theorem*}
    Let\/ $R$ be regular, local and\/ $F$--finite. Let\/ $I$ be an ideal
    such that $A=R/I$ is analytically irreducible. Then\/ $H^c_I(R)$ is\/
    $D_R$--simple if and only if the tight closure of zero in\/ $H^d_m(A)$
    is\/ $F$--nilpotent.
\end{theorem*}

In particular, if $A$ is $F$--rational (has only rational singularities if
obtained by reduction from characteristic zero), then $H^c_I(R)$ is simple
as a $D_R$--module. More precisely, if $A$ is analytically irreducible and
$F$--injective, then $A$ is $F$--rational if and only if $H^c_I(R)$ is
$D_R$--simple. Interestingly, an outcome of the results that $\Ll(R,A)$ is
well behaved under completion is that the parameter test module also is
well behaved under completion (Corollary \ref{cor.TestCommCompl}).
\begin{corollary*}
    Let\/ $A = R/I$ be a normal domain,where\/ $R$ is regular, local and\/
    $F$--finite. Then the parameter test module commutes with completion.
    \ie $\tau_{\omega_{\wh{A}}} = \wh{A} \tensor_A \tau_{\omega_A}$.
\end{corollary*}
This was one property of $\tau_{\omega_A}$ left unsettled in
\cite{Smith.test} where it is proven that $\tau_{\omega_A}$ commutes with
localization for complete Cohen--Macaulay domains $A$.

We finish Chapter \ref{chap.IntHomD_R--mod} with some examples of graded
complete intersections where the construction of $\Ll(A,R)$ can be made
much more explicit. In the examples we give, $\Ll(A,R)$ arises as the
intersection of kernels of natural maps $H^c_I(R) \to H^n_m(R)=E_R$. This
hints at a similar description of $\Ll(A,R)$ due to \idx{Vilonen}
\cite{Vil} in the characteristic zero setting.

As a last application, we show a characteristic $p$ analog of a result of
S.P.Smith showing that $H^1_{(f)}(k[x,y])$ is a simple
$D_{k[x,y]}$--module if the plane curve defined by $f$ is a cusp.

\subsection*{Chapter \ref{chap.Problems}} Since this is only the
beginning in a study of the finite characteristic analog of the
intersection cohomology $D_R$--module there are many questions left
unanswered. Some of them are discussed in Chapter \ref{chap.Problems}.
First, we must ask whether this approach globalizes. We showed that if $Y
\subseteq X$ is a closed irreducible subvariety of the smooth $k$--variety
$X$, then at each point $x \in X$ the stalk $H^c_{I_x}(\Oo_{x,X})$ of
$\Hh^c_{[Y]}(\Oo_X)$ at this point has a unique simple
$D_{\Oo_{x,X}}$--submodule $\Ll(\Oo_{x,Y},\Oo_{x,X})$. Is there a (unique
simple) $\Dd_X$--submodule of $H^c_{I_x}(\Oo_{x,X})$ which localizes at
each point to the unique simple $D_{\Oo_{x,X}}$--submodules of the stalks?
We describe a possible line of attack by reducing the question to a
question about the behavior of a global version of the parameter test
module $\tau_{\omega_R}$ under localization.

The next question is whether the characteristic zero module reduces modulo
$p$ to the finite characteristic version we produce here. For this one
first has to develop a framework for reducing $D_R$--modules to finite
characteristic. If this is achieved the answer of this question should be
in reach.

One can further ask if there is a direct proof of the existence of
$\Ll(R,A)$, \ie one without the detour through $R[F^e]$--modules. An
answer to this question might be found by investigating the extensive work
on $D_R$--modules in finite characteristic due to Berthelot. Even though
his goal of establishing the theory is more arithmetic it is likely that
it offers great insights about questions like the one we are considering.
It might even be that this approach will give a proof somewhat resembling
the one in characteristic zero. Furthermore, Berthelot's setup is global,
so a study of his theory is likely to improve the understanding of the
preceding two questions, too.

In \cite{Em.Kis,Em.Kis2}, Emerton and Kisin construct an (anti)
equivalence between the derived category of bounded complexes of unit
$R[F^e]$--modules with finitely generated cohomology and the derived
category of bounded complexes of \'etale $\FF_{P^e}$--sheaves with
constructible cohomology. This can be viewed as a finite characteristic
analog of the Riemann--Hilbert correspondence. Since in characteristic
zero the intersection cohomology complex, which corresponds to $\Ll(Y,X)$,
is such an important object we hope that the corresponding object of
$\Ll(Y,X)$ in finite characteristic under Emerton and Kisin's
correspondence is equally significant.

%%%%%%%%%%%%%%%%%%%%%%%%%%%%%%%%%%%%%%%%%%%%%%%%%%%%%%%%%%%%%%
%%                                                          %%
%%   This is file: chapter2.tex                             %%
%%   It contains the Second Chapter of my                   %%
%%   dissertation: dissertation.tex                         %%
%%                                                          %%
%%%%%%%%%%%%%%%%%%%%%%%%%%%%%%%%%%%%%%%%%%%%%%%%%%%%%%%%%%%%%%

\chapter{Modules with Frobenius action}\label{chap.Fmod}

In this chapter we give, after the notation is set up, a fairly self
contained introduction to the theory of $R[F^e]$--module, \ie $R$--modules
on which a power of the Frobenius morphism acts. Except for notation we
follow Lyubeznik \cite[Section 1-3]{Lyub} fairly closely. Even though all
the key results can be found there we give many of the proofs. We hope
that our somewhat different viewpoint, and some alternative arguments will
provide the reader with some new insights.

There are two differences in our treatment which must be mentioned: First,
we generalize to powers of the Frobenius. This is straightforward (cf.\
\cite[Remark 5.6]{Lyub} and \cite{Em.Kis}). Secondly, we do not assume
that the ring $R$ is regular but instead indicate when this is really a
necessary assumption. Those are both fairly minor adjustments.

\section{Notation and Generalities}\label{sec.Notation}
Unless otherwise specified, all objects will be of finite characteristic
$p$ whenever this notion applies. All rings and \idx{ring} homomorphisms
are unitary. Modules over a ring are left modules and morphisms are of
left modules unless we explicitly say otherwise. If $A$ is a ring, we
denote by $A$--\mod the category of left modules over $A$ and by \mod--$A$
the category of right $A$--modules.

The symbol $R$ will always denote a commutative noetherian ring over a
field of finite characteristic. We denote by $F=F_R$\index{$F$} the
Frobenius \index{Frobenius!map} map on $R$ which raises each element of
$R$ to its $p$TtH power. Thinking of $F$ as a map $F: \Spec R \to \Spec
R$, it is given as the identity on the underlying topological space, and
the $p$TtH power map on the structure sheaf. We denote by
$R^e$\index{$R^e$} the $R$--$R$--bimodule which, as a left $R$-module, is
just $R$, with the structure on the right given by the $e$TtH iterate of
the Frobenius map: for $m \in R^e$ and $r \in R$, $m \cdot r=r^{p^e}m$. In
order to cause no confusion with the standard notation of a direct sum we
will denote the $n$-fold direct sum of $R$ always by $R^{\dirsum
n}$\index{$R^{\dirsum n}$}. To reintroduce confusion we denote by
$R^{p^e}$\index{$R^{p^e}$} the subring of $R$ consisting of the $p^e$TtH
powers of elements of $R$. If $R$ is a domain, the inclusion $R^{p^e}
\subseteq R$ is isomorphic to the inclusion of $R \subseteq R^{1/p^e}$,
where $R^{1/p^e}$\index{$R^{1/p^e}$} is the overring of $p^e$TtH roots of
$R$.

The most important tool is \idx{Peskine} and \idx{Szpiro}'s Frobenius
\index{Frobenius!functor} functor \cite{Pes.Szp}:

\begin{definition}
    The \emph{Frobenius functor} is the right exact functor from
    $R$--modules to $R$--modules given by
    \[
        \F_RM \defeq R^1 \tensor_R M.
    \]
    Its $e$TtH iterate is $\F[e]_RM=R^e\tensor_R M$. If the danger of confusion
    is low, we omit the subscript and write just $\F$\index{$\F$} for $\F_R$\index{$\F_R$}.
\end{definition}

Note that the notation chosen is suggestive: the Frobenius functor is
indeed the pullback functor for the Frobenius map $F: \Spec R \to \Spec
R$. The right exactness of $\F[e]$ is clear by the right exactness of
tensor. Other properties $\F[e]$ inherits from the fact that it is nothing
but tensoring with $R^e$ are:
\begin{itemize}
    \item $\F[e]$ commutes with direct sums.
    \item $\F[e]$ commutes with direct limits.
    \item \label{x.FcommutesLoc} $\F[e]$ commutes with localization. Indeed, for any multiplicatively
        closed subset $S \subseteq R$ one has an isomorphism of $R$--bimodules
        $S^{-1}R \tensor_R R^e \cong R^e \tensor_R S^{-1}R$ via the map
        $\frac{r}{s} \tensor r'\mapsto rr's^{p^e-1} \tensor \frac{1}{s}$.
\end{itemize}
Throughout this dissertation we use the notation $\F[e]$\index{$\F[e]$}
and $R^e \tensor_R \usc$\index{$R^e \tensor_R \usc$} interchangeably.
Sometimes the brevity of $\F[e]$ appeals, at other times emphasizing the
structure as a tensor product is of advantage.

The crucial property which makes the Frobenius functor such a powerful
tool when $R$ is regular is the fact that it is also left exact in this
case. By a theorem of \idx{Kunz} \cite{Kunz}, this even characterizes
regularity in finite characteristic, \ie $R$ is regular if and only if the
\index{Frobenius!functor!flat} Frobenius functor is flat (\ie left exact).
Thus, under the additional assumption that $R$ is \emph{regular}, we have:
\begin{itemize}\label{x.Fregcommutes}
    \item $\F[e]$ is exact.
    \item $\F[e]$ commutes with arbitrary (not necessarily direct)
    sums.\note{Brian Conrad remarks that this might require
    $F$--finiteness, not sure though}
    \item $\F[e]$ commutes with finite intersections.
\end{itemize}
For details, the reader should refer to Lyubeznik \cite[Remarks
1.0]{Lyub}. To get a feeling for the Frobenius functor we consider a
simple but instructive example.
\begin{example}\label{ex.FonIdeal}
    Let $I$ be an ideal in a regular ring $R$. Then $\F[e]R=R^e \tensor R$
    is canonically isomorphic to $R$ itself. Under this isomorphism
    $\F[e]I$ corresponds to the ideal $I^{[p^e]}$ of $R$, where
    $I^{[p^e]}$\index{$I^{[p^e]}$}
    denotes the ideal generated by the $p^e$TtH powers of the elements of
    $I$. The quotient $A=R/I$ can be viewed as a module over itself as
    well as a module over $R$. Note that $\F[e]_A A \cong A$ and $\F[e]_R
    A \cong R/I^{[p^e]}$ are different. Worse, most times $\F[e]_R A$ is
    not even an $A$--module in any obvious way. This points out that we
    have to pay attention to the ring we are working over when using the
    Frobenius functor.
\end{example}

\section{Definition and basic properties of $R[F]$--modules} Here, the
basic properties of our fundamental objects are reviewed. These are the
$R[F]$--modules, \ie $R$--modules $M$ endowed with an action of the
Frobenius.

As it turns out, many of the objects and properties defined in \cite{Lyub}
are well defined for non--regular $R$ too. The real power of the theory
though derives from the flatness of the Frobenius in the regular case.
Thus, when we get deeper into the theory we will almost exclusively work
over a regular ring $R$. But for now only assume $R$ to be commutative and
noetherian.

\begin{definition}\label{def.RFmod}
    An \emph{$R[F^e]$--module} \index{RFmodule@$R[F^e]$--module}\index{$\theta^e$} is an $R$-module $M$
    together with an $R$-linear map
    \[
        \theta^e_M: \F[e]_R M \to{} M.
    \]
    A morphism between two $R[F^e]$--modules $(M,\theta^e_M)$ and
    $(N,\theta^e_N)$ is an $R$--linear map $\phi: M \to N$ such that the
    following diagram commutes:
\begin{equation} \label{eqn.RFhom}
\begin{split}
\xymatrix@C=3pc{ {\F[e]M} \ar[r]^{\F[e](\phi)} \ar[d]_{\theta^e_M} & {\F[e]N} \ar[d]^{\theta^e_N} \\
                {M} \ar[r]^{\phi}                            & {N} }
\end{split}
\end{equation}
    The category of all $R[F^e]$--modules together with these morphisms we
    denote by \emph{$R[F^e]$--\mod.}\index{RF-mod@$R[F^e]$--\mod}
\end{definition}
To give an alternative description of the category $R[F^e]$--\mod we apply
the adjointness of pushforward and pullback for a map of schemes (cf.\
\cite[page 110]{Hartshorne}) to the Frobenius on $\Spec R$:
\[
    \Hom_R(\F[e]M,M) \cong \Hom_R(M,F^e_*M)
\]
This shows that the $\theta^e$'s in the definition are in one-to-one
correspondence with maps $F^e_M: M \to F^e_*M$. If $\theta^e_M$ and
$F^e_M$ are identified under this correspondence we get a commutative
diagram:
\begin{equation}\label{eqn.FandTheta}
\begin{split}
\xymatrix{ {} & {R^e \tensor M} \ar[d]^{\theta_M^e} \\
           {M} \ar[ru]^{F^e_R\tensor \id_M} \ar[r]_{F_M^e} & {M}  }
\end{split}
\end{equation}
Therefore, $F^e_M=\theta^e_M \circ (F^e_R \tensor \id_M)$ and
$\theta^e_M(r \tensor m)=rF^e(m)$. Thus we can equivalently think of an
$R[F^e]$--module as an $R$--module together with such an $R$--linear map
$F^e_M: M \to F^e_*M$ (we can think of $F^e_M$ is as a $p^e$--linear map
from $M \to M$; as such it is not $R$--linear but we have
$F^e(rm)=r^{p^e}F^e(m)$ as one would expect from a Frobenius action).
Consequently, in order to specify an $R[F^e]$--module we either indicate
the Frobenius structure\index{Frobenius!structure} $\theta^e_M: \F[e]M \to
M$ or a Frobenius action\index{Frobenius!action} $F^e_M:M \to F^e_*M$. To
simplify notation we often drop the subscripts and describe an
$R[F^e]$--module as a tuple $(M,F^e)$, $(M,\theta^e)$ or, somewhat
redundantly, as a triple $(M,\theta^e,F^e)$.

If we define the ring $R[F^e]$ by formally adjoining the non-commutative
variable $F^e$ to $R$ and forcing the relations $r^{p^e}F^e=F^er$ for all
$r \in R$, then an $R[F^e]$--module as defined above is in fact just a
module over this (non-com\-muta\-ti\-ve) ring $R[F^e]$\index{$R[F^e]$} ;
more precisely:
\begin{proposition}
    The category\/ $R[F^e]$--\mod is equivalent to the category of left modules
    over the ring\/ $R[F^e]$.
\end{proposition}
\begin{proof}
Observe that an action of $R[F^e]$ on an $R$--module $M$ is nothing but a
$p^e$-linear map $F^e:M \to M$. Thus the modules over $R[F^e]$ are pairs
$(M,F^e_M)$ where $M$ is an $R$--module and $F^e_M: M \to F^e_{R*}M$ is
$R$--linear. By the adjointness discussed above these pairs $(M, F^e_M)$
are in one-to-one correspondence with $R[F^e]$--modules (\ie the pairs
$(M, \theta^e)$).

It remains to show that under this correspondence an $R$--linear $\phi: M
\to N$ is linear over the ring $R[F^e]$ if and only if $\phi$ is a map of
$R[F^e]$--modules. $R[F^e]$--linearity comes down to the commutation of
the following diagram:
\begin{equation}
\begin{split}
\xymatrix{ {M} \ar[r]^{\phi} \ar[d]_{F^e_M} &{N} \ar[d]^{F^e_N} \\
           {M} \ar[r]^{\phi}                &{N}   }
\end{split}
\end{equation}
With the help of Diagram \eqnref{eqn.FandTheta} this is easily seen to be
equivalent to the commutativity of the Diagram \eqnref{eqn.RFhom} in
Definition \ref{def.RFmod}.
\end{proof}

\begin{remark}[cf.\ \protect{\cite{Em.Kis}}]\label{rem.RFsubEnd}\note{this might require
the additional assumption that all maximal ideals have positive depth as
in \cite{Em.Kis}} There is a natural map of rings $R[F^e] \to
\End_{\FF_{p^e}}(R)$ which sends $F$ to the Frobenius map on $R$. This is
an inclusion whenever $R$ contains a non zero divisor that is not a unit
in $R$: Suppose $\phi = \sum a_iF^{ei}_R$ is the zero endomorphism on $R$.
If $r \in R$ is a non invertible non zero divisor consider the equations
$0=\phi(r^t)$ for all $t \geq 0$. Thus, after dividing this equation by
$r^t$ we get
\[
-a_0 = a_1r^{t(p^e-1)}+a_2r^{t(p^{2e}-1)}+\ldots+a_dr^{t(p^{de}-1)}
\]
which means that $a_0$ is divisible by arbitrary high powers of $r$ (since
the right hand side is so for increasing $t$). Thus $a_0 = 0$. Inductively
one concludes that all $a_i$ are in fact zero. Thus we can think of
$R[F^e]$ as the subring of $\End_{\FF_{p^e}}(R)$ generated by $R$ and the
$e$TtH iterate of the Frobenius map on $R$ in this case.

As an example of the case that the map $R[F^e] \to \End_{\FF_{p^e}}(R)$ is
not injective assume $R = \FF_{p^e}$. Then $-1+F^e$ is the zero
endomorphism on $R$ thus $R[F^e]$ is not a subring of
$\End_{\FF_{p^e}}(R)$.
\end{remark}

Since the category $R[F^e]$--\mod is just the category of left modules
over the ring $R[F^e]$ (justifying the notation $R[F^e]$--\mod) the
natural inclusion of rings $R[F^{ne}] \subseteq R[F^{e}]$ represents the
category of $R[F^{e}]$--modules as a subcategory of $R[F^{ne}]$--modules.
Note that it is not a full subcategory since a submodule which is stable
under the action of $F^{ne}$ might not be stable under the action of $F^e$
(cf.\ Section \ref{sec.ExAndCounterEx}, Example \ref{ex.DRnorRFinfty}).
All the categories of
$R[F^e]$--modules\index{RF-mod@$R[F^e]$--mod!abelian} for various $e$ are
abelian categories (as module categories over associative rings are
abelian \cite[Section 1.1]{Weib.hom}). They form a directed system and we
call the direct limit of these categories the category of
$R[F^\infty]$--modules. The objects of
$R[F^{\infty}]$--\mod\index{RF-mod@$R[F^{\infty}]$--\mod} are just the
directed union of the objects of $R[F^e]$ for all $e$, \ie an
$R[F^{\infty}]$--module is just an $R[F^e]$--module for some, not
explicitly specified $e$. It is straightforward to check that the category
$R[F^{\infty}]$ is also abelian.

\begin{definition}
    An $R[F^e]$--module $(M, \theta^e)$ is called a \emph{unit}
    $R[F^e]$--module\index{unit $R[F^e]$--module} if
    \[
        \theta^e: \F[e]M \to M
    \]
    is an isomorphism.
\end{definition}

We denote the category of unit $R[F^e]$--modules by
$uR[F^e]$--\mod\index{uRF-mod@$uR[F^e]$--\mod}. Since the unit property is
preserved by the inclusions $R[F^e]\text{--\mod} \subseteq
R[F^{ne}]\text{--\mod}$ we can form the directed union of all the
categories $uR[F^e]$--\mod and obtain the category of unit
$R[F^\infty]$--modules, denoted by
$uR[F^\infty]$--\mod\index{uRF-mod@$uR[F^\infty]$--\mod}, of course.

In the case that $R$ is regular, $uR[F^e]$--\mod is an
abelian\index{uRF-mod@$uR[F^e]$--mod!abelian} subcategory
\label{x.uRFisAbelian} of $R[F^e]$--\mod, which ensures that the kernel of
a map of unit $R[F^e]$--modules is also unit. This follows from exactness
of $\F[e]_R$ for $R$ regular. This also implies that $uR[F^\infty]$--\mod
is an abelian subcategory of $R[F^\infty]$--\mod. These statements are
false if $R$ is not regular as then kernels of maps of unit
$R[F^e]$--modules are generally not unit.

The following notational conventions have proved very practical and will
be used through\-out: Let $N$ be an $R$ submodule of an $R$--module $M$.
\begin{enumerate}
    \item $N^{[p^e]}$\index{$N^{[p^e]}$} denotes the image of $R^e \tensor N$ in $R^e \tensor M$.
    \item If $(M,\theta^e,F^e)$ is an $R[F^e]$--module, then
    the set $F^e(N)$\index{$F^e(N)$} is the subset of $M$ consisting of the elements
    $\set{F^e(n) |\ n \in N}$. By abuse of notation, $F^e(N)$ will
    primarily denote the $R$--module generated by this set.
\end{enumerate}
Note that with this notation $F^e(N)=\theta^e(N^{[p^e]})$. On the few
occasions where this abusive notation may lead to confusion we will
explicitly say ``the set $F^e(N)$'' if we do \emph{not} refer to the
$R$--module generated by this set. One easily verifies that $N$ is an
$R[F^e]$--submodule of $M$ if and only if $F^e(N)\subseteq N$.
Furthermore, if $M$ is unit and $R$ is regular then $N$ is a unit
$R[F^e]$--submodule of $N$ if and only if $F^e(N)=N$ and $\theta^e$ is in
fact an isomorphism from $\F[e]N$ to $F^e(N)$.

\begin{examples}\label{ex.BasicRFmods}
As seen in Example \ref{ex.FonIdeal} $\F[e]R$ is canonically isomorphic to
$R$ and thus $R$ is a unit $R[F^e]$--module. An ideal $I \subseteq R$ is
an $R[F^e]$--submodule of $R$ since $F(I)=I^{[p^e]}\subseteq I$. In
general $I$ is however not a unit $R[F^e]$--submodule since the inclusion
$I^{[p^e]} \subseteq I$ is normally strict.

Let $S \subseteq R$ be a multiplicatively closed subset of $R$. Then the
localization $S^{-1}R$ is naturally a unit $R[F^e]$--module. Its
structural map $\theta^e_{S^{-1}R}: R^e \tensor S^{-1}R \to S^{-1}R$ is
given by sending $r'\tensor\frac{r}{s}$ to $\frac{rr^{p^e}}{s^{p^e}}$. Its
inverse is the map $\frac{r}{s} \mapsto s^{p^e-1}r \tensor \frac{1}{s}$.
Furthermore the natural localization map $R \to S^{-1}R$ is a map of unit
$R[F^e]$--modules.

The local cohomology modules $H^i_I(R)$\index{$H^i_I(R)$} of $R$ with
support in $I$ can be calculated as the cohomology modules of the \Czech
complex\index{Czech complex@\Czech complex}
\[
    \Cz{R;\xn} = R \to \dirsum R_{x_i} \to \dirsum R_{x_ix_j} \to \ldots \to
    R_{x}
\]
where $\xn$ are a set of generators of $I$ and $x$ denotes the product of
the $x_i$'s. The modules of $\Cz{R;\xn}$ are localizations of $R$ and
therefore unit $R[F^e]$--modules. The maps of the \Czech complex are just
signed sums of localization maps and therefore maps of $R[F^e]$--modules.
Thus the \idx{local cohomology} modules $H^i_I(R)=\Cz{R;\xn}$ are
$R[F^e]$--modules as the category of $R[F^e]$--modules is abelian. If $R$
is regular then the local cohomology modules $H^i_I(R)$ are unit
$R[F^e]$--modules for the same reason ($uR[F^e]$--\mod is abelian for
$R$--regular).

Also note that if $(R,m)$ is local of dimension $n$, then the top local
cohomology module $H^n_m(R)$ is a unit $R[F^e]$--module even if $R$ is not
regular. This follows since $H^n_m(R)$ is the cokernel of the last map of
an appropriate \Czech complex arising from a system of parameters of $R$.
As $\F[e]_R$ is always right exact this is enough to conclude that the
cokernel is in fact unit.
\end{examples}

\begin{remark}
We have to comment on the choice of notation. The reader familiar with
Lyubeznik's\index{Lyubeznik} paper \cite{Lyub} will notice that the
notation chosen here is very different from Lyubeznik's. For easy
reference a dictionary between the notation in \cite{Lyub} and the
notation used here is given in Table \ref{tbl.Notation}.

\begin{table}[h]
\begin{center}
\begin{tabular}{|c|c|}\hline
  Notation in \cite{Lyub}\rule[-2mm]{0mm}{7mm} & Notation here \\ \hline
  $R\{f\}$--module\rule{0mm}{5mm} &               $R[F]$--module \\
  $F$--module      &   unit $R[F]$--module\\
  \hspace{1cm}$F$--finite module\hspace{1cm} &  \hspace{1cm}finitely generated unit $R[F]$--module\hspace{1cm} \\
  $F(M)$               &  $\F M$ \\
  $\theta^{-1}$            &  $\theta$ \\
  $\theta^{-1}(F(N))$ & $F(N)$ \\
  $\Hh_{R,A}$\rule[-2mm]{0mm}{5mm} & \Dd \\ \hline
\end{tabular}
  \caption{Comparing notation with Lyubeznik's
\cite{Lyub}.}\label{tbl.Notation}
\end{center}
\end{table}
This notation is essentially an adaptation of the notation used by
\idx{Emerton} and \idx{Kisin} \cite{Em.Kis,Em.Kis2} which in turn
originates from \idx{crystalline cohomology}. Another place where modules
with a Frobenius action appear and are systematically studied is in
\idx{Hartshorne} and \idx{Speiser} \cite{HaSp}. They call a unit
$R[F^e]$--module an $(R,F)$--module with \idx{level structure}.
\end{remark}

\subsection{Tensor product for
$R[F^e]$--modules}\label{sec.PropRFmod}\index{RF-mod@$R[F^e]$--mod!tensor
product} Let $R \to A$ be a map of rings, the functor $A \tensor_R
\usc$ can be extended to a functor from $R[F^e]$--\mod to
$A[F^e]$--\mod preserving unit modules. By this we mean that given
an $R[F^e]$--module $(M,\theta^e$), the $A$-module $A \tensor_R M$
carries a natural $A[F^e]$--module structure. To define the
Frobenius structure ${\theta'}^e$ on $A \tensor_R M$, consider the
following isomorphism of $A$--$R$ bimodules
\[
  \pi: \ A^e \tensor_A A \to[\cong] A^e \to[\cong] A \tensor_R R^e.
\]
which is the composition of the two indicated natural isomorphisms. Then
we define
\[
{\theta'}^e: A^e \tensor_A A \tensor_R M \to[\pi \tensor_R \id_M]
A \tensor_R R^e \tensor_R M  \to[\id_A \tensor \theta^e] A
\tensor_R M.
\]
Functoriality of the definition of the $A[F^e]$--structure on $A
\tensor_R M$ is clear by naturality of the isomorphism $\pi$.
Furthermore, ${\theta'}^e$ is an isomorphism if $\theta^e$ is an
isomorphism since $\pi$ is always an isomorphism. If $A$ is
faithfully flat then ${\theta'}^e$ is an isomorphism if and only
if $\theta^e$ is an isomorphism.

Since, by restriction, every $A[F^e]$--module is also an $R[F^e]$--module,
the tensor product $A \tensor_R M$ is also an $R[F^e]$--module. But it is
important to keep in mind that even if $M$ is a unit $R[F^e]$--module, $A
\tensor_R M$ is normally \emph{not} a unit $R[F^e]$--module (as we just
argued it is a unit $A[F^e]$--module in this case). For example take $A =
R/I$ and $M=R$. Then $M$ is a unit $R[F^e]$--module but $A \tensor_R
M=R/I=A$ is certainly not a unit $R[F^e]$--module as $R^e \tensor M \cong
R/I^{[p^e]}$, and this is almost never isomorphic to $R/I$.

There are exceptions to this behavior; the most important for us are if
$A$ is a localization of $R$, if $A$ is the completion of the local ring
$R$ along its maximal ideal, or if $A$ arises from $R$ by extending a
perfect field contained in $R$. In all these three cases, if $M$ is a unit
$R[F^e]$--module, then $A \tensor_R M$ is also a unit $R[F^e]$--module. In
fact, every unit $A[F^e]$--module $\wt{M}$ is a unit $R[F^e]$--module
regardless if it arises from a unit $R[F^e]$--module via base change. The
remainder of this section is a discussion of this phenomenon.

\subsubsection{Localization}\index{RF-mod@$R[F^e]$--mod!localization} The special case of $A=S^{-1}R$, a localization of $R$ at some
multiplicatively closed set $S \subseteq R$ deserves extra attention (cf.\
Example \ref{ex.BasicRFmods}). If $(M,\theta^e)$ is a (unit)
$R[F^e]$--module, then $S^{-1}M = S^{-1}R \tensor_R M$ is a (unit)
$S^{-1}R[F^e]$--module as just discussed. Since $\F[e]_R$ commutes with
localization (cf.\ page \pageref{x.FcommutesLoc}) we see that, as
$R$--modules, $\F[e]_R(S^{-1}M) \cong \F[e]_{S^{-1}R}(S^{-1}M)$ and the
natural map $M \to S^{-1}M$ is a map of $R[F^e]$--modules. This implies
that if $M$ was a unit $R[F^e]$--module, then $S^{-1}M$ is not only a unit
$S^{-1}R[F^e]$--module but also a unit $R[F^e]$--module. Concretely, the
structure is given by
\[
    {\theta'}^e: R^e \tensor S^{-1}M \to[r \tensor \frac{m}{s} \mapsto
    \frac{r}{s^{p^e}}\theta^e(1 \tensor m)] S^{-1}M.
\]
\label{x.FonLoclzn}In terms of the Frobenius action ${F'}^e$ on
$S^{-1}M$ corresponding to ${\theta'}^e$ this can be expressed
quite simply as ${F'}^e(\frac{m}{s})=\frac{F^e(m)}{s^{p^e}}$,
where $F^e$ denotes the Frobenius action on $M$ corresponding to
$\theta^e$. This implies that the induced $R[F^e]$--structure on
$S^{-1}M$ is in fact uniquely determined by the requirement that
the localization map $M \to S^{-1}M$ be $R[F^e]$--linear
\label{x.RFuniqLoc}.

\subsubsection{Completion}\index{RF-mod@$R[F^e]$--mod!completion} Let $R$ be $F$--finite (\ie $R^e$ is a finitely
generated right $R$--module) and $R \to \wh{R}$ is the natural map to the
$I$-adic completion along some ideal $I$ of $R$. The natural map
\[
    R^e \tensor_R \wh{R} \to \wh{R}^e
\]
sending $r' \tensor r$ to $r'r^{p^e}$ is an isomorphism of
$R$--$\wh{R}$--bimodules. First, it is easy to observe that the indicated
map is linear as claimed. Secondly, observe that
\[
    R^e \tensor_R \wh{R}    = R^e \tensor_R \invlim \frac{R}{I^t}
                            \cong \frac{R^e}{R^eI^t} =
                            \frac{R^e}{I^{t[p^e]}R^e} \cong \wh{R}^e
\]
where we used the fact that for finitely generated $R$--modules $\wh{M}
\cong \wh{R} \tensor M$ \cite[Theorem 7.2]{Eisenbud.CommAlg}, and then
also that the sequence $I^{t[p^e]}$ is cofinal within the powers of $I$.
Now, let $(N,\theta^e)$ be a unit $\wh{R}[F^e]$--module. Then
\[
    R^e \tensor_R N = R^e \tensor_R \wh{R} \tensor_{\wh{R}} N = \wh{R}^e
\tensor_{\wh{R}} N \to[\theta^e] N
\]
is an isomorphism. Therefore, $N$ is naturally a unit $R[F^e]$--module.

\subsubsection{Field Extension}\index{RF-mod@$R[F^e]$--mod!field extension}
\note{Might not need the fields to be perfect?} Let $k \subseteq R$ be a
perfect field and let $K$ be a perfect extension field of $k$ (possibly
infinite). Let $R \to K \tensor_k R = R_K$ be the natural inclusion. As in
the complete case, we claim that the multiplication map $R^e \tensor_R R_K
\to R_K^e$ given by sending $r' \tensor r$ to $r'r^{p^e}$ is an
isomorphism of $R$--$R_K$--bimodules. The best way to observe this is to
give an inverse of this map as follows:
\[
    R_K^e  = (R \tensor_k K)^e \to R^e \tensor_k K = R^e \tensor_R (R
\tensor_k K).
\]
The equal signs are by definition or by the natural isomorphism $R
\tensor_k K \cong K \tensor_k R$. The arrow is given by sending $r \tensor
t$ to $r \tensor t^{1/p^e}$. This is easily checked to be a well defined
map ($K$ is assumed perfect), inverse to the multiplication map above.
Thus we conclude as in the cases of localization and completion above that
a unit $R_K[F^e]$--module $N$, indeed carries a natural unit
$R[F^e]$--module structure.

These three observations are worth being summarized in a separate
proposition as they play a distinguished role in our later treatments.
\begin{proposition}\label{prop.LocCompFieldex}
    Let\/ $R \to S$ be a map of rings which is either a localization of\/ $R$,
    a completion of\/ $R$ (in this case\/ $R$ is $F$--finite)
    or\/
    $S=K \tensor_k R$ with\/ $K \supset k$ a perfect extension field of the
    field perfect field\/ $k \subseteq R$. Then the forgetful functor
    \[
        S[F^e]\text{--mod} \to R[F^e]\text{--mod}
    \]
    restricts to a functor from unit\/ $S[F^e]$--modules to unit\/
    $R[F^e]$--modules. That is, every unit\/ $S[F^e]$--module is naturally a
    unit\/ $R[F^e]$--module.
\end{proposition}

\section{Finitely generated $R[F^e]$--modules}
The most powerful concept introduced in \cite{Lyub} is that of an
$F$--finite module. In our language, these correspond to $R[F]$--modules
which are finitely generated\index{RF-mod@$R[F^e]$--mod!finitely
generated} modules over the ring $R[F]$. Most of this section can be found
in \cite[Section 2]{Lyub} in the case that $R$ is regular. This seems to
be the case where the theory is most useful, thus we will frequently
assume that $R$ is regular in this section. Our definitions, however, are
set up to work fairly well in the non regular case; for regular $R$, they
coincide with the ones found in \cite[Section 2]{Lyub}.

If $M_0$ is an $R$--module with an $R$--linear map $\beta: M_0 \to
\F[e]M_0$, then we can use the Frobenius powers of this map to obtain a
directed system:
\[
   M_0 \to[\beta] \F[e]M_0 \to[\protect{\F[e]}\beta] \F[2e]M_0
                          \to[\protect{\F[2e]}\beta] \F[3e]M_0 \to \ldots \quad
\]
The direct limit $M$ of this system we call \emph{the unit
$R[F^e]$--module generated by
$\beta$}\label{x.DefGenerator}\index{$\beta$}. The $R$--module $M_0$ (or
the map $\beta$) we call a \emph{generator}\index{generator} for $M$. To
justify this notation we have to show that this limit $M$ does in fact
carry a natural unit $R[F^e]$--structure. Using that $\F[e]$ commutes with
direct limits, the following diagram of directed systems indicates this
natural unit structure on $M$. The first row represents $\F[e]M$ and the
second is $M$.
\[
\xymatrix@C=3pc{ {}&{\F[e]M_0} \ar^{\id}[d]\ar^{\F[e]\beta}[r]
&{\F[2e]M_0}
                          \ar^{\id}[d]\ar^{\F[2e]\beta}[r] &{\F[3e]M_0} \ar^{\id}[d]\ar[r]
                          &{\ldots}\\
   {M_0} \ar^{\beta}[r] &{\F[e]M_0} \ar^{\F[e]\beta}[r] &{\F[2e]M_0}
                          \ar^{\F[2e]\beta}[r] &{\F[3e]M_0} \ar[r]
                          &{\ldots}}
\]
These two limit systems are canonically isomorphic and therefore $M$ is a
unit $R[F^e]$--module. Obviously, every unit $R[F^e]$--module
$(M,\theta^e)$ has a generator, namely the inverse of its structural
morphism ${\theta^e}^{-1}: M \to \F[e]M$.

If the unit $R[F^e]$--module $M$ has a generator $(M_0,\beta)$ which is
finitely generated as an $R$--module, then $M$ is finitely generated as a
module over the ring $R[F^e]$. Indeed, the images of a set of $R$--module
generators of $M_0$ in the limit $M$ are $R[F^e]$--module generators for
$M$. Naturally one might ask whether the converse is true too; \ie if a
finitely generated unit $R[F^e]$--module has an $R$--finitely generated
generator. This is in fact the case as the corollary of the next
proposition shows.

\begin{definition}\label{def.root}
    Let $R$ be regular. A \emph{root}\index{root} of a unit $R[F^e]$--module $M$ is a
    generator $\beta: M_0 \to \F[e]M_0$ such that $\beta$ is injective and
    $M_0$ is a finitely generated $R$--module.
\end{definition}

The next task is to determine the unit $R[F^e]$--modules which have a
root. For this we observe the following general proposition.

\begin{proposition}\label{prop.fgIsStabfg}
    Let\/ $(M,\theta^e)$ be an\/ $R[F^e]$--module such that\/ $\theta^e$ is
    surjective. Then\/ $M$ is finitely generated as an\/ $R[F^e]$--module if
    and only if there is a finitely generated\/ $R$--submodule\/ $M_0$ such
    that\/ $M_0 \subseteq F^e(M_0)$ and\/ $M = \Union_n F^{en}(M_0)$.
\end{proposition}
\begin{proof}
    Let's begin with the ``only if'' direction: First note that the images
    of $R$--module generators of $M_0$ under $F^{ne}$ are $R$--generators
    of $F^{ne}(M_0)$. This implies that $M = \Union F^{ne}(M_0)$ is
    generated as an $R[F^e]$--module by the $R$--generators of $M_0$.

    Conversely, assume that $M$ is a finitely generated $R[F^e]$--module. Let $M'$ be
    the $R$--module generated by some finitely many $R[F^e]$--module
    generators of $M$. In other words,
    \begin{equation}\label{l.1}\tag{\ltag}
        R[F^e]M'=\sum_{n=0}^{\infty} F^{ne}(M') = M.
    \end{equation}
    Since $\theta^e$ is surjective $M = F^e(M)$. Applying $F^e$ to (\ref{l.1})
    we get $M = F^e(M) = F^e(\sum_{n=0}^\infty F^{ne}(M')) = \sum_{n=1}^\infty
    F^{ne}(M')$. Since $M'$ was finitely generated it is contained in a finite
    part of the above sum, say $M' \subseteq \sum_{n=1}^m F^{ne}(M')$. Now set
    $M_0 = \sum_{n=0}^{m-1} F^{ne}(M')$ and we see right away that $M_0
    \subseteq F^e(M_0)$. Iterating we get the following sequence of inclusions
    \[
        M_0 \subseteq F^e(M_0) \subseteq F^{2e}(M_0) \subseteq F^{3e}(M_0)
        \subseteq \cdots
    \]
    whose union is $M$ since $M_0$ contains $M'$.
\end{proof}

In the case that $R$ is regular, this is exactly what we need to show that
every finitely generated unit $R[F^e]$--module has a root:

\begin{corollary}
    Let\/ $R$ be regular and let\/ $(M,\theta^e)$ be a unit\/ $R[F^e]$--module.
    Then\/
    $M$ is finitely generated as an\/ $R[F^e]$--module if and only if\/ $M$ has a
    root.
\end{corollary}
\begin{proof}
    We apply the previous proposition and obtain $M$ as the increasing
    union $\Union F^{ne}(M_0)$ for some finitely generated $R$--submodule
    $M_0$ of $M$. Since for $R$ regular, $F^{ne}(M_0)$ is isomorphic to
    $\F[ne]M_0$ via $\theta^{ne}$, the union $M$ is isomorphic to the
    limit of the following directed system.
    \[
        M_0 \to[\beta] \F[e]M_0 \to[\protect{\F[e]\beta}] \F[2e]M_0
        \to[\protect{\F[2e]\beta}] \cdots
    \]
    Here $\beta$ is the restriction of ${\theta^e}^{-1}$ to $M_0$ thus in
    particular $\beta$ is injective. As $M_0$ is finitely generated we see that
    $M$ has a root.

    Conversely, if $\beta: M_0 \to \F[e]M_0$ is a root of $M$ by the
    exactness of $\F[e]$, all maps $\F[ne]\beta$ are injective as well;
    \ie all the maps in the directed system arising form $\beta$ whose
    limit is $M$ are injective. Thus we can identify $M_0$ with its image
    in $M$. Then obviously, $M_0 \subseteq F^e(M_0)$ and $M=\Union
    F^{ne}(M_0)$. With the previous proposition we conclude that $M$ is a
    finitely generated $R[F^e]$--module.
\end{proof}

This notion of root seems only useful when $R$ is regular, as then the
injectivity of $\beta$ implies the injectivity of all the maps in the
direct system arising from $\beta$, by the exactness of $\F[e]$. This is
crucial for reducing questions about a possibly very big $M$ to the
finitely generated root $M_0$ (cf.\ \cite[Section 1-3]{Lyub} where this
philosophy is applied with great success to various finiteness properties
of local cohomology modules). A candidate for a generalization of a root
in a possibly non--regular setting are finitely generated submodules $M_0$
of a unit $R[F^e]$--module $M$ such that $M_0 \subseteq F^{e}(M_0)$ and $M
= \Union F^{ne}(M_0)$ as in Proposition \ref{prop.fgIsStabfg}.
\begin{definition}\label{def.RootinGeneral}
    A finitely generated $R$--submodule $M_0$ of an $R[F^e]$--module $M$ is called a \emph{root}
    \index{root!non--regular} if $M_0 \subseteq F^e(M_0)$ and $M = R[F^e]M_0$.
\end{definition}
We don't know if this generalization of Definition \ref{def.root} will
lead to anything interesting but we will keep using it whenever the
regularity assumption is not needed in our discussions. With this new
notation Proposition \ref{prop.fgIsStabfg} shows that if $(M,\theta^e)$ is
a finitely generated $R[F^e]$-module such that $\theta^e$ is surjective
then $M$ has a root.

Now, let $R$ be regular. The last proof shows that a root of a unit
$R[F^e]$--module $M$ can be thought of as an $R$--finitely generated
submodule $M_0$ such that $M=\Union F^{ne}(M_0)$ and $M_0 \subseteq
F^{ne}(M_0)$. For $R$--submodules of a unit $R[F^e]$--module
$(M,\theta^e)$ over a regular ring $R$ we have that $\theta^e$ gives an
isomorphism between $\F[e]N$ and $F^e(N)$. Thus $F^e$ and $\F[e]$ are
isomorphic functors on submodules of $M$. The advantage of working with
$F^e$ instead of $\F[e]$ is that by doing so we never leave the ambient
module $M$. This makes many arguments much more transparent as it
simplifies notation considerably. We will take advantage of this
alternative point of view frequently.

\begin{remark}
Lyubeznik calls a unit $R[F]$--module $F$--finite if it is generated by
$\beta: M_0 \to \F M_0$ for some finitely generated $R$--module $M_0$. By
taking the image of $M_0$ in the unit $R[F^e]$--module it generates, it is
clear that every $F$--finite module has a root. Thus we just showed that
for regular rings $R$, an $R[F]$--module $M$ is $F$--finite (in the sense
of Lyubeznik \cite{Lyub}) if and only if $M$ is a finitely generated
$R[F]$--module (this is the case $e=1$ of the last corollary). In the
following, we will not use this meaning of $F$--finite since it clashes
with the standard usage of $F$--finiteness meaning that $R$ is a finitely
generated module over $R^p$, its subring of $p^{th}$--powers.
\end{remark}

The concept of a root is very useful to prove the basic properties of
finitely generated $R[F^e]$--modules.

\begin{proposition}\label{prop.fguRFSubs}
    Let\/ $R$ be regular and\/ $M$ be a finitely generated unit\/ $R[F^e]$--module\/ $M$ with
    root\/ $M_0$. If\/ $N$ is a unit\/ $R[F^e]$--submodule of\/ $M$, then\/ $N_0=N
    \intsec M_0$ is a root of\/ $N$. In particular,\/ $N$ is a finitely
    generated\/ $R[F^e]$--module.

    The unit\/ $R[F^e]$--submodules\/ $N$ of\/ $M$ are in an inclusion preserving
    one-to-one correspondence
    with the submodules\/ $N_0 \subseteq M_0$ such that\/ $F^e(N_0) \intsec M_0 =
    N_0$.
\end{proposition}
\begin{proof}
    As $M_0$ is a root, we have $M_0 \subseteq F^e(M_0)$ and $M = \Union
    F^{ne}(M_0)$. The same we have to show for $N$ and $N_0$. Since $N$ is
    unit we have $F^{ne}(N)=N$ for all $n>0$. Thus we get
    \[
        N_0=N \intsec M_0  \subseteq F^e(N) \intsec F^e(M_0)=F^e(N \intsec
        M_0)=F^e(N_0)
    \]
    while also using that for $R$ regular $F^e$ commutes with finite
    intersections. Similarly,
    \[
        N=N \intsec M=N \intsec \Union F^{ne}(M_0) = \Union F^{ne}(M_0 \intsec
        N) = \Union F^{ne}(N_0)
    \]
    which shows the first part of the proposition.

    For the second part note that $N_0 = N \intsec M_0$ has the desired
    property. Using that $M_0 \subseteq F^e(M_0)$ we get
    \[
        M_0 \intsec F^e(N_0) = M_0 \intsec F^e(M_0) \intsec N = M_0 \intsec N =
        N_0.
    \]
    Conversely, given any submodule $N_0$ of $M_0$ with the property $M_0 \intsec
    F^e(N_0) = N_0$, then especially $N_0 \subseteq F^e(N_0)$. Thus $N_0$
    generates a unit $R[F^e]$--submodule $N=\Union F^{ne}(N_0)$ of $M$. To see
    that $N \intsec M_0=N_0$, we first observe that by induction on $n$ one
    obtains from the case $n=1$, as assumed above, for all $n$ that $F^{ne}(N_0)
    \intsec M_0 = N_0$. Thus $N \intsec M_0=\Union F^{ne}(N_0) \intsec M_0 =
    \Union N_0 = N_0$.
\end{proof}

\begin{corollary}
    Let\/ $R$ be regular. The category of finitely generated unit\/
    $R[F^e]$--modules is a full abelian subcategory of u$R[F^e]$--\mod
    which is closed under extensions.
\end{corollary}
\begin{proof}
    We have to show that finite generation is passed to sub--objects and
    quotients. The submodule case is covered by the last Proposition
    \ref{prop.fguRFSubs}. The quotient case is obvious as, in general,
    quotients of finitely generated modules are finitely generated. This
    shows that finitely generated unit $R[F^e]$--modules form an abelian category.

    To check closedness under extension let $0 \to M' \to M \to M'' \to 0$ be
    an exact sequence with $(M',{\theta'}^e)$ and $(M'',{\theta''}^e)$
    $R[F^e]$--modules. By the Five Lemma \cite[page 169]{Lang.alg} applied to
    the diagram
    \[
        \xymatrix{
            {0}\ar[r]&{\F[e]M'}\ar[r]\ar_{{\theta'}^e}[d]&{\F[e]M}\ar[r]\ar@{-->}[d]
                         &{\F[e]M''}\ar[r]\ar^{{\theta''}^e}[d]&{0} \\
            {0}\ar[r]&{M'}\ar[r]&{M}\ar[r]&{M''}\ar[r]&{0}
                 }
    \]
    we see that the dotted arrow exists and is an isomorphism if $M'$ and
    $M''$ are unit; \ie $M$ is a unit $R[F^e]$--module. If $M'$ and
    $M''$ are finitely generated $R[F^e]$--modules, then so is $M$ as this
    property of modules over a ring passes to extensions in general.
\end{proof}
\begin{corollary}
    Let\/ $R$ be regular. A finitely generated unit\/ $R[F^e]$--module has the
    ascending chain condition in the category of unit\/ $R[F^e]$--modules.
\end{corollary}
\begin{proof}
Ascending chains of unit $R[F^e]$--submodules of $M$ correspond to
ascending chains of $R$--submodules of a root $M_0$ of $M$ by the second
part of Proposition \ref{prop.fguRFSubs}. As $R$ is noetherian and $M_0$
is finitely generated, these chains stabilize. Thus, also, chains of unit
$R[F^e]$--submodules stabilize.
\end{proof}
\begin{remark}
This noetherian property of the category of finitely generated unit
$R[F^e]$--modules makes it possible to transfer the standard proofs of
finiteness of invariants attached to a finitely generated $R$--module to
the much bigger class of finitely generated unit $R[F^e]$--modules. In
particular Lyubeznik shows \cite[Chapter 2]{Lyub} the finiteness of the
set of associated primes and Bass numbers for any unit $R[F^e]$--module.
\end{remark}

A crown jewel of the theory of unit $R[F^e]$--modules is Theorem 3.2 from
\cite{Lyub}:
\begin{theorem}[\protect{\cite[3.2]{Lyub}}]\label{thm.fguRFhaveDCC}
    Let\/ $R$ be a finitely generated algebra over a regular local ring.
    Then every finitely generated unit\/ $R[F^e]$--module has finite length
    as a unit\/ $R[F^e]$--module.
\end{theorem}
This means that the category of finitely generated unit $R[F^e]$--modules
also has the descending chain condition. If $R$ is complete and local the
proof follows with similar techniques as indicated below in the proof of
Proposition \ref{prop.MinimalRoot}. In particular, Lemma \ref{lem.Chev}
below is the key player. The main difficulty consists of reducing to the
local case. For this, the analog of Kashiwara's equivalence for
$R[F^e]$--modules \cite[Proposition 3.1]{Lyub} is important. One
interesting consequence of Theorem \ref{thm.fguRFhaveDCC} is that over a
complete local ring finitely generated unit $R[F^e]$--modules have a
\emph{unique} minimal root:
\begin{proposition}\label{prop.MinimalRoot}
    Let\/ $(R,m)$ be regular. The intersection of finitely many roots of a
    unit\/
    $R[F^e]$--module\/ $M$ is again a root of\/ $M$.

    If\/ $(R,m)$ is also local and complete and\/ $M$ is a finitely generated
    unit\/
    $R[F^e]$--module, then the intersection of all roots of\/ $M$ is also a root
    of\/ $M$; \ie $M$ has a unique minimal root\index{root!unique minimal}.
\end{proposition}
This is proved as Theorem 3.5 in \cite{Lyub}. We give a slightly different
argument here still using the key Lemma 3.3 of Lyubeznik which is a
straightforward generalization of Chevalley's\index{Chevalley} theorem
\cite[VIII, Theorem 13]{ZarSam2}:
\begin{lemma}[\protect{\cite[3.3]{Lyub}}]\label{lem.Chev}
    Let\/ $(R,m)$ be complete and local. Given a collection\/ $\cal{N}$ of
    submodules of a finitely generated\/ $R$--module\/ $M$ which is closed under
    finite intersections let\/ $N$ be the intersection of all the modules
    in\/
    $\cal{N}$. For every\/ $s \in \NN$ there is some\/ $N'\in \cal{N}$ such
    that\/
    $N' \subseteq N + m^sM$.
\end{lemma}
\begin{proof}[Proof of Proposition \ref{prop.MinimalRoot}]
For the first part, it is enough to show that the intersection of two
roots is again a root. We use the point of view of Definition
\ref{def.RootinGeneral} and think of roots $M_1$ and $M_2$ of $M$ as
finitely generated $R$--submodules of $M$ such $M_i \subseteq F^e(M_i)$
and $M = \bigcup F^{er}M_i$. Then clearly $M_1 \cap M_2 \subseteq F^e(M_1)
\cap F^e(M_2) = F^e(M_1 \cap M_2)$ since $F^e$ commutes with finite
intersections (see page \pageref{x.Fregcommutes}). Since $M_2$ is finitely
generated for some $r_0$, we have $M_2 \subseteq F^{er_0}M_1$. Thus,
$F^{er}(M_2) \subseteq F^{e(r+r_0)}(M_1 \cap M_2)$. Therefore $M=\bigcup
F^{er}(M_2) \subseteq \bigcup F^{er}(M_1 \cap M_2)$. This shows that $M_1
\cap M_2$ is a root of $M$.

We have to show that the intersection $N$ of all roots of $M$ is also a
root of $M$. Instead of taking the intersection over \emph{all} roots we
take an appropriate subset with the same intersection, namely all roots
inside of some given root $M_0$ of $M$. This collection of submodules of
the finitely generated $R$--module $M_0$ has the finite intersection
property by the first part. Thus we can apply Lemma \ref{lem.Chev} to pick
for each $i$ from the roots contained in $M_0$ one such that $N_i
\subseteq N + m^iM_0$. Obviously, $\bigcap N_i = N$. Using that each $N_i$
is a root of $M$, and therefore $N_i \subseteq F^e(N_i)$, we get
\begin{equation}\label{eqn.localMinRoot}
\begin{split}
    N = \bigcap N_i &\subseteq \bigcap F^{er}(N_i) \subseteq \bigcap F^{er}(N+m^iM_0) \\
                    &= F^{er}(N) + \bigcap m^{i[p^{er}]}F^{er}(M_0) \subseteq  F^{er}(N) + \bigcap m^{ip^{er}}F^{er}(M_0)\\
                    &= F^{er}(N).
\end{split}
\end{equation}
This shows that $N \subseteq F^{e}(N)$ and thus $N$ is the root of some
unit $R[F^e]$--submodule of $M$. This submodule is, of course, $N'
\defeq \bigcup F^{er}(N)$ and for it to be equal to $M$ we must show that
$F^{er}(N)$ contains a root of $M$ for some $r >0$. This isproveen by
induction on the length of $M$ as a unit $R[F^e]$--module which is finite
by Theorem \ref{thm.fguRFhaveDCC} as follows.

If $M$ is a simple unit $R[F^e]$--module (\ie $\length M = 1$), then
besides $N' = M$ the only other possibility is $N'=0$. In the first case
we are done so we assume that $N'=0$. Then $N$ is also zero which implies
that some $N_i=0$ as follows: By the Artin-Rees Lemma we find an integer
$t$ such that for all $s>t$
\[
    M_0\cap m^sF^e(M_0) \subseteq m^{s-t}(M_0 \cap m^tF^e(M_0)) \subseteq
    m^{s-t}M_0.
\]
Thus for some $s \gg 0$ we can find one of the $N_i$'s, such that $N_i
\subseteq m^sM_0$ but $N_i \not\subseteq m^{s+1}M_0$ (we are assuming that
$N=0$ then $N_i \subseteq m^iM_0$ by construction). Now
\[
    N_i \subseteq F^e(N_i) \cap M_0 \subseteq F^e(m^sM_0) \cap M_0
\subseteq m^{sp^e}F^e(M_0) \cap M_0 \subseteq m^{sp^e-t}M_0
 \]
which is a contradiction since $sp^e-t > s+1$ for $s \gg 0$. Thus $N_i=0$
but this is absurd since $N_i$ is a root of $M$. This finishes the case
$\length M=1$.

For the induction step, let
\[
    0 \to M' \to M \to[\pi] M'' \to 0
\]
be an exact sequence of nonzero unit $R[F^e]$--modules. Then the length of
$M'$ and $M''$ as unit $R[F^e]$--modules is strictly smaller than the
length of $M$. By induction we can assume that $M'$ as well as $M''$ have
unique minimal roots $N'$ and $N''$ respectively. Let $r>0$ be such that
$M_0 \cap M' \subseteq F^{er}(N')$ and $\pi(M_0) \subseteq F^{er}(N'')$.
This is possible since $N'$ and $N''$ are roots and $M_0$ is finitely
generated. For an arbitrary root $M_1$ of $M$ it is easy to check that
$M_1 \cap M'$ and $\pi(M_1)$ are roots of $M'$ and $M''$ respectively.
Since $N'$ and $N''$ are \emph{unique minimal} roots we get $N' \subseteq
M_1 \cap M'$ and $N'' \subseteq \pi(M_1)$. Thus, in particular, we have
\[
    M_0 \cap M' \subseteq F^{er}(M_1) \cap M'\ \text{ and }\
    \pi(M_0)\subseteq \pi(F^{er}(M_1)).
\]
This implies that $M_0 \subseteq F^{er}(M_1)$ for all roots $M_1$ of $M$.
Note that $r$ does not depend on $M_1$ and therefore this inclusion holds
simultaneously for all roots $N_i$ in the definition of $N$; therefore
\[
    M_0 \subseteq \bigcap F^{er}(N_i) \subseteq F^{er}(N)
\]
where the second inclusion is part of equation \eqnref{eqn.localMinRoot}.
This shows that, in fact, the $R[F^e]$--module generated by $N$ is all of
$M$, \ie $N$ is the unique minimal root of $M$.
\end{proof}

\subsection{Generation and base change}
We briefly summarize some results on finite generation and roots under
base change. The following basic observation is very useful.
\begin{proposition}
    Let\/ $R \to S$ be a map of rings. Let\/ $M$ be a finitely generated
    unit\/
    $R[F^e]$--module with root\/ $M_0$. Then the image of\/ $S \tensor M_0$ in\/ $S
    \tensor M$ is a root of the finitely generated unit\/ $S[F^e]$--module\/ $S
    \tensor M$.

    If\/ $R \to S$ is flat, then\/ $S \tensor M_0$ itself is the root of\/ $S \tensor
    M$. If\/ $R \to S$ is faithfully flat, then a submodule\/ $M_0$ of\/ $M$ is a
    root of\/ $M$ if and only if\/ $S \tensor M_0$ is a root of\/ $S \tensor M$.
\end{proposition}
\begin{proof}
The first statement follows from the fact that tensor commutes with
Frobenius and direct limits. Using our alternative definition of root in
Definition \ref{def.RootinGeneral}, we have $M = \bigcup F^{er}(M_0)$.
Applying $S \tensor \usc$ we get $S \tensor M = \dirlim F^{er}(S \tensor
M)$. All members of this direct limit come with a natural map to $S
\tensor M$, obtained from the inclusion $F^{er}(M_0) \subseteq M$. Thus,
taking images in $S \tensor M$, and denoting the image of $S \tensor M_0$
in $S \tensor M$ by $\bar{M}_0$, the direct limit becomes a union over all
$F^{er}(\bar{M}_0)$. In particular, $\bar{M}_0 \subseteq F^e(\bar{M}_0)$,
thus $\bar{M}_0$ is a root of $S \tensor M$.

The flatness of $S$ ensures that $S \tensor M_0$ is a submodule of $S
\tensor M$. Thus, by the first part, $S \tensor M_0$ is a root. If $S
\tensor M_0$ is a root, then $S \tensor M_0 \subseteq S \tensor F^e(M_0)$.
Assuming faithfully flatness of $R \to S$ this ensures that $M_0 \subseteq
F^e(M_0)$ and $M = \bigcup F^{er}(M)$. Therefore $M_0$ is a root of $M$.
\end{proof}

With this pleasant behavior of roots under completion we can derive the
following corollary of the theorem on the existence of a unique minimal
root to the non--complete case.
\begin{corollary}
    Let\/ $R$ be a regular local ring. Let\/ $M$ be a finitely generated
    unit\/ $R[F^e]$--module. If the unique minimal root of the finitely
    generated unit\/ $\wh{R}[F^e]$--submodule\/ $\wh{R} \tensor M$ is extended
    from a submodule\/ $N \subseteq M$, then\/ $N$ is the unique minimal root of
    $M$.
\end{corollary}
\begin{proof}
    First, by faithfully flat descent $N$ is a finitely generated $R$--module.
    Since Frobenius commutes with tensor, $F^e(\wh{R} \tensor
    N) = \wh{R} \tensor F^e(N)$ and therefore $N \subseteq F^e(N)$. Thus
    $N$ is the root of some $R[F^e]$--submodule of $M$. For the same
    reason (completion commutes with Frobenius), the unit $R[F^e]$--module
    generated by $N$ must be all of $M$, and therefore, $N$ is a root of
    $M$.

    Secondly, if $N_0$ is some root of $M$, then $\wh{R} \tensor N_0$ is a
    root of $\wh{R} \tensor M$. Thus $\wh{R} \tensor N \subseteq \wh{R} \tensor
    N_0$ since the first one is the unique minimal root. But, by faithfully
    flatness of completion, this implies that $N \subseteq N_0$. Therefore
    $N$ is contained in any root of $M$ and is a root itself. Thus $N$ is
    the unique minimal root of $M$.
\end{proof}
This proof works more generally whenever $R \to S$ is a faithfully flat
extension. If for a unit $R[F^e]$--module $M$, the unit $S[F^e]$--module
$S \tensor M$ has a unique minimal root, which is extended from a
submodule $N \subseteq M$, then $N$ is the unique minimal root of
$M$.\footnote{This can be compared to the fact proved later that the
existence of a unique simple $D_R$--submodule descends down from a
faithfully flat extension, see Lemma \ref{lem.descend}}

Next we focus on the behavior of finite generation for $R[F^e]$--modules
under localization. Let $M$ be a (unit) $R[F^e]$--module. We saw before
that for $x \in R$ the localization $M_x$ naturally is a (unit)
$R[F^e]$--module. For any $R$--submodule $N \subseteq M_x$ we have that
\[
    F^{ne}(x^{-1}N)=x^{-p^{ne}}F(N).
\]
Here, we denote the $R$--submodule generated by the set $\set{\, yn\, |\,
n \in  N\, }$, for some $y \in R_x$, with $yN$. The equation follows from
the description of the Frobenius action on a localization as on page
\pageref{x.FonLoclzn}.

Now let $M$ be such that it has a root $M_0$, \ie $M$ satisfies the
assumptions of Proposition \ref{prop.fgIsStabfg}. Then $M_0$ is a finitely
generated $R$--submodule of $M$ with $M_0 \subseteq F^{e}(M_0)$ and $M =
\Union F^{ne}(M_0)$. Then $M'_0
\defeq x^{-1}M_0$ is a root of the $R[F^e]$--module $M_x$. First,
\[
    M'_0 = x^{-1}M_0 \subseteq x^{-p^e}F^{e}(M_0) = F^e(x^{-1}M_0) =
    F^e(M'_0).
\]
Secondly, $M_x = \Union x^{-p^{te}}M = \Union x^{-p^{te}} (\Union
F^{ne}(M_0))$. By a diagonal trick the last union is just $ \Union
x^{-p^{ne}}F^{ne}(M_0) = \Union F^{ne}(M'_0)$. Thus $M_x = \Union
F^{ne}(M'_0)$. We just showed:
\begin{lemma}
    Let\/ $M$ be an\/ $R[F^e]$--module which has a root\/ $M_0$ (\eg $M$ is
    finitely generated and unit) and let\/ $x \in R$. Then\/ $x^{-1}M_0$ is a root of
    the\/
    $R[F^e]$--module\/ $M_x$.
\end{lemma}
This lemma immediately implies that for a regular ring $R$ all local
cohomology modules $H^i_I(M)$ of a finitely generated unit
$R[F^e]$--module $M$ are also finitely generated unit $R[F^e]$--modules.
Thus, within the category of unit $R[F^e]$--modules, the local cohomology
modules satisfy the ascending chain condition and even have finite length
by Theorem \ref{thm.fguRFhaveDCC} in this category.

If we compare the last lemma with the results at the beginning of this
chapter on the behavior of unit $R[F^e]$--modules under localization and
completion, we are led to ask whether the localizations/completion of a
finitely generated unit $R[F^e]$--module is again finitely generated as
such. This is probably not the case for the completion by cardinality
reasons, for example: If $(R,m)$ is essentially of finite type over a
field $k$, then $R[F^e]$ is countably generated over $k$. Thus, every
finitely generated $R[F^e]$--module is countably generated over $k$, but
$\wh{R}$ is uncountably generated over $k$. In the case of localization at
more than just finitely many elements of $R$ one also, quite likely,
looses finitely generatedness over $R[F]$.

%%%%%%%%%%%%%%%%%%%%%%%%%%%%%%%%%%%%%%%%%%%%%%%%%%%%%%%%%%%%%%
%%                                                          %%
%%   This is file: chapter3.tex                             %%
%%   It contains the Second Chapter of my                   %%
%%   dissertation: dissertation.tex                         %%
%%                                                          %%
%%%%%%%%%%%%%%%%%%%%%%%%%%%%%%%%%%%%%%%%%%%%%%%%%%%%%%%%%%%%%%

\chapter{$D_R$--submodules of unit $R[F]$--modules}\label{chap.DRsub}

This chapter is dedicated to the proof of the first Main Theorem of this
dissertation. This in turn is one of the key ingredients for the result in
the title of the dissertation.

\mainthmonecontent \index{Main Theorem 1}

In order to achieve this, several techniques have to be developed. First,
we review the basic properties of differential operators in finite
characteristic. Most importantly this is the filtration of $D_R$, the ring
of differential operators on $R$, by subrings $D_R^{(e)}$ consisting of
differential operators which are linear over $R^{p^e}$ (cf.\
\cite{SmithSP.diffop, Yeku}). The connection between unit
$R[F^e]$--modules and $D_R$--modules is explored through a version of
Frobenius descent \cite{Ber.FrobDesc} in the special case of submodules of
a unit $R[F^e]$--modules. This Frobenius descent also appears in the proof
of Lyubeznik's Lemma 5.4 although somewhat implicitly. Our explicit
description turns out to be quite useful as a way to guide the intuition
as well as in the proof of several important steps of Main Theorem 1.

Next we describe the endomorphisms of a simple $D_R$--module. We are able
to show that under the assumption that $k$ is uncountable, $\End_{D_R}(N)$
is algebraic over $k$ for a simple $D_R$--module $N$. We believe that this
is true even without the somewhat annoying assumption that the field $k$
be uncountable. If this were known we could immediately drop the same
assumption in Main Theorem 1.\footnote{During my dissertation defense
Brian Conrad suggested that the smearing out techniques of EGA4.3
\cite{EGA4.3} could be used to avoid the assumption of uncountability of
the field $k$ in Main Theorem 1. It seems that this, in fact, is the case
and a careful treatment of this is in preparation.}

The last section of this chapter is devoted to an example showing the
necessity of the assumption of algebraic closure in Main Theorem 1. This
provides a counterexample to Lyubeznik's \cite[Remark 5.6a]{Lyub} where he
claims Main Theorem 1 without assuming that $k$ be algebraically closed.

\section{$D_R$--modules in finite characteristic}\index{DR-module@$D_R$--module}

An important technique in the study of the ring of differential operators
$D_R$ of a regular ring $R$ in characteristic zero is to use the fact that
$D_R$ is noetherian. In characteristic $p > 0$ this is not true and
different techniques have to be developed. Some kind of finiteness is
imposed by employing the Frobenius. One first observes that quite
generally $D_R$ admits a filtration by $R$--algebras $\D[e]_R$ which
consists of operators linear over the subring $R^{p^e}$ of $R$. This
observation, and a resulting \idx{Morita equivalence} of all the
categories of $\D[e]_R$--modules for all $e$ in the case that $R$ is
regular and $F$--finite, was successfully applied (cf.\ S.P. Smith
\cite{SmithSP.diffop,SmithSP:DonLine} and Haastert
\cite{Haastert.DiffOp,Haastert.DirIm} for example) to show properties of
$D_R$ in finite characteristic.

We will review the definition and basic properties of differential
operators and show how unit $R[F^e]$--modules naturally carry a
$D_R$--module structure; all the material on $D_R$--modules can be found
in \cite{EGA4}, alternatively we recommend \cite{Trav.Phd} for a pleasant
introduction. When we get to the connection with unit $R[F^e]$--modules,
\cite{Lyub} is the source we follow. We begin with the definition of the
ring of differential operators.
\begin{definition}
Let $k$ be a commutative ring and let $R$ be a $k$--algebra (of arbitrary
characteristic). The \emph{ring of $k$--linear \idx{differential
operator}s $D_{R/k}$} is a subring of the $k$--linear
endomorphisms $\End_k(R)$ of $R$. \idx{$D_{R/k}$} is the
union $\Union D^n_{R/k}$ where the \idx{$D^n_{R/k}$} are defined
inductively by
\begin{eqnarray*}
 D^{-1}_{R/k} &=& 0 \\
 D^{n+1}_{R/k} &=&\{\phi \in \Hom_k(R,R) | [\phi,r] \in D^n_{R/k} \text{ for
all } r \in
 R\}.
\end{eqnarray*}
The elements of $D^n_{R/k}$ are called ($k$--linear) differential operators of
degree at most $n$.
\end{definition}
Here $[\phi,\psi]=\phi\circ\psi-\psi\circ\phi$\index{$[\phi,\psi]$}
denotes the commutator. If $k$ is understood it is often omitted from the
notation and we just write $D_R$ for $D_{R/k}$. We identify $R$ with
$D^0_R=\End_R(R)$. Also observe that if we have $k \subseteq k' \subseteq
R$, then $D_{R/k'} \subseteq D_{R/k}$ since $\End_{k'}(R) \subseteq
\End_k(R)$.

The following example of the polynomial ring over a field demonstrates the
different behavior of $D_R$ in characteristic zero and characteristic
$p>0$.

\begin{example}\label{ex.DofPoly}
If $R=k[\xn]$ is the polynomial ring over a perfect field $k$ of arbitrary
characteristic, then the ring of $k$--linear differential operators $D_R$
is generated as an $R$--module by the divided power operators
$\frac{\partial^j}{j!\partial x_i^j}$ (cf.\ \cite[Chapter 5]{Lyub} or
\cite{Trav.Phd}). We have to make a remark on the meaning of the
expressions $\frac{\partial^{j}}{j!\partial x_i^{j}}$. These are
understood to first carry out $\frac{\partial^{j}}{\partial x_i^j}$, then
see if anything cancels with $\frac{1}{j!}$ and only then remember that we
were in finite characteristic (if this was the case). So, for example,
$\frac{\partial^p}{p!\partial
x_i^p}(x_jx_i^p)=\text{``}\frac{1}{p!}p!x_j\text{''}=x_j$ even in
characteristic $p
> 0$. Whereas in characteristic $p$ one has that
$\frac{\partial^j}{\partial x_i^j}$ is the zero operator for $j \geq p$.

If the characteristic of $k$ is zero, $D_R$ is a finitely generated
$R$--algebra with generators $\frac{\partial}{\partial x_i}$, the partial
derivations. This is no longer true if the characteristic is $p>0$ since,
for example, $\frac{\partial^{p}}{{p}!\partial x_1^{p}}$ is a differential
operator which is not in the algebra of derivations. In fact, there is an
increasing chain of $R$--subalgebras of $D_R$
\[
    \D[0]_R \subseteq \D[1]_R \subseteq \D[2]_R \subseteq \D[3]_R \subseteq
\cdots
\]
whose union is $D_R$. $\D[0]_R$ is just $R$, $\D[1]_R$ is the $R$--algebra
generated by $1$ and the derivations; \ie $\D[1]_R = R \oplus \op{Der}_{R/k}$.
More generally, $\D[e]_R$ is the
$R$--algebra generated by the divided power operators
$\frac{\partial^{p^a}}{p^a!\partial x_i^{p^a}}$ with $a < e$.

Similar observations can be made for any regular ring $R$ essentially of
finite type over a perfect field by observing that locally $R$ is essentially
\'etale over a polynomial ring, and that $D_R$ is well behaved under essentially
\'etale covers (cf. \cite[Theorem 3.2.5] {Trav.Phd}). That is, if $R'$ is
essentially \'etale over $R$, then the map of left $R$--algebras
\[
    R' \tensor_{R} D_{R} \to D_{R'}
\]
sending $r \tensor \delta to r\delta$ is an isomorphism.
\end{example}

Next we give an alternative but completely equivalent description of
$D_{R/k}$, see \cite{EGA4,Trav.Phd} for proofs. If we denote by
$J=J_{R/k}$ the ideal in $R \tensor_k R$ which is the kernel of the
multiplication map $R \tensor_k R \to R$, then
\[
    D^n_R = \Ann_{\End_k(R)}(J^{n+1}).
\]
Here we consider the natural action of $R \tensor_k R$ on $\End_k(R)$ by
pre- and post-composing with elements in $R$; i.e $(r \tensor r')(\phi)(x)
= r\phi(r'x)$. As an ideal, $J$ is generated by $\set{\,r \tensor 1 - 1
\tensor r \,|\,r \in R\,}$ where it is in fact enough to range over a set
of $k$--algebra generators of $R$.

Let us now assume that $R$ is of finite characteristic $p>0$. An important
observation for us is that $\phi \in \End_k(R)$ is linear over $R^{p^e}$
if and only if $J^{[p^e]}$--kills $\phi$, \ie
\[
    \End_{R^{p^e}}(R) = \Ann_{\End_{k}(R)}(J^{[p^e]}).
\]
This is clear since $J^{[p^e]}$ is generated by $\set{r^{p^e} \tensor 1 -
1 \tensor r^{p^e} \ |\ r \in R}$ and such elements kill $\phi$ if and only
if $\phi$ commutes with $r^{p^e}$ by definition of the action of $R
\tensor_k R$ on $\End_k(R)$. Since $J^{[p^e]} \subseteq J^{p^e}$ this
shows that for $n < p^e$ we get $D^n_R \subseteq \End_{R^{p^e}}(R)$. Thus
\[
   D_R \subseteq \bigcup_{q=p^e} \End_{R^q}(R)
\]
expresses $D_R$ as a union of subrings $D_R \intsec \End_{R^{p^e}}(R)$. In
analogy with the above Example \ref{ex.DofPoly} we define $\D[e]_R$ as
$D_R \intsec \End_{R^{p^e}}(R)$, \ie $\D[e]_R$ are those differential
operators linear over $R^{p^e}$. The elements of $\D[e]_R$ we call
differential operators of \emph{level} $e$. With an additional hypotheses we
get (cf.\ \cite{Yeku}):

\begin{proposition}
    Let\/ $R$ be a finitely generated algebra over\/ $R^p$ and\/ $k \subseteq R$
    be perfect. Then\/ $\D[e]_{R/k} = \End_{R^{p^e}}(R)$ and therefore
    \begin{equation}\label{eqn.prop.DrEnd}
        D_{R/k} = \Union \End_{R^{p^e}}(R).
    \end{equation}
\end{proposition}
\begin{proof}
It remains to show that $\End_{R^{p^e}}(R)$ is contained in $D_{R/k}$. We
show the stronger condition that $\End_{R^{p^e}}(R) \subseteq
D_{R/R^{p^e}}$ (here we use that $k$ is perfect as this ensures that $k
\subseteq R^{p^e}$ and thus $D_{R/R^{p^e}} \subseteq D_{R/k}$). Let
$\xn[r]$ be the finitely many algebra generators of $R$ over $R^{p}$.
These generate also $R$ as an algebra over $R^{p^e}$. The ideal
$J=J_{R/R^{p^e}}$ is therefore generated by $\set{r_i \tensor 1 - 1
\tensor r_i\ |\ i=1 \ldots n\ }$ which is a finite set of $n$ elements.
This implies that $J^{np^e} \subseteq J^{[p^e]}$ thus we have the reverse
inclusion for their annihilators: $\End_{R^{p^e}}(R) \subseteq
D^{np^e-1}_{R/R^{p^e}}$.
\end{proof}
Note that in the right hand side of Equation \eqnref{eqn.prop.DrEnd} in
this proposition $k$ does not appear. This implies that $D_{R/k}$ is  the
same for every perfect $k \subset R$. Consequently we write just $D_R$
instead of $D_{R/k}$ for $k$ perfect. In particular this shows the well
known fact that for $k$, a perfect field, $D_{k/\FF_p}=D_{k/k}=k$.

%Using the Frobenius descent techniques developed later in this chapter,
%one can show that the ad hoc notation of $\D[e]$ in Example
%\ref{ex.DofPoly} was correct. One shows that, if $R$ is regular and of
%finite type over a perfect field, then $\D[e]_R$ is the $R$--algebra
%generated by the operators of degree strictly smaller than $p^e$.
%\note{proof this later and then refer to it from here}

Left modules over $D_R$ are called $D_R$--modules. The category of
$D_R$--modules is as usual denoted by $D_R$--\mod. $D_R$--modules behave
nicely under localization as the following example shows.

\begin{example}\label{ex.DRuniqLoc}
The ring $R$ itself is a $D_R$--module; the action of $\delta \in D_R$ on
$R$ is via its action as an element of $\delta \in \End_k(R)$. If $M$ is a
$D_R$--module, then any localization $S^{-1}M$ carries a unique
$D_R$--module structure such that the natural localization map $M \to
S^{-1}M$ is $D_R$--linear. The operator $\delta \in D_R$ is linear over
$R^{p^e}$ for $e \gg 0$. For such $e$, we define
\begin{equation*}\label{l.DRonLoc1}
\delta\left(\frac{m}{s}\right) =\delta\left(\frac{s^{p^e-1}m}{s^{p^e}}\right)=
\frac{\delta(s^{p^e-1}m)}{s^{p^e}}.
\end{equation*}
This is independent of $e$ since for $c>e$ we have $\delta(s^{p^c-1}m)=
\delta((s^{p^{c-e}-1})^{p^e}s^{p^e-1}m)=
(s^{p^e})^{p^{c-e}-1}\delta(s^{p^e-1}m)$. The independence of the
representatives $s \in S$ and $m \in M$ for $\frac{m}{s}\in S^{-1}M$ is
equally easy to show. To see uniqueness, observe that
\begin{equation}\tag{\ltag}\label{l.DRonLoc2}
    s\delta\left(\frac{m}{s}\right) =
{[s,\delta]\left(\frac{m}{s}\right)-\delta(m)}
\end{equation}
and therefore the action of $\delta$ on $\frac{m}{s}$ is determined by the
lower degree operator $[s,\delta]$. Since the $R$ module structure on
$S^{-1}M$ is uniquely determined we conclude by induction on the degree of
a differential operator that the $D_R$--structure is also unique. Note
that one can alternatively define the $D_R$--structure on $S^{-1}M$ by
(\ref{l.DRonLoc2}) which is a characteristic independent point of view.
\end{example}

\subsection{Unit $R[F^e]$--modules are $D_R$--modules}\label{sec.UnitRFareD}
We developed all insights to construct a $D_R$--module structure on unit
$R[F^e]$--modules; most importantly we showed that $D_R \subseteq \Union
\End_{R^{p^e}}(R)$. Note that we can alternatively think of
$\End_{R^{p^e}}(R)$ as the right $R$--module endomorphisms of $R^e$; \ie
\[
    \End_{R^{p^e}}(R) = \End_\modR(R^e).
\]
With this point of view, modules of the form $R^e \tensor M$ for an
$R$--module $M$ carry a canonical $\D[e]_R$--module structure:
\begin{proposition}
    The Frobenius functor\/ $\F[e]$ is naturally a functor from\/ $R$--modules
    to\/ $\D[e]_R$--modules. If\/ $M$ is an $R$--module, then\/ $\delta \in
    \D[e]_R$ acts on\/ $\F[e]M=R^e \tensor M$ as the map\/ $\delta \tensor
    \id_M$, where we interpret\/ $\delta$ as an element of\/
    $\End_{\modR}(R^e)$ via the inclusion\/ $\D[e]_R \subseteq
    \End_{\modR}(R^e)$.
\end{proposition}
This also makes it clear how one endows a unit $R[F^e]$--module
$(M,\theta^e)$ with a natural $D_R$--module structure. If we define
inductively \label{pg.deftheta}
\[
    \theta^{e(r+1)} \defeq \theta^e \circ \F[e](\theta^{er}) = \theta^{er}
    \circ \F[er](\theta^e)
\]
we see that $\theta^{er}:R^{er} \tensor M \to[\cong] M$ is an isomorphism
of left $R$--modules. By the previous proposition $R^{er} \tensor M$ is
naturally a $\D[er]_R$--module and thus, via the isomorphism
$\theta^{er}$, so is $M$. Concretely, $\delta \in \D[er]_R$ acts as the
map $\theta^{er} \circ (\delta \tensor \id_M) \circ ({\theta^{er}})^{-1}$
on $M$. This is illustrated by the following diagram.
\[
\xymatrix@C=4.5pc{
    {R^{er} \tensor M} \ar^(.6){(\theta^{er})^{-1}}@{<-}[r]\ar_{\delta \tensor \id_M}[d] &{M} \ar^{\delta}@{-->}[d] \\
    {R^{er} \tensor M} \ar^(.55){\theta^{er}}[r]       &{M} \ }
\]
The dotted vertical arrow stands for the action of $\delta$ on $M$. It is
defined by following the diagram along the path to the left. Since this
works for all $r$, \ie for all $\delta \in \D[er]_R$ we get:
\begin{proposition}\label{prop.wellDefDRonUnit}
    Let\/ $(M,\theta^e)$ be a unit\/ $R[F^e]$--module. The unit\/
    $R[F^e]$--structure\/ $\theta^e$ of\/ $M$ induces a natural\/ $D_R$--module
    structure on\/ $M$.
\end{proposition}
\begin{proof}
In the preceding discussion we describe how this structure is constructed.
It remains to show the well definedness of this definition. We have to
show that the action of $\delta \in \D[er]_R$ is the same when we think of
$\delta$ as an element of $\D[er']_R$ for some $r' \geq r$. Inductively it
is enough to assume that $r'=r+1$. The following diagram illustrates the
scene.
\begin{equation}\label{eqn.WellDefDRonUnit}
\begin{split}
    \xymatrix@C=4pc{
        {R^{e(r+1)}\tensor M}\ar^(.56){\F[er]({\theta^e})^{-1}}@{<-}[r]\ar_{\delta \tensor \id_{R^e\tensor M}}[d]  &{R^{er}\tensor M}\ar^(.6){({\theta^{er}})^{-1}}@{<-}[r]\ar_{\delta \tensor \id_M}[d]    &{M}\ar@{-->}[d] \\
        {R^{e(r+1)}\tensor M}\ar^(.55){\F[er](\theta^e)}[r]
    &{R^{er}\tensor M}\ar^(.6){\theta^{er}}[r]          &{M} }
\end{split}
\end{equation}
Following the dotted arrow by going around the right square corresponds to
the action of $\delta$ when thinking of it as an element of $\D[er]_R$.
Walking around all the way corresponds to $\delta \in \D[e(r+1)]_R$. Thus,
it is enough to observe that the left hand square is commutative. This is
clear if we note that the identification $R^{e(r+1)} \cong R^{er} \tensor
R^e$ exhibits $\delta \in \End_\modR(R^{er})$ as the element $\delta
\tensor \id_{R^e} \in \End_\modR(R^{e(r+1)})$. This justifies the notation
$\delta \tensor \id_{R^e \tensor M}$ for the leftmost vertical map of
\eqnref{eqn.WellDefDRonUnit} and since $\F[er]\theta^e = \id_{R^{er}}
\tensor \theta^e$ it is clear that the left square commutes.
\end{proof}

Note that this proof also shows that we get the same $D_R$--module
structure on $M$ whether we view $(M,\theta^e)$ as a unit $R[F^e]$--module
or as the unit $R[F^{er}]$--module $(M,\theta^{er})$ for some $r > 0$.
Therefore, this process works in fact for unit $R[F^{\infty}]$--modules;
thus every unit $R[F^\infty]$--module carries a natural $D_R$--module
structure.

The examples of $R[F^e]$--modules as well as $D_R$--modules we considered
so far are the ring $R$ itself, various localizations of $R$ and kernels
and cokernels of maps of such modules, the local cohomology modules among
them. These objects all carry \emph{natural} $D_R$ and $R[F^e]$--module
structures, essentially because $R$ carries a natural such structure. We
have to make sure that the natural unit $R[F^e]$--module structure, say on
local cohomology modules, for example, induces the natural $D_R$--module
structure. In order to show this, it is convenient to temporarily denote
by $\xi(M)$, the $D_R$--module $M$ with structure induced by $\theta^e$.
We start by observing that the natural $R[F]$--structure $\theta$ on $R$
induces its natural $D_R$--structure. Let us denote the action of $\delta
\in \D[e]_R \subseteq \End_k(R)$ induced by $\theta$ by $\delta_\theta$.
It is given by the dotted arrow making the diagram commute:
\[
    \xymatrix@C=3pc{
        {R^e \tensor R} \ar@{<-}^(.65){({\theta^{e}})^{-1}}[r]\ar_{\delta\tensor \id_R}[d] &{R} \ar@{-->}^{\delta_\theta}[d] \\
        {R^e \tensor R} \ar^(.6){\theta^{e}}[r]       &{R} }
\]
Therefore $\delta_\theta(r)=\theta^e(\delta\tensor \id_R(r \tensor
1))=\delta(r)$. This is indeed the natural action of $D_R$ on $R$. With the
next proposition we get the desired compatibility.
\begin{proposition}
    The association\/ $\ \xi: uR[F^e]\text{--\mod} \to D_R\text{--\mod}$ is
    an exact functor which commutes with localization.
\end{proposition}
\begin{proof}
First we have to show that a map $f: M \to N$ of unit $R[F^e]$--modules is
in fact $D_R$--linear, \ie it is a map $f: \xi(M) \to \xi(N)$ of
$D_R$--modules. This means that for all $\delta \in \D[er]_R$ we have to
show $\delta_N \circ f = f \circ \delta_M$ where $\delta_N$ (resp.
$\delta_M)$ denotes the action of $\delta$ on $\xi(N)$ (resp. $\xi(M)$).
But this is equivalent to the obvious equality $(\delta \tensor \id_N)
\circ (\id_{R^e} \tensor f) = (\id_{R^e} \tensor f) \circ (\delta \tensor
\id_M)$ as the corresponding diagrams are isomorphic via $\theta^{er}_M$
and $\theta^{er}_N$ (the reader should sketch these diagrams). Thus $\xi$
is a functor.

It is exact since $\xi$ does not change the underlying $R$--structure and
a sequence of unit $R[F^e]$--modules is exact if and only if it is exact
as $R$--modules, which is also the criterion for exactness as
$D_R$--modules. To see that $\xi$ commutes with localization just observe
that the induced $R[F^e]$--module (resp. $D_R$--module) structure on
$S^{-1}M$ is, by page \pageref{x.RFuniqLoc}, the unique
$R[F^e]$--structure (resp. $D_R$--structure by Example \ref{ex.DRuniqLoc})
for which the localization map $M \to S^{-1}M$ is $R[F^e]$--linear (resp.
$D_R$--linear).
\end{proof}
Since $\xi$ is so well behaved we will now drop it from the
notation. Thus, as before, if $(M,\theta^e)$ is an
$R[F^e]$--module, then we refer to the $D_R$--module $\xi(M)$ just
by $M$. After all, they are the same $R$--modules. Later we will
consider the $D_R$--structures induced from different
$R[F^e]$--structures $\theta^e$ and ${\theta'}^e$ on $M$. In this
case we will denote the action of $\delta \in D_R$ induced by
$\theta^e$ with $\delta_{\theta^e}$ and the action induced by
${\theta'}^e$ with $\delta_{{\theta'}^e}$.

We have seen that $uR[F^e]$--\mod is an (abelian if $R$ is regular)
subcategory of $D_R$--\mod. It is important to keep in mind that it is
\emph{not} a full subcategory. Even $uR[F^\infty]$--\mod is not full in
$D_R$--\mod as Example \ref{ex.DRnorRFinfty} below indicates. It is not
full since $D_R$--submodules of unit $R[F^e]$--modules need not be fixed
by any power of $F^{er}$.

\section{Frobenius Descent and Morita equivalence}\label{sec.FrobDesc}

As we saw in the last section, one can express the ring of differential
operators $D_R$ as the increasing union of subrings $\D[e]_R$ consisting
of the operators that are linear over $R^{p^e}$. This is true quite
generally (whenever $R$ is a finitely generated $R^p$--algebra) and proved
useful for endowing unit $R[F^e]$--modules with a $D_R$--module structure.
If $R$ is regular, further techniques along these lines are possible. Most
importantly, this is the fact that the Frobenius functor gives an
equivalence between the categories $\D[e]$--mod for various $e$. This
result is obtained as a basic Morita equivalence; for convenience we will
review the setup. This technique was successfully used by S.P. Smith
\cite{SmithSP.diffop,SmithSP:DonLine} B. Haastert
\cite{Haastert.DiffOp,Haastert.DirIm} and R. B{\o}gvad
\cite{Bog.DmodBorel} to study the ring of differential operators in finite
characteristic.

As so often when working with the Frobenius map $F^e:R \to R$, the fact
that the source and target are denoted by the same symbol easily leads
to confusion. Thus we temporarily generalize to an injective map of rings
$f: R \to A$ which makes $A$ a locally finitely generated free
$R$--module. Our application will be the Frobenius map which follows this
pattern if $R$ is $F$--finite and regular.
\begin{proposition}
    Let\/ $R$ be a ring and\/ $f: R \to A$ a map of rings such that\/ $A$ is a
    locally finitely generated free\/ $R$--module. Then\/ $f^* = A \tensor_R
    \usc$ is an equivalence between the category of\/ $R$--modules and\/
    $\End_R A$--modules. The inverse functor is\/ $\Hom_R(A,R)
    \tensor_{\End_R(A)} \usc[m]$.
\end{proposition}
\begin{proof}
We view $\Hom_R(A,R)$ as an $R$--$(\End_R A)$--bimodule. The action is by
post-- and pre--composition respectively. Similarly, we can view $A$
itself as a $(\End_R A)$--$R$--bimodule. Then it is clear that the
described associations are functors between the claimed categories as they
are just tensoring with an appropriate module.

It remains to show that they are canonically inverse to each other. For
this observe that the natural map
\[
    \Phi:\ A \tensor_R \Hom_R(A,R) \to \End_R A
\]
given by sending $a \tensor \phi$ to $a \cdot (f \circ \phi)$ is an
isomorphism by the assumption that $A$ is a locally finitely generated
\emph{free} $R$--module. With this in mind, checking that the above is an
isomorphism comes down to the fact that $\Hom\ $ commutes with finite
direct sums in the second argument. Thus $\Phi$ is a natural
transformation of $A \tensor_R \Hom_R(A,R) \tensor_{\End_R A} \usc$ to the
identity functor on $\End_R A$--mod.

Conversely it is equally easy to see that the map
\[
    \Psi:\ \Hom_R(A,R) \tensor_{\End_R(A)} A \to R
\]
given by sending $\phi \tensor a$ to $\phi(f(a))$ is also an isomorphism.
After a local splitting $\pi$ of $f$ is chosen (it exists by local
freeness of $A$ over $R$), its inverse is given by $a \mapsto \pi \tensor
a$. This finishes the proof of the proposition as $\Psi$ is now the
natural transformation from $\Hom_R(A,R) \tensor_{\End_R(A)} A \tensor_R
\usc$ to the identity functor on $R$--mod.
\end{proof}

We aim for a more explicit description of $\Hom_R(A,R)\tensor_{\End_R(A)}
M$ for an $\End_R(A)$--module $M$. Let $J$ denote the left ideal of
$\End_R(A)$ consisting of all $\phi \in \End_R(A)$ such that $\phi \circ
f$ is the zero map. In other words, $J$ is the kernel of the evaluation
map $\End_R(A) \to A$, sending $\phi$ to $\phi(1)$. We have a natural
isomorphism of $R$--modules
\begin{equation}\label{eqn.loc10}
    \Hom_R(A,R)\tensor_{\End_R(A)} M \to \Ann_{M} J
\end{equation}
given by sending $g \tensor m$ to $(f \circ g)(m)$. In fact, this map is
just the composition
\[
    \Hom_R(A,R) \tensor_{\End_R A} M \to[f \circ \tensor \id M] \End_R A
    \tensor_{\End_R A} M \to[\cong] M
\]
and we claim that its image is exactly $\Ann_{M} J$, which, by definition,
we identify with $\Hom_{\End_R A}(A, M)$ ($m$ corresponds to the map
$\phi_m$ sending $1 \in A$ to $m$). As in the proof of the last
proposition, any local splitting $\pi$ of $f$ enables us to write down an
inverse to \eqnref{eqn.loc10} explicitly. It is given by sending $m$ to
$\pi \tensor m$. For this to be inverse to \eqnref{eqn.loc10}
\begin{equation}\label{eqn.loc12}
    m \mapsto \pi \tensor m \mapsto (f \circ \pi)(m)
\end{equation}
must be the identity, \ie $(f \circ \pi)(m)=m$ for $m \in \Ann_{M} J$.
Indeed, we have
\[
    m = \phi_m(1) = \phi_m((f \circ \pi)(1)) = (f \circ \pi)\phi_m(1) = (f
    \circ \pi)(m)
\]
where we applied the  $\End_R A$--linearity of $\phi_m$ to $(f \circ
\pi)\in \End_R A$ to obtain the middle equality. Conversely, $\phi \tensor
m \mapsto (f \circ \phi) \tensor m \mapsto (\pi \circ f \circ \phi)
\tensor m = \phi \tensor m$ since $\pi \circ f = \id_R$. Thus we proved:
\begin{corollary}
    Let\/ $f: R \to A$ be a map of rings, such that\/ $A$ is a finitely
    generated locally free\/ $R$--algebra. Let\/ $J$ be the left ideal of\/
    $\End_R(A)$ consisting of all endomorphisms\/ $\phi$ such that\/
    $\phi(1)=0$. Then, for any\/ $\End_R(A)$--module\/ $M$ we have\/
    $\Hom_R(A,R) \tensor_{\End_R(A)}M \cong \Ann_M J$.

    Furthermore, if we have a splitting\/ $\pi$ of\/ $f$ (this exists locally,
    for example) then\/ $f \circ \pi \in \End_R(A)$ and\/ $Ann_M J = (f \circ
    \pi)(M)$. That is, the functors\/ $\Ann_{(\usc)}J$ and\/ $(f\circ
    \pi)(\usc)$ are inverse to the functor\/ $f^*(\usc)=A \tensor_R \usc$.
\end{corollary}
\begin{proof}
The first part of this proposition was proved in the preceding discussion.
For the second part, using the same reasoning as in \eqnref{eqn.loc12}, we
see that $\Ann_{M} J = (f \circ \pi)(M)$. The right is contained in the
left as we showed that $(f \circ \pi)(m)=m$ for $m \in \Ann_{M} J$.
Conversely, if $\phi \in J$ we have to show that $\phi((f \circ
\pi)(m))=0$. This is clear as $\phi((f \circ \pi)(m))=(\phi\circ f \circ
\pi)(m)$ and $\phi \circ f = 0$ by definition of $J$. Thus we get that $(f
\circ \pi)(M)=\Ann_M(J)$. We point out that the set $(f \circ \pi)(M)$ is
independent of the splitting $\pi$ since it is equal to $\Ann_M(J)$.
\end{proof}

Obviously this all applies to the Frobenius map $F^e: R \to F^e_*R$ if $R$
is regular and $F$--finite. In this case $F^e_*R$ (equivalently, $R^e$) is
a finitely generated and locally free left $R$--module (resp.\ right
$R$--module). In this case, $\End_R (F^e_*R) \cong \End_{R^{p^e}}(R) =
\D[e]_R$ and the ideal $J=J_e$ consists of the operators $\delta \in
\D[e]_R$ of level $e$ such that $\delta(1)=0$. Thus we get:
\begin{corollary}[Frobenius Descent]\label{cor.FrobDesc}
    Let\/ $R$ be regular and\/ $F$--finite. Then\/ $\F[e]$ is an equivalence of
    categories between the category of\/ $R$--modules and
    $\D[e]_R$--modules.

    If\/ $J$ denotes the ideal consisting of all\/ $\delta \in \D[e]_R$ with
    $\delta(1)=0$, then the inverse functor\/ $T^e$ of $\F[e]$ is given by\/
    $\Ann_{\usc} J$.

    If\/ $\pi_e$ is a splitting of\/ $F^e_R$, then the inverse functor\/ $T^e$
    can be realized as the projection operator\/ $(F^e_R \circ \pi_e) \in
    \D[e]_R$ which sends a\/ $\D[e]_R$--module\/ $M$ to its image under this
    differential operator.
\end{corollary}
\begin{remark}
This result should be compared to the so-called Cartier descent as
described, for example, in Katz \cite[Theorem 5.1]{Katz}. It states that
$F^*$ is an equivalence between the category of $R$--modules and the
category of modules with integrable connection and $p$--curvature zero.
The inverse functor of $\F$ on a module with connection $(M,\nabla)$ is in
this case given by taking the horizontal sections $\ker \nabla$ of $M$. As
an $R$--module with integrable connection and $p$--curvature zero is
nothing but a $\D[1]_R$--module, Cartier descent is just the case $e=1$ of
Corollary \ref{cor.FrobDesc}. Thus, in analogy with Cartier descent, we
can think of $T^eM$ as the module of \emph{horizontal sections} of level
$e$. Note that $T^e M=(F^e_R \circ \pi_e)(M)$ is an $R$--submodule of
$F^e_*M$ or equivalently an $R^{p^e}$--submodule of $M$.
\end{remark}

By transitivity, the last corollary implies that the categories of
$\D[e]_R$--modules for all $e$ are equivalent since each single one of
them is equivalent to $R$--mod. The functor giving the equivalence between
$\D[f]$--mod and $\D[f+e]$--mod is, of course, $\F[e]$. Concretely, to
understand the $\D[f+e]$--module structure on $\F[e]M$ for some
$\D[f]$--module $M$, we write $M \cong \F[f]N$ for $N=T^f(M)$. Then
$\F[e]M=\F[(f+e)]N=R^{(f+e)}\tensor N$ carries obviously a
$\D[f+e]$--module structure with $\delta \in \D[f+e]$ acting via $\delta
\tensor \id_N$.

Since $D_R =\bigcup \D[e]_R$ this shows that $\F[e]$ is in fact an
auto--equivalence on the category of $D_R$--modules:
\begin{proposition}\label{prop.FrobDesc4Dmod}
    Let\/ $R$ be regular and\/ $F$--finite. Then\/ $\F[e]$ is an equivalence of
    the category of\/ $D_R$--modules with itself. The inverse is given by\/
    $T^e \defeq \Ann_{(\usc)}J$ where\/ $J$ consists of all operators\/
    $\delta \in \D[e]_R$ such that\/ $\delta(1)=0$. If we have a splitting\/
    $\pi^e$ of the Frobenius, then\/ $T^e=\Image ( F^e \circ \pi_e )$ and
    this is independent of the choice of such splitting.
\end{proposition}
\begin{proof}
In the last paragraph we showed how $\F[e]$ is an equivalence between
$\D[f]$--mod and $\D[f+e]$--mod for all $f$. By viewing the $D_R$--module
$M$ as a $\D[f']$--module for all $f'$, we observe that $\F[e]M$ is a
$\D[f]$--module for all $f=f'+e$. These structures of various levels are
compatible with the inclusion $\D[f] \subseteq \D[g]$ for $g \geq f$ and
thus we get a $D_R$--structure on $\F[e]M$ naturally. Checking this
claimed compatibility comes down to remembering that $\F[e] \circ T^e$ is
isomorphic to the identity functor. Then the argument is very similar to
the proof of Proposition \ref{prop.wellDefDRonUnit}. It comes down to
observing that for $\delta \in \D[f]$ the following diagram commutes.
\[
\xymatrix{
    {\F[(f+1)](T^{(f+1)-e}(M))}\ar@{=}[r]\ar_{\delta\tensor\id_{T^{(f+1)-e}(M)}}[d] &{\F[f](T^{f-e}(M))}
                                \ar@{=}[r]\ar_{\delta\tensor \id_{T^{f-e}(M)}}[d] &{M} \ar@{-->}^{\delta}[d]\\
    {\F[(f+1)](T^{(f+1)-e}(M))} \ar@{=}[r]       &{\F[f](T^{f-e}(M))}
                                \ar@{=}[r]       &{M}
    }
\]
The equal signs indicate the natural isomorphism of functors $\F[e] \circ
T^e \cong \id$, possibly applied repeatedly. The dashed arrow represents
the action of $\delta$, and it is well defined if and only the left square
commutes. To see this we observe once more that $\delta \in \D[f]$, under
the inclusion of $\D[f]_R \subseteq \D[f+1]_R$, is represented by $\delta
\tensor \id \in \End_{\modR}(R^f \tensor R^1)$. The second part now
follows from the construction of $T^e$ as in Corollary \ref{cor.FrobDesc}.
\end{proof}

\subsection{Frobenius Descent for unit $R[F]$--modules}
In this section we develop a concrete version of the Frobenius descent for
$D_R$--submodules of a unit $R[F^e]$--module. Just as the Frobenius
functor $\F[e]N$ of some submodule $N$ of a unit $R[F^e]$--module
$(M,\theta^e)$ has its concrete description as the submodule
$F^e(M)=\theta^e(\F[e](N))$ of $M$, there is an analog description of its
inverse $T^e$ for $D_R$--submodules of $M$.
\begin{proposition}
    Let\/ $R$ be regular and\/ $F$--finite and let\/ $(M,\theta^e)$ be a unit\/
    $R[F^e]$--module. If\/ $N$ is a\/ $D^e_R$--submodule, then, if\/ $\pi_e$ is
    a splitting of the Frobenius\/ $F^e_R$ (exists locally), we have\/ $T^e(N)
    \cong (\pi_e \tensor \id_M) \circ (\theta^e)^{-1}(N)$. The resulting left
    $R$--linear injection $T^e(N) \subseteq M$ is independent of $\pi_e$.
\end{proposition}
\begin{proof}
We use the description of $T^e$ as in \ref{cor.FrobDesc} which was also
independent of the chosen splitting $\pi_e$. The proof comes down to
observing that the operation of the differential operator $F^e_R \circ
\pi_e$ on a unit $R[F^e]$--module $(M,\theta^e)$ is by definition given by
the lower path of this diagram.
\[
\xymatrix@C=4pc{
    {M}\ar^{(F^e_R \circ \pi_e)}[rr]\ar_{{\theta^e}^{-1}}[d]&&{M} \\
    {\F[e]M}\ar^(.55){\pi_e \tensor \id_M}[r]&{M}\ar^(.45){F^e_R \tensor \id_M}[r]&{\F[e]M}\ar_{\theta^e}[u] }
\]
Now just note that $\theta^e \circ (F^e_R \tensor \id_M)=F^e_M$ is a
$p^e$--linear isomorphism from $(\pi_e \tensor \id_M) \circ
(\theta^e)^{-1}(N)$ to $T^e(N)=(F^e_R \circ \pi_e)(N)$. This establishes
$T^e(N)$ naturally as an $R$--submodule of $M$.
\end{proof}
It is important to note the difference from the description of $T^e$ in
Corollary \ref{cor.FrobDesc}. There $T^e(N)=(F^e \circ \pi_e)(N)$ is an
$R^{p^e}$--submodule of $M$. Here, on the other hand, we show that
$T^e(N)$ is isomorphic to the $R$--submodule $(\pi_e \tensor \id_M) \circ
(\theta^e)^{-1}(N)$ of $M$. In fact it is easy to see that if $N$ is a
$D_R$--submodule of $M$, then $T^e(N)$ in this description is also a
$D_R$--submodule of $M$. \note{Include explanation}

We indicate yet another description of $T^e(N)$ in this case. Let
$(M,\theta^e,F^e)$ be a unit $R[F^e]$--module. We claim that for any
$D_R$--submodule $N$ of $M$ we have
\[
    T^e(N) \cong (\pi_e \tensor \id_M) ({\theta^e}^{-1}(N)) = (F^e)^{-1}(N).
\]
To see this, first note that $(F^e)^{-1}(N) = (F^e_R \tensor
\id_M)^{-1}({\theta^e}^{-1}(N))$, since $\theta^e$ and $F^e$ correspond to
each other by adjointness of $\F[e]$ and $F^e_*$. Setting $N' =
{\theta^e}^{-1}(N)$, it remains to see that $(F^e_R \tensor
\id_M)^{-1}(N') = (\pi_e \tensor \id_M) (N')$. Since $\pi_e \circ F^e_R =
\id_R$ we have $(\pi_e \tensor \id_M) (N') \subseteq (F^e_R \tensor
\id_M)^{-1}(N')$. The reverse inclusion follows from the fact that the
differential operator $F^e_R \circ \pi_e$ fixes $N'$ (as it is a
$D_R$--submodule of $\F[e]M$) and therefore $(F^e_R \tensor
\id_M)^{-1}(N')\subseteq  (F^e_R \tensor \id_M)^{-1}((F^e_R \circ
\pi_e)(N')) = (\pi_e \tensor \id_M)(N')$ which finishes the argument. The
last equality follows from the fact that the map $F^e_R$ is pure and thus
$F^e_R \tensor \id_M$ is injective.
\begin{corollary}
    Let\/ $M$ be a unit\/ $R[F^e]$--module. Then, for every\/ $D_R$--submodule\/
    $N$ we have
    \[
        F^e((F^e)^{-1}(N)) = N\ \text{ and } \ (F^e)^{-1}(F^e(N))=N.
    \]
    That is,\/ $T^e=(F^e)^{-1}$ as functors on\/ $D_R$ submodules of\/ $M$.
\end{corollary}
\begin{remark}
The version of Frobenius descent we encounter in Lyubeznik uses that a
$D_R$--submodule $N$ of a finitely generated unit $R[F^e]$--module $M$ is
described by a certain infinite collection $\{N_i\}$ of $R$--submodules
$N_i = F^{ei}(N) \cap M_0$ of a root $M_0$ of $M$. With this
characterization, the $D_R$--submodule of $M$ corresponding to $T^e(N)$ is
represented by the infinite collection $\{N'_i\}$ where $N'_i=N_{i+1}$
(see \cite[Proposition 5.4]{Lyub}). This point of view is related to the
notion of an $X^{\infty}$--module as in \cite{Haastert.DiffOp} or
\cite{SmithSP.diffop}.
\end{remark}
\subsection{$D_R[F^e]$--modules}
In Section \ref{sec.UnitRFareD} we describe the natural $D_R$--module
structure on a unit $R[F^e]$--module $M$. A consequence of Proposition
\ref{prop.FrobDesc4Dmod} is that $\F[e]M$ carries a natural
$D_R$--structure as well. And in fact with this $D_R$--structure on $M$
and $\F[e]M$ respectively, we have:

\begin{lemma}\label{lem.ThetaDlinear}
    Let\/ $R$ be regular and\/ $F$--finite and\/ $(M,\theta^e)$ be a unit\/
    $R[F^e]$--module. Then\/ $\theta^e$ is a map of\/ $D_R$--modules.
\end{lemma}
\begin{proof}
Just as in the proof of Proposition \ref{prop.wellDefDRonUnit} this comes
down to observing the commutativity of the left hand square of the diagram
\eqnref{eqn.WellDefDRonUnit} in that same proof.
\end{proof}

This motivates us to define a $D_R[F^e]$--module in analogy with
$R[F^e]$--modules as a $D_R$--module $M$ together with a $D_R$--linear map
$\theta^e: \F[e]M \to M$.

\begin{definition}
    A \emph{$D_R[F^e]$--module} \index{DRFmodule@$D_R[F^e]$--module} is a
    $D_R$-module $M$ together with a $D_R$-linear map
    \[
        \theta^e_M: \F[e]_R M \to{} M.
    \]
    In other words a $D_R[F^e]$--module is an $R[F^e]$--module
    $(M,\theta^e)$ that carries a $D_R$--structure such that $\theta^e$ is
    $D_R$--linear.
\end{definition}
As just remarked, unit $R[F^e]$--modules are $D_R[F^e]$--modules.
Conversely, by forgetting the $D_R$--structure, every unit
$D_R[F^e]$--module is a unit $R[F^e]$--module. Thus, for the unit case
this doesn't lead to anything new. One reason for working with the
seemingly more complicated category of $D_R[F^e]$--modules is the
following proposition.
\begin{proposition}\label{prop.DRfulluRF}
    Let $R$ be regular and $F$--finite. The category of finitely generated
    unit\/ $D_R[F^e]$--modules is a \emph{full} subcategory of the
    category of\/ $D_R[F^e]$--modules.
\end{proposition}
Note that the category of finitely generated unit $R[F^e]$--modules is not
full in $R[F^e]$--mod. The example of an ideal $I$ of $R$ which is an
$R[F^e]$--submodule but not unit shows this nicely (cf.\ Example
\ref{ex.BasicRFmods}).

\begin{proof}[Proof of \ref{prop.DRfulluRF}]
We have to show that a $D_R[F^e]$--submodule of a finitely generated unit
$R[F^e]$--module is in fact also unit. If $N$ is such a $F^e$--stable
$D_R$--submodule of a unit $R[F^e]$--module $M$, then $F^e(N) \subseteq
N$. Applying $T^e$ and using its defining property we see that $N
\subseteq T^e(N)$. Iterating we get an increasing chain of $R$--modules
\[
        N \subseteq T^e(N) \subseteq T^{2e}(N) \subseteq \ldots.
\]
Intersecting this chain with a root $M_0$ of $M$ yields a chain of
submodules of $M_0$ which, as $M_0$ is a finitely generated $R$--module,
must stabilize. Let $N_0 \defeq T^{er}(N) \cap M_0=T^{e(r+1)}(N)\cap M_0=
\ldots$ be the stable member. Using $F^{er}(T^{er}(N))=N$, we see that
$F^e(N_0)=F^e(T^{e(r+1)}(N_0) \cap M_0)=T^{er}(N_0)\cap F^e(M_0) \supseteq
T^{er}(N_0) \cap M_0 = N_0$. Applying $F^{er}$ for all $r$ to this
inclusion we get another increasing sequence
\[
    N_0 \subseteq F^e(N_0) \subseteq F^{2e}(N_0) \subseteq \ldots
\]
Let $L$ be its limit. As $L$ arises as the increasing union of the
Frobenius powers of a single submodule, it is obviously a unit submodule
of $N$, \ie $F^e(L)=L \subseteq N$. For the converse inclusion let $n \in
N$. For all sufficiently large $r \geq 0$ we have $n \in F^{er}(M_0)$. For
such $r$ also $N_0=T^{er}(N) \cap M_0$ and thus $F^{er}(N_0)=N \cap
F^{er}(M_0)$. Thus $n \in F^{er}(N_0)$, and therefore $n \in L$ as $L$ is
the increasing union of all $F^{er}(N_0)$.
\end{proof}
\begin{remark}
    Analogously to the case of $R[F^e]$--modules one can show that
    $D_R[F^e]$--modules are just modules over an appropriate ring
    $D_R[F^e]$. In the case that the ring $R[F^e]$ is in fact a subring of
    $\End_k(R)$ (cf.\ Remark \ref{rem.RFsubEnd}) we can think of this ring
    $D_R[F^e]$ as the subring of $\End_k(R)$ generated by $R[F^e]$ and
    $D_R$. In general one can define the ring $D_R[F]$ to be $R[F^e]
    \tensor_R D_R$ and then equip this tensor product with an appropriate
    ring structure. This is done in \cite{Em.Kis2} where many other
    interesting properties in this context are shown. For example, they
    show the following (at first) surprising result:
    \begin{proposition}
        Let $R$ be regular essentially of finite type over a perfect field
        $k$. Let $M$ be a unit $R[F^e]$--module. The following are
        equivalent.
        \begin{enumerate}
        \item $M$ is finitely generated as a $D_R[F^e]$--module.
        \item $M$ is finitely generated as an $R[F^e]$--module.
        \item $M$ is finitely generated as a $D_R$--module.
        \end{enumerate}
    \end{proposition}
    The proof of this is a clever application of Frobenius descent
    together with Lyubeznik's Theorem 5.6 \cite{Lyub}.
\end{remark}

\subsection{Simple finitely generated unit $R[F^e]$--modules are semisimple $D_R$--modules}
Another main result of \cite{Lyub} is his Theorem 5.6 which shows that
finitely generated unit $R[F^e]$--modules have finite length as
$D_R$--modules for $R$ an $F$--finite and regular ring. Before we state
the result we first recall the definition of semisimplicity.
\begin{definition}
A $D_R$--module $M$ is called semisimple if every submodule is a direct
summand.
\end{definition}
There are other characterizations of semisimplicity. If $M$ has finite
length as a $D_R$--module, then the following are equivalent:
\begin{enumerate}
    \item $M$ is semisimple.
    \item Every simple submodule of $M$ is a direct summand.
    \item $M \cong N_1 \oplus \ldots \oplus N_t$ where the $N_i$ are simple $D_R$--modules.
\end{enumerate}
The third characterization is the one we will make use of in our
applications. Another rather trivial but important observation is that
every submodule of a semisimple module is itself semisimple.
\begin{theorem}\label{thm.unitRFsemisimpleDR}
    Let\/ $R$ be\/ $F$--finite and regular and\/ $M$ be a finitely generated
    simple unit\/ $R[F^e]$--module. Then, as a\/ $D_R$--module,\/ $M$ is
    semisimple and of finite length.
\end{theorem}
\begin{proof}
    By \cite[Theorem 5.6]{Lyub}, $M$ has finite length as a $D_R$--module
    and thus it remains to check that every simple $D_R$--submodule is a
    direct summand. Let $N$ be a simple $D_R$--submodule of $M$. By
    Frobenius descent, $F^{er}$ is an equivalence of categories; in
    particular $F^{er}(N)$ is also a simple $D_R$--module.

    The simplicity implies that $(N+F^e(N)+\cdots+F^{e(r-1)}(N))\cap
    F^{er}(N)$ is either zero or all of $F^{er}(N)$. Let $r$ (allow
    $r=\infty$) be the first time the second case happens and set
    $N'\defeq N+F^e(N)+\cdots+F^{e(r-1)}(N)$ noticing that the sum is in
    fact direct. Then $F^e(N')=F^e(N)+\cdots+F^{er}(N) \subseteq N'$ as
    $F^{er}(N) \subseteq N'$. Thus $N'$ is a $D_R[F^e]$--submodule of the
    finitely generated unit $R[F^e]$--module $M$ and therefore also
    finitely generated and unit by Proposition \ref{prop.DRfulluRF}. As
    $M$ is a simple $R[F^e]$--module, $N'=M$.

    The case $r=\infty$ contradicts the finite length of $M$ as a
    $D_R$--module by \cite[Theorem 5.6]{Lyub}. To illustrate some
    techniques we include an alternative argument. Assuming $r=\infty$
    will contradict the finite generation of $M$ as an $R[F^e]$--module:
    Let $M_0$ be a root of $M$, in particular $M_0 \subseteq F^{et}(M_0)$
    for all $t$. As $M'$ is finitely generated we have $M_0 \subseteq
    (N+F^e(N)+\cdots+F^{e(s-1)}(N))
    \defeq N_0$ for some $s$ and $F^{es}(N_0) \cap N_0 = 0$ by assumption
    of $r=\infty$. But this implies that $M_0 \subseteq F^{es}(M_0)
    \subseteq F^{es}(N_0)$ and thus $M_0 \subseteq F^{es}(N_0) \cap N_0=0$
    which is a contradiction.

    Thus we showed that $M \cong N \oplus (F^e(N) \oplus \ldots \oplus
    F^{e(r-1)}(N))$ and therefore $M$ is semisimple and of finite length
    $r$ as a $D_R$--module.
\end{proof}
This proof would be independent of \cite[Theorem 5.6]{Lyub} if we could
find an alternative argument for the fact that $M$ has the descending
chain condition as a $D_R$--module, or at least if we could ensure the
existence of some nonzero simple $D_R$--submodule of $M$.

\subsection{Semisimple $D_R$--module with Frobenius action}
The $D_R$--semi\-sim\-pli\-ci\-ty of simple unit $R[F^e]$--modules has
some simple but interesting consequences.

We have a canonical decomposition $M \cong N_1^{\oplus n_1} \oplus \ldots
\oplus N_t^{\oplus n_t}$ where $N_i$ are \emph{distinct} simple
$D_R$--modules. $M_i=N_i^{\oplus n_i}$ is called the $N_i$--isotypic
component of $M$. This decomposition of $M$ into its isotypic components
is unique, contrary to the decomposition into simple modules. Modules
which have only one isotypic component, \ie which are isomorphic to a
direct sum of copies of the same simple $D_R$--module $N$ are called
$N$--isotypic, or just isotypic. We first observe the following lemma.

\begin{lemma}
    Let\/ $R$ be regular. Let\/ $M \cong M_1 \oplus \ldots \oplus M_t$ be the
    decomposition of the\/ $D_R$--module\/ $M$ into its isotypic components. Then\/
    $\F[e]M \cong \F[e]M_1 \oplus \ldots \oplus \F[e]M_r$ and the\/ $\F[e]M_i$
    are the isotypic components of\/ $\F[e]M$.
\end{lemma}
\begin{proof}
    As the Frobenius commutes with finite direct sums the first statement
    is clear. If $M_i \cong N_i^{\oplus n_i}$ for a simple $D_R$--module
    $N$, then $\F[e]M_i \cong (\F[e]N_i)^{\oplus n_i}$. By Frobenius
    descent ($\F[e]$ is an equivalence of categories, Corollary
    \ref{prop.FrobDesc4Dmod}) $\F[e]N_i$ is also a simple $D_R$--module
    and for $i \neq j$ we have $\F[e]N_i \not\cong \F[e]N_j$. This proves
    all assertions.
\end{proof}
As a consequence of this lemma we get a simple but important proposition
which roughly says that $D_R$--isotypic components of unit
$D_R[F^e]$--modules are indeed $D_R[F^{er}]$--submodules for some $r>0$.
\begin{proposition}\label{prop.IsotypicareRFsub}
    Let\/ $M \cong M_1 \oplus \ldots \oplus M_t$ be a\/ $D_R$--module,
    where\/ $M_i$ are the finitely many isotypic components. Let\/
    $\theta^e: \F[e]M \to M$ be an isomorphism of\/ $D_R$--modules
    (\/$(M,\theta^e)$ is a unit\/ $D_R[F^e]$--module). Then there is\/ $r
    > 0$ such that for all\/ $i$, the restriction of\/ $\theta^{er}$ to
    $\F[er]M_i$ is an isomorphism onto\/ $M_i$.
\end{proposition}
\begin{proof}
    Note that $\theta^e$ maps isotypic components of $\F[e]M$ to isotypic
    components of $M$. As $\theta^e$ is an isomorphism it maps $\F[e]M_i$
    isomorphically onto $M_{\sigma(i)}$ for some permutation $\sigma \in
    S_t$, the group of permutations of $\set{1,\ldots,t}$. By induction on
    $r$ it follows that $\theta^{er}$ maps $\F[er]M_i$ isomorphically onto
    $M_{\sigma^r(i)}$. With $r=\op{order}(\sigma)$ the proposition
    follows.
\end{proof}
This proposition in particular shows that the direct sum decomposition of
$M$ into its $D_R$--isotypic components is in fact a decomposition of $M$
as an $R[F^{er}]$--module for some $r$. Next we focus our attention on
$D_R$--modules with just one isotypic component.
\begin{lemma}\label{lem.SubofIsotypareUnit}
    Let\/ $M$ be\/ $N$--isotypic for a simple\/ $D_R$--module\/ $N$. Let
    $\theta^e: \F[e]M \to M$ be an isomorphism of\/ $D_R$--modules. Then
    as\/ $D_R$--modules,\/ $\F[e]N$ is (non canonically) isomorphic to\/
    $N$.
\end{lemma}
\begin{proof} $M$ being $N$--isotypic implies that $\F[e]M$ is $\F[e]N$--isotypic.
But $\F[e]M \cong M$ implies that $M$ is $N$--isotypic. Therefore $\F[e]N \cong
N$.
\end{proof}
As a corollary of all this we get a rather general statement about
$D_R$--simple submodules of finitely generated unit $R[F^e]$--modules:
\begin{proposition}\label{prop.SimpSubofUnitisUnit}
    Let\/ $R$ be regular. Let\/ $M$ be a finitely generated unit\/
    $R[F^e]$--module and let\/ $N \subseteq M$ be a simple\/
    $D_R$--submodule. Then\/ $N$ is a simple unit\/ $R[F^{er}]$--module
    for some\/ $r>0$.
\end{proposition}
\begin{proof} By Theorem \ref{thm.fguRFhaveDCC} $M$ has finite length as a unit
    $R[F^e]$--module. Thus $M$ has a finite filtration by finitely
    generated unit $R[F^e]$--modules whose quotients are simple unit
    $R[F^e]$--modules. Since $N$ is simple as a $D_R$--module it is a
    submodule of one of the quotients of this filtration. Replace $M$ by
    this quotient and assume therefore that $M$ is a simple, finitely
    generated unit $R[F^e]$--module. By Theorem
    \ref{thm.unitRFsemisimpleDR} $M$ is semisimple of finite length as a
    $D_R$--module and by Proposition \ref{prop.IsotypicareRFsub} each
    isotypic component of $M$ is a unit $R[F^{re}]$--module for some
    common $r>0$. Thus, again by $D_R$--simplicity of $N$ we see that $N$
    is in fact a submodule of one of the isotypic components. Replace $M$
    by this isotypic component and apply Lemma
    \ref{lem.SubofIsotypareUnit} to obtain that $N$ is a unit
    $R[F^{re}]$--module. Note that $N$ is then automatically simple as a
    unit $R[F^{er}]$--module since it is $D_R$--simple.
\end{proof}
A similar proof shows the equivalent proposition for $D_R$--simple
quotients of finitely generated unit $R[F^e]$--modules. It is important to
emphasize that with this $R[F^{re}]$--module structure, $N$ normally is
not an $R[F^{re}]$--submodule of $M$. That is, the inclusion $N \subseteq
M$ is normally \emph{not} a map of $R[F^{re}]$--modules. The reason for
this that there are cases where no power of the Frobenius on $M$
stabilizes the $D_R$--submodule $N$, thus $N$ cannot be an
$R[F^{er}]$--submodule of $N$. See Section \ref{sec.ExAndCounterEx} below
for such an example.

\section{Endomorphisms of simple $D_R$--modules}
The proof of Main Theorem 1, which will follow shortly, relies on a
seemingly natural expected property of the endomorphism ring of a simple
$D_R$--module.
\begin{claim}\label{cl.EndDsimp}
    Let\/ $R$ be a regular ring, essentially of finite type over a field\/
    $k$. If\/ $N$ is a simple\/ $D_R$--module, then\/ $\End_{D_R}(N)$ is
    algebraic over\/ $k$.
\end{claim}
In characteristic zero this follows from the fact that $D_R$ is almost
commutative, \ie it has a commutative associated graded which is a
finitely generated (commutative) algebra over $k$. Then a form of
Quillen's Lemma \cite{Quill} implies that, in fact, endomorphisms of
simple $D_R$--modules satisfy an algebraic equation. In this proof the
finitely generatedness of $D_R$ over $R$ is critically used in the form of
a noncommutative generic freeness argument. As already mentioned, in
finite characteristic, $D_R$ is not finitely generated. Thus the
characteristic zero argument miserably fails. Nevertheless, we have the
following result hinting that Claim \ref{cl.EndDsimp} might still be true
in finite characteristic. It goes back to Dixmier \cite{Dix} and we recall
the argument for convenience.
\begin{lemma}\label{lem.EndNisAlgebraic}
    Let\/ $R$ be essentially of finite type over an \emph{uncountable} field\/
    $k$. If\/ $N$ is a simple\/ $D_R$--module, then\/ $\End_{D_R}(N)$ is
    algebraic over\/ $k$.
\end{lemma}
\begin{proof}
As $N$ is simple, we have $D_Rn=N$ for some (every nonzero) $n \in N$.
Thus every $\phi \in \End_{D_R}(N)$ is determined by its value on $n$. As
$D_R$ is countably dimensional over $k$ so is $N=D_Rn$ and thus
$\End_{D_R}(N)$ is also countably dimensional over $k$. Therefore, for any
fixed $\phi \in \End_{D_R}(N)$ (say $\phi \not\in k$) the uncountable set
$\set{(\phi + \lambda)^{-1}\ |\  \lambda \in k\ }$ must be linearly
dependent (we use that $\End_{D_R}(N)$ is a division ring by Shur's
lemma). A relation of linear dependence among some finitely many $(\phi +
\lambda_i)^{-1}$ gives, after clearing denominators, an algebraic relation
for $\phi$. Clearing denominators works just like in the commutative case
since all $(\phi + \lambda_i)^{-1}$ commute with each other.
\end{proof}
\begin{remark}\label{rem.EndDNgeneral}
The same proof works whenever the cardinality of $k$ is strictly
bigger then the cardinality of a $k$--vectorspace basis of $R$. Thus we get
that $\End_{D_R}(N)=k$ whenever $k$ is an algebraically closed field
contained in $R$ such that the cardinality of $k$ is strictly bigger then
the cardinality of a $k$--vectorspace basis of $R$. In the above lemma the
latter is countable and the former is uncountable.
\end{remark}
The proof is characteristic independent and even works for non--regular
$R$. Therefore it gives little indication whether the statement is true
without the awkward assumption of uncountability of $k$; \ie whether Claim
\ref{cl.EndDsimp} holds in finite characteristic. Nevertheless we have a
strong, though somewhat unfounded, belief that this is in fact the case.

In the proof of Main Theorem 1 the fact that $\End_{D_R}(N)$ is algebraic
over $k$ is an essential ingredient. More precisely, we use that if $k$ is
algebraically closed, then $\End_{D_R}(N)=k$ whenever $N$ is a simple
$D_R$--module. Thus we are bound to the assumption of $k$ being
uncountable; if Claim \ref{cl.EndDsimp} would be established this
assumption would no longer be needed for Main Theorem 1.

\section{Change of $R[F]$--module structure}
This section investigates the question on how different $R[F^e]$--module
structures on the same $R$--module are related to each other. For a given
$R$--module $M$ we want to describe the set of all $R[F^e]$--structures on
$M$. By definition, this is just the set $\Hom_R(\F[e]M,M)$. We mainly
focus on the case that $M$ does indeed have some unit $R[F^e]$--structure
$\theta^e: \F[e]M \to[\cong] M$. In this case, pre-composing with
${\theta^e}^{-1}$ induces an isomorphism of $\Hom_R(\F[e]M,M) \cong
\End_R(M)$, where the unit $R[F^e]$--structures of the left hand side
correspond to $\Aut_R(M)$, the automorphisms of $M$.

As an example for the case where $M$ is not necessarily unit we discuss in
Section \ref{sec.InjHullEx} the injective hull of the residue field of a
local ring $R$. This leads to a characterization of the Gorenstein
property of $R$ in terms of the unit structures on $E_{R/m}$.\note{still
not proven maybe not true} We begin with some basic observations.
\begin{lemma}\label{lem.thetaPhi}
    Let\/ $(M,\theta^e)$ be a unit\/ $R[F^e]$--module. The\/ $R[F^e]$--module
    structures\/ $\Hom_R(\F[e]M,M)$ of\/ $M$ are in one-to-one correspondence
    with\/ $\End_R(M)$. If\/ $R$ is\/ $F$--finite, then\/ $\End_R(M)$ carries a
    natural unit\/ $R[F^e]$--module structure.
\end{lemma}
\begin{proof}
Pre--composing with the map ${\theta^e}^{-1}$ gives an isomorphism
$\Hom_R(\F[e]M,M) \cong \End_R(M)$. The natural unit $R[F^e]$--structure
on $\End_R(M)$ is given by
\[
    R^e \tensor \End_R(M) \to \End_R(R^e \tensor M) \to \End_R(M)
\]
where the first map is the natural one which is an isomorphism if $R$ is
$F$--finite (cf.\ Proposition \ref{prop.HomAndTensor}). The second map is
obtained by pre--composing with ${\theta^e}^{-1}$ and post--composing with
$\theta^e$. Therefore it is an isomorphism since $\theta^e$ is an
isomorphism.
\end{proof}
The situation can be illustrated nicely with the following
diagram. We denote the Frobenius action corresponding to
$\theta^e$ by $F^e$. Let $({\theta'}^e,{F'}^e)$ be another
$R[F^e]$--structure on $M$. Then the diagram
\[
  \xymatrix{  & {\F[e]M} \ar[d]_{\theta^e} \ar[rd]^{{\theta'}^e} & {} \\
             {M} \ar ^{F^e_R \tensor \id} [ru] \ar^{F^e}[r] \ar ^{{F'}^e}
            @/_2pc/ [rr] & {M} \ar[r]^{\phi} & {M} }
\]
is commutative, and $\phi$ is the endomorphism of $M$ which
corresponds to ${\theta'}^e$ under the identification
$\Hom_R(\F[e]M,M) \cong \End_R(M)$. That is, we have ${\theta'}^e
= \phi \circ \theta^e$ and ${F'}^e=\phi \circ F^e$. In this
context we say that $\phi$ represents ${\theta'}^e$.

If we denote, somewhat abusively, the Frobenius action on $\End_R(M)$
induced by $\theta^e$ also by $F^e$, then we have that
\begin{equation*}\tag{\lltag}\label{eqn.Fphi}
    F^e(\phi)=\theta^e\circ \F[e]\phi \circ {\theta^e}^{-1}.
\end{equation*}
This description makes it clear that $\phi \in \End_R(M)$ is a map of
$R[F^e]$--modules $(M,\theta^e)$ if and only if $\phi = F^e(\phi)$, \ie
$\End_{R[F^e]}(M)$ are those $\phi \in \End_R(M)$ that are fixed by this
Frobenius action.

Another consequence is that ${\theta'}^e$ is an isomorphism if and only if
$\phi$ is an isomorphism by the formula ${\theta'}^e = \phi \circ \theta^e$.
Thus we have:
\begin{lemma}
    Let\/ $(M,\theta^e)$ be a unit\/ $R[F^e]$--module. The unit\/
    $R[F^e]$--module structures on\/ $M$ are in one-to-one correspondence
    with\/ $\Aut_R(M)$.
\end{lemma}
Before turning to some examples we note a lemma explaining how our notation
behaves under higher powers of the Frobenius.
\begin{lemma}\label{lem.highertheta}
    Let\/ $(M,\theta^e,F^e)$ be a unit\/ $R[F^e]$--module. Let\/ ${\theta'}^{e}=\phi
    \circ \theta^e$ be another\/ $R[F^e]$--structure on the\/ $R$--module\/ $M$. Then
    \[
        {\theta'}^{er} = \phi_r \circ \theta^{er}   \text{\quad and \quad}
        {F'}^{er} = \phi_r \circ F^{er}.
    \]
    where\/ ${F'}^e$ is the Frobenius action corresponding to\/ ${\theta'}^e$ and\/
    $\phi_r$ is inductively defined by\/ $\phi_{r+1} = \phi_r \circ
    F^{er}(\phi)$.
\end{lemma}
\begin{proof}
This follows by induction on $r$ remembering that $\theta^{er}$
and ${\theta'}^{er}$ were also defined inductively (cf.\ page
\pageref{pg.deftheta}) as $\theta^{e(r+1)}=\theta^{er} \circ
\F[er]\theta^e$, and analogously for ${\theta'}^e$. The base case
${\theta'}^e = \phi \circ \theta^e$ is by assumption. The
inductive step reads:
\[
\begin{split}
    {\theta'}^{e(r+1)} &= {\theta'}^{er} \circ \F[er]{\theta'}^e = \phi_r \circ
    \theta^{er} \circ \F[er] \phi \circ \F[er] \theta^e = \phi_r \circ
    F^{er}(\phi)
    \circ \theta^{er} \circ \F[er]\theta^e \\
    &= \phi^{r+1} \circ \theta^{e(r+1)}
\end{split}
\]
Here we used equation \eqnref{eqn.Fphi} above with $e$ substituted by $er$
to replace $\theta^{er} \circ \F[er] \phi$ by $F^{er}(\phi) \circ
\theta^{er}$.
\end{proof}

To get a feeling for what is behind all these notationally complicated
but, in fact, rather trivial abstract observations we look at the case we
are most interested in, that of free $R$--modules.

\subsection{Free $R$--modules}\label{sec.FreeRF} We consider the case that $M$ is a finitely
generated free $R$--module. We describe a formalism to represent Frobenius
actions (\ie $p^e$--linear maps) on $M$ by matrices analogous to how
linear maps are represented by matrices. This is done with the help of
Lemma \ref{lem.thetaPhi} which allows us to represent any Frobenius
structure $\theta^e$ by a linear map $\phi \in \End_R(M)$ after fixing
some unit $R[F^e]$--structure on $M$. The choice of a free basis of $M$
fixes such a unit $R[F^e]$--structure as well as it allows to represent
the map $\phi$ by a matrix $A$. This is, of course, the matrix we have in
mind to represent $\theta^e$.

Concretely, if $M \cong R^{\oplus n}$ is identified with the space of
columnvectors of length $n$ and entries in $R$, then a Frobenius action
$F^e$ on $M$ is determined by its values on the basis $e_i =
(0,\ldots,1,\ldots,0)^t$. With this notation, if $F^e(e_j)=\sum a_{ij}e_i$
and if we denote the matrix with entries $(a_{ij})$ by $A$, then
\[
    F^e(\begin{pmatrix} r_1 \\ r_2 \\ \vdots \\ r_n \end{pmatrix}) =
    A\begin{pmatrix} r^{p^e}_1 \\ r^{p^e}_2 \\ \vdots \\ r^{p^e}_n
    \end{pmatrix}.
\]

What we investigate here is how this representing matrix $A$ behaves under
change of basis and under taking higher powers of the Frobenius. In this
context, the matrix $A_r \defeq AA^{[p^e]} \cdots A^{[p^{ei}]} \cdots
A^{[p^{er}]}$ appears naturally, where $A^{[p^{ei}]}$ is the matrix with
entries the $p^{ei}$TtH powers of the entries of $A$.

The following proposition summarizes the results. It is the crucial
ingredient in the construction of the examples in Section
\ref{sec.ExAndCounterEx} below.

\begin{proposition}\label{prop.bcForRF}
    Let\/ $R$ be a ring. Let\/ $M$ be a finitely generated free\/ $R$--module of
    rank\/ $n$. After fixing a basis\/ $(\xn[e])$ of\/ $M$, the\/ $R[F^e]$--module
    structures of\/ $M$ are represented by square matrices of size\/ $n$.

    If\/ $A$ denotes the matrix representing some Frobenius structure\/ $F^e$,
    then\/ $A_r =  AA^{[p^e]} \cdots A^{[p^{er}]}$ is the matrix
    representing\/ $F^{er}$.

    If\/ $(\xn[f])$ is a new basis of\/ $M$ and\/ $C$ is the base change matrix
    (\ie $f_j=\sum c_{ij}e_i$), then with respect to this new basis\/ $F^e$
    is represented by\/ $B \defeq C^{-1}AC^{[p^e]}$.
\end{proposition}

For the rest of this section the statements of the proposition are
discussed and proved. This is just a modification of well known results
from linear algebra. Just as in linear algebra, the difficulty lies in not
getting lost in the notation. The reader with faith in the above
proposition can skip the rest of this section without harm.

\begin{proof}[Proof of Proposition \ref{prop.bcForRF}]
The choice of a free basis $\Aa=(\xn[e])$ carries along a natural unit
$R[F^e]$--structure $\theta^e_{\Aa}$ on $M$. Since the basis fixes an
isomorphism $M \cong R^{\oplus n}$ and each direct summand $R$ carries a
natural unit $R[F^e]$--structure we get a unit $R[F^e]$--structure on $M$. It
is given by sending $1 \tensor e_i$ to $e_i$. The corresponding Frobenius
action $F^e_{\Aa}$ is defined by $F^e_{\Aa}(e_i)=e_i$ and then extended
$p^e$--linearly.

Thus, the Frobenius structure on $M$ was obtained from the natural
Frobenius structure on the columnspace which a choice of basis
identifies with $M$. More generally note that all matrix spaces
carry a natural Frobenius action in which the Frobenius $F^e_R$ of
$R$ acts on each entry separately. This action on a matrix $M$ we
also denote by $F^e_R$, \ie we have that $F^e_R(M)=M^{[p^e]}$ is
the matrix whose entries are the $p^e$th powers of the entries in
$M$. Since matrices themselves are maps between matrix spaces (of
the appropriate sizes) the composition $F^e_R \circ M$ as maps is
different from $F^e_R(M)$. In fact we have $F^e_R \circ
M=F^e_R(M)F^e_R=M^{[p^e]}F^e_R$ as maps between matrix spaces. As
is general custom we often omit the composition symbol from the
notation and write $F^e_R M$ for $F^e_R \circ M$.

We fix the unit $R[F^e]$--structure $\theta^e_\Aa$ on $M$ which is
obtained from the basis $\Aa$. As seen in Lemma \ref{lem.thetaPhi} above
any $R[F^e]$--structure $(\theta^e,F^e)$ on $M$ is represented by an
endomorphism $\phi$ of $M$. Now, with respect to the basis $\Aa$, $\phi$
is given by some $n\times n$ matrix $A$ with entries in $R$. Thus $F^e$ is
represented by $AF^e_R$ with respect to the same basis. By abuse of
notation we just say that $A$ represents $\theta^e$ (or $F^e$) in this
context. Of course, $F^e_{\Aa}$ itself is represented by the identity
matrix $I$ with respect to the basis $\Aa$.

What is the matrix representing a power $F^{er}$ of the given Frobenius
action $F^e$? For this, first observe that the liner map
$F^{er}_{\Aa}(\phi) = \theta^{er}_{\Aa} \circ \F[er]_{\Aa}\phi \circ
{\theta^{er}_{\Aa}}^{-1}$ is represented by the matrix
$A^{[p^{er}]}=F^e_R(A)$, again with respect to the basis $\Aa$. This is
easily verified by hand; one observes that
$F^{er}_\Aa(\phi)(e_i)=\theta_{\Aa}^{er}(\id \tensor \phi(1 \tensor
e_i))=\theta_{\Aa}^{er}(1 \tensor \sum a_{ji}e_j)=\sum
a_{ji}^{p^{er}}e_j$. Thus the map $\phi_r$ (defined in Lemma
\ref{lem.highertheta}) is represented by the matrix $A_r = A A^{[p^e]}
\cdots A^{[p^{er}]}$ with respect to the basis $\Aa$.  Now, using Lemma
\ref{lem.highertheta}, we see that $F^{er}=\phi_r \circ F^{er}_\Aa$ and
therefore $F^{er}$ is represented by $A_rF^{er}_R$.

The next question is what happens when we change the basis. Let $\Bb =
(\xn[f])$ be another free basis of $M$ giving us the Frobenius structure
$\theta^e_{\Bb}$ and action $F^e_{\Bb}$. What is the matrix $B$
representing the given $F^e$ with respect to the basis $\Bb$? Let $C$ be
the matrix responsible for the base change; \ie $f_j = \sum c_{ij} e_i$.
Then, with respect to this new basis, $F^e$ is represented by
$C^{-1}AF^e_RC$. By definition of $F^e_R$, this is the same as
$C^{-1}AC^{[p^e]}F^e_R$. What this means is that if we write $F^e = \psi
\circ F^e_{\Bb}$, then $\psi \in \End_R(M)$ is represented by the matrix
$C^{-1}AC^{[p^e]}$ with respect to the basis $\Bb$. This finishes the
proof of Proposition \ref{prop.bcForRF}.
\end{proof}

\subsection{Injective hull: An Example}\label{sec.InjHullEx}

So far we only considered different $R[F^e]$--structures on an $R$--module
under the assumption that $M$ carries some \emph{unit}
$R[F^e]$--structure. If we drop this assumption the situation get more
complicated. One has to describe the set $\Hom_R(\F[e]M,M)$ which now is
no longer isomorphic to $\End_R(M)$. In general this is hard. As an
example we study the case that $M$ is the injective hull of the residue
field of a local ring.

Let $A$ be a local ring which is the quotient $A=R/I$ or a regular local
ring $R$. Let $E_A$ denote the injective hull of the residue field of $A$.
We denote the injective hull of the residue field of $R$ by $E=E_R$. First
we recall some basic properties of the injective hull.

\begin{lemma}\label{lem.HomErEQHomEa}
    With the notation just introduced, the injective hull\/ $E_A$ can be
    identified with\/ $\Ann_{E_R} I$. Furthermore, the functors\/
    $\Hom_A(\usc,E_A)$ and\/ $\Hom_R(\usc,E_R)$ are isomorphic on\/ $A$--mod.
\end{lemma}
\begin{proof}
By adjointness of $\Hom$ and tensor we have a natural isomorphism of
functors
\begin{equation*}\tag{\ltag}\label{eqn.loc.ErEa}
    \Hom_R(\usc,E_R) \cong  \Hom_A(\usc,\Hom_R(R/I,E_R)) =
    \Hom_A(\usc,\Ann_{E_R} I).
\end{equation*}
Thus, as $E_R$ is injective $\Hom_R(\usc,E_R)$ is exact and therefore
$\Ann_{E_R} I$ is injective (by definition of injectivity). It remains to
show that the inclusion of the residue field $k$ into $\Ann_{E_R}I$ is
essential. This also follows easily from the fact that $k \subseteq E_R$
is an essential map. Thus $\Ann_{E_R} I$ is an injective, essential
extension of $k$ and consequently it is an injective hull $E_A$. Equation
\eqnref{eqn.loc.ErEa} now shows also the last claim of the lemma.
\end{proof}
The next lemma shows a natural compatibility of annihilators with the
Frobenius. Even though we only apply it here in the case of the injective hull
it is true quite generally.
\begin{lemma}\label{lem.FcommutesAnn}
    Let\/ $(M,\theta^e,F^e)$ be an\/ $R[F^e]$--module. Then for any ideal\/ $I$ of\/
    $R$ we have an inclusion of submodules of\/ $M$
    \[
        F^e(\Ann_M I) \subseteq \Ann_M F^e_R(I).
    \]
    If\/ $R$ is regular and\/ $M$ is unit this is, in fact, an equality.
\end{lemma}
\begin{proof}
If $F^e(m)$ is in the left side, then $I^{[p^e]}F^e(m) = F^e(Im)=0$. Thus
$F^e(m)$ is annihilated by $F^e_R(I)=I^{[p^e]}$. The more interesting part
is that we have equality if $M$ is unit and $R$ regular. In this case it
is easy to verify that the inclusion above is in fact obtained as the
following sequence of maps
\begin{equation*}
\begin{split}
    F^e(\Ann_M I) &\cong \F[e]\Ann_M I \cong \F[e]\Hom(R/I,M) \cong
    \Hom(\F[e](R/I),\F[e]M) \\
                  &\cong \Hom(R/I^{[p^e]},M) \cong \Ann_M F^e(M).
\end{split}
\end{equation*}
The first map is the isomorphism ${\theta^e}^{-1}$. The third is the
natural isomorphism of $S \tensor \Hom(N,M) \cong \Hom(S \tensor N, S
\tensor M)$ for finitely presented $N$ and $R$--flat $S$ as furter
discussed below in Proposition \ref{prop.HomAndTensor}. Here $N = R/I$ is
finitely presented and $S = R^e$ is flat as $R$ is assumed regular. The
fourth map is obtained again from the isomorphism $\theta^e$ and the
natural identification of $\F[e](R/I)$ with $R/I^{[p^e]}$. If one follows
all these explicitly defined maps one sees that it is in fact just the
inclusion $F^e(\Ann_M I) \subseteq \Ann_M F^e(I)$. Thus this inclusion is,
in fact, equality.
\end{proof}
Armed with these two lemmata we can easily derive a satisfying
description of $\Hom_A(\F[e]_A E_A,E_A)$, the set of
$A[F^e]$--module structures on $E_A$. We tensor the natural
projection $R^e \to[\pi] A^e$ with $E_A$ and use Lemma
\ref{lem.HomErEQHomEa} to get
\[
    \Hom_A(A^e \tensor E_A,E_A) \cong \Hom_R(A^e \tensor E_A,E_R)
    \hookrightarrow \Hom_R(R^e \tensor E_A,E_R).
\]
This establishes $\Hom_A(A^e \tensor E_A,E_A)$ as the annihilator of $I$
in $\Hom_R(R^e \tensor E_A,E_R)$. To see this let $\phi: A^e \tensor E_A
\to E_A$ be an $A$--linear map. Its image in $\Hom_R(R^e \tensor E_A,E_R)$
is the map
\[
    R^e \tensor E_A \to[\pi \tensor \id] A^e \tensor E_A \to[\phi] E_A
\subseteq E_R.
\]
and already the first factor $\pi \tensor \id$ is annihilated by $I$.
Conversely, if a map $\psi: R^e \tensor E_A \to E_R$ is killed by $I$ its
image $\Image \psi$ is in $E_A = \Ann_{E_R} I$ and $I^e \tensor E_A$ is in
its kernel, which exactly means that $\psi$ factors thru a map $\phi$ as
above.

The next step is to find a better description of $\Hom_R(R^e \tensor
E_A,E_R)$. The natural inclusion $E_A \subseteq E_R$ gives a surjection
$\Hom_R(R^e \tensor E_R, E_R) \to \Hom_R(R^e \tensor E_A, E_R)$. Now we
use that for a regular local ring $R$ the injective hull $E_R$ is
isomorphic to the top local cohomology module of $R$, which, by Examples
\ref{ex.BasicRFmods}, is a unit $R[F^e]$--module. We fix such a unit
$R[F^e]$--structure $\theta^e:R^e\tensor E_R \to E_R$. Then applying
$\theta^e$ to the first arguments of the surjection just constructed we
get an isomorphic surjection
\[
    \Hom_R(E_R,E_R) \to \Hom_R(F^e(E_A),E_R).
\]
The left hand side is canonically isomorphic to $\hat{R}$, the completion
of $R$ along its maximal ideal. By Lemma \ref{lem.FcommutesAnn} the right
hand side is $\Hom_R(F^e(\Ann_{E_R}I),E_R) = \Hom_R(\Ann_{E_R}
I^{[p^e]},E_R) \cong \hat{R}/I\hat{R}$ and the map is just the natural
projection $\hat{R} \to \hat{R}/I^{[p^e]}\hat{R}$. Now, $\Hom_A(A^e
\tensor E_A,E_A)$ is just the annihilator of $I$ in $\hat{R}/I\hat{R}$ as
we just found out. Summarizing in a proposition we have:
\begin{proposition}
    Let\/ $A=R/I$ be a quotient of the regular local ring\/ $(R,n)$. The\/
    $A[F^e]$--module structures on\/ $E_A$ are in one-to-one correspondence
    with
    \[
        \frac{(I^{[p^e]}:_{\hat{R}}I)}{I^{[p^e]}\hat{R}}.
    \]
    Concretely, for\/ $w \in (I^{[p^e]}:I)\hat{R}$ the corresponding
    Frobenius action on\/ $E_A=\Ann_{E_R} I$ is given by\/ $\eta \mapsto
    w\eta^{p^e}$ where\/ $\eta^{p^e}$ is the standard Frobenius action on\/
    $E_R$ under the identification\/ $E_S \cong H^n_m(R)$.
\end{proposition}
\begin{proof}
Since the injective hull of $R$ is a unit $R[F^e]$--module, the module of
$R[F^e]$--structures on $E_R$ is just $\End_R(E_R)\cong \hat{R}$. With
this identification we showed in the preceding discussion that the set of
$A[F^e]$--structures on $E_A$ (\ie $\Hom_A(A^e \tensor E_A, E_A)$) is the
annihilator of $I$ in $\hat{R}/I^{[p^e]}\hat{R}$. It remains to observe
that an element in $(I^{p^e}:I)\hat{R}$ represents the claimed
$A[F^e]$--structure. This is also clear from the preceding discussion.
\end{proof}
\note{here might be a good place to include an example such that $E_R$ is
not unit, maybe this is characterized by gorensteinness}

It is well known that if $(A,m)$ is Gorenstein, then $E_A$ is in fact
isomorphic to the top local cohomology module $H^d_m(A)$, where $d$
denotes the dimension of $A$. According to Example \ref{ex.BasicRFmods}
this is always a $\emph{unit}$ $R[F^e]$--module. Thus in this case
$\Hom_A(A^e \tensor E_A,E_A) \cong \End_A(E_A) \cong \hat{A}$ and
therefore $\frac{(I^{[p^e]}:I)\hat{R}}{I^{[p^e]}\hat{R}}$ is a free
$\hat{A}$--module of rank one. This can also be observed with different
methods as is done in \cite{Fedder83}, for example.

\section{Change of $D_R[F]$--module structure}
So far we have fixed an $R$--module $M$ and studied the $R$--module of
$R[F^e]$--module structures on $M$. Specializing further, we now fix a
$D_R$--module $M$ and study the set of $D_R[F^e]$--module structures on
$M$. At first, this is fairly analogous to the case of
$R[F^e]$--structures and we will just comment on it. As we go along it
turns out that the additional rigidity we impose by fixing the
$D_R$--structure, forces this set of $D_R$--module structures to be fairly
small and controllable. In the cases we care about (\ie $M$ is semisimple
as a $D_R$--module) $\End_{D_R}(M)$ can be identified with a product of
matrix algebras over the algebraically closed field $k \subseteq R$. For
this we assume Claim \ref{cl.EndDsimp}. As a consequence we are able to
write every $D_R[F^e]$--module structure on a finitely generated
semisimple $D_R$--module $M$ as coming from a $k[F^e]$--structure on an
appropriate $k$--vector space. This will be the key ingredient for the
reduction step in the proof of Main Theorem 1.

For simplicity we assume here that $R$ is regular. Few of the things we
discuss are true more generally and not assuming regularity would
complicate the proofs. Assuming regularity makes things pleasant as
Frobenius descent can then be employed very efficiently. Let $M$ be a
$D_R$--module. By definition, a $D_R[F^e]$--module structure on $M$ is a
$D_R$--linear map $\theta^e: \F[e]M \to M$ where $\F[e]M$ carries its
natural $D_R$--structure given by Frobenius descent (cf.\ Proposition
\ref{prop.FrobDesc4Dmod}). Thus, the set of $D_R[F^e]$--module structures
on $M$ can be identified with $\Hom_{D_R}(\F[e]M,M)$. Naturally, these are
just the $R[F^e]$--structures on the $R$--module $M$ which are compatible
with the given $D_R$--structure, \ie those, for which the map $\theta^e:
\F[e]M \to M$ is $D_R$--linear. As before, if $M$ admits a unit
$D_R[F^e]$--structure, then $\Hom_{D_R}(\F[e]M,M)\cong \End_{D_R}(M)$. The
analog of Lemma \ref{lem.thetaPhi} in this $D_R$--module context now
reads.
\begin{lemma}\label{lem.thetaPhiDmod}
    Let\/ $(M,\theta^e)$ be a unit\/ $D_R[F^e]$--module. Every\/
    $D_R[F^e]$--structure\/ ${\theta'}^e$ on\/ $M$ can be written as\/
    ${\theta'}^e=\phi \circ \theta^e$ for some unique\/ $\phi \in
    \End_{D_R}(M)$.
\end{lemma}
Of course, the map $\phi$ in this proposition can be constructed
from ${\theta'}^e$. It is just $\phi = {\theta'}^e \circ
{\theta^e}^{-1}$ and thus ${\theta'}^e$ defines a unit structure
if and only if $\phi \in \Aut_{D_R}(M)$. The corresponding
Frobenius actions are also related by ${F'}^e = \phi \circ F^e$.

We now assume that $R$ is essentially of finite type over an algebraically
closed field $k$. A finitely generated semisimple $D_R$--module $M \cong
N_1^{\oplus n_1}\oplus \ldots \oplus N_r^{\oplus n_r}$ can then
equivalently be written as $M \cong \bigoplus (N_i \tensor_k V_i)$ where
$V_i$ is a $n_i$--dimensional $k$--vectorspace. The action of $\delta \in
D_R$ on $N_i \tensor_k V_i$ is given by $\delta(n \tensor v)=\delta(n)
\tensor v$. This is well defined since $k$ is algebraically closed and
therefore in particular perfect, thus $D_R$ acts trivially on $k$. This
setup leads to a few observations which we summarize in a lemma.

\begin{lemma}\label{lem.EndDEndk}
    With notation as indicated (and also assuming Claim \ref{cl.EndDsimp})
    we have\/ $\End_{D_R}(N \tensor V) \cong \End_k(V)$. Tensoring with\/ $N
    \tensor \usc$ gives a one-to-one correspondence between the\/
    $k$--vectorsubspaces of\/ $V$ and the\/ $D_R$--submodules of\/ $N \tensor
    V$.
\end{lemma}
\begin{proof}
After the choice of a basis for $V$, the ring $\End_{D_R}(M)$ is the
matrix algebra over $\End_{D_R}(N)=k$ of size $\dim_k V$ (here we use
Claim \ref{cl.EndDsimp}). This identifies $\End_{D_R}(M)$ with
$\End_k(V)$. Given a $\phi_k \in \End_k(V)$ the corresponding map in
$\End_{D_R}(N)$ is $\id_N \tensor \phi_k$.

Let $M'$ be a $D_R$--submodule of $M$. Since $M$ is semisimple we find a
$D_R$--submodule $M''$ such that $M' \oplus M'' \cong M$. Then $M'$ is the
kernel of the endomorphism $\phi: M \to[\pi_{M''}] M'' \subseteq M$. Thus
by the first part $\phi = \id_N \tensor \phi_k$ for some $\phi_k$ in
$\End_k(V)$. Then clearly $M' = N \tensor V'$ with $V' = \ker \phi_k$.
\end{proof}
Thus we see that if $M=N_1^{\oplus n_1} \oplus \ldots \oplus N_r^{\oplus
n_r}$ is semisimple, then $\End_{E_R}(M)$ is a finite direct sum of matrix
algebras of size $n_i$ over $k$.

We want to extend the last lemma so that it incorporates Frobenius
operations. This way we will be able to reduce questions about unit
$R[F^e]$--submodules of $M$ (\eg are there any nontrivial ones?) to the
equivalent questions about $k[F^e]$--subspaces of $V$. The main
observation is:
\begin{proposition}\label{prop.MreducetoV}
    Let\/ $M=N \tensor_k V$ be an isotypic\/ $D_R$--module (\ie $N$ is\/
    $D_R$--simple and\/ $V$ is a finite dimensional\/ $k$--vectorspace) that
    is a unit\/ $D_R[F^e]$--module with Frobenius action\/ $F^e_M$. Then there
    is a Frobenius action\/ $F^e_V$ on\/ $V$ such that the\/ $F^e_M$--stable\/
    $D_R$--submodules of\/ $M$ are in one-to-one correspondence with the\/
    $F_V$--stable subspaces of\/ $V$.
\end{proposition}
\begin{proof} First observe that by Proposition
\ref{prop.SimpSubofUnitisUnit}, $N$ is itself a unit
$D_R[F^e]$--module and we denote the corresponding Frobenius
action by $F^e_N$. The choice of a basis of $V$ equips $V$ with a
unit $k[F^e]$--structure as seen in Section \ref{sec.FreeRF}. This
Frobenius action we denote by ${F'}^e_V$. Then ${F'}^e_M \defeq
F^e_N \tensor {F'}^e_V$ defines a unit $D_R[F^e]$--structure. The
corresponding $R[F^e]$--structure on $M$ is given by
\[
    (R^e \tensor_R N) \tensor_k V \to N^e \tensor_k V \cong N \tensor_k
(k^e \tensor_k V) \to N \tensor_k V
\]
where the first map is the unit $R[F^e]$--structure on $N$ and the last is
the unit $k[F^e]$--structure on $V$. Since these both are isomorphisms, so
is the composition.

By Lemma \ref{lem.thetaPhiDmod}, we can express $F^e_M=\phi \circ
{F'}^e_M$ for some $\phi \in \Aut_{D_R}(M)$. By Lemma
\ref{lem.EndDEndk} we can write $\phi = \id_N \tensor \phi_V$ for
some $k$--vectorspace automorphism $\phi_V$ of $V$. Denoting the
corresponding unit $k[F^e]$--structure on $V$ by $F^e_V=\phi_V
\circ {F'}^e_V$, one easily verifies that $F^e_M = F^e_N \tensor
F^e_V$. Indeed,
\[
    F^e_N \tensor F^e_V = F^e_N \tensor (\phi_V \circ {F'}^e_V) = \phi \circ
    (F^e_N \tensor {F'}^e_V) = \phi \circ {F'}^e_M = F^e_M.
\]

With $F^e_V$ we have constructed the desired Frobenius action on $V$.
Since $F^e_M=F^e_N \tensor F^e_V$ we have $F^e_M(N \tensor V')=N \tensor
F^e_V(V')$. Therefore, a $D_R$--submodule $N \tensor V'$ of $M$ is stable
under $F^e_M$ if and only if the subspace $V'$ of $V$ is stable under
$F^e_V$. By Lemma \ref{lem.EndDEndk} every $D_R$--submodule of $M$ is of
this form. This finishes the proof.
\end{proof}
This proof is really just an reformulation of the correspondence
$\End_{D_R}(M) \cong \End_k(V)$ in terms of Frobenius structures. As $M$
is assumed to carry unit a $D_R[F^e]$--structure we get
\[
    \Hom_{D_R}(R^e \tensor M,M) \cong \End_{D_R}(M) \cong \End_k(V) \cong
\Hom_k(k^e\tensor V,V)
\]
where we apply Lemma \ref{lem.thetaPhiDmod} to both $M$ and $V$. Thus the
$D_R[F^e]$--structures of $M$ correspond to the $k[F^e]$--structures on
$V$. All the last proposition says is that this correspondence can be
achieved in way compatible with \ref{lem.EndDEndk}, \ie specializing to a
correspondence between the $F^e$--stable submodules of $M$ and
$F^e$--stable subspaces of $V$.

\section{Proof of Main Theorem 1}
As already indicated the proof of Main Theorem 1 will be a reduction to
the case that $M$ is a finite dimensional vector space over an
algebraically closed field $k$. Thus we treat this case first.

\subsection{The case of a vector space}
Finite dimensional $k$--vectorspaces with a Frobenius action are well
studied as a part of the theory of finite dimensional $p$--Lie algebras
(see for example the many publications of Dieudonn\'e on the subject). For
our purposes, we only need one result from this theory \cite[Proposition
3, page 233]{DieudonLie}. For convenience we review a slightly altered
version here and also recall the proof.
\begin{proposition}\label{prop.main1VScase}
    Let\/ $k$ be an algebraically closed field and\/ $V$ a finite
    dimensional\/ $k$--vectorspace. Let\/ $F^e$ be a Frobenius action on
    $k$. Then\/ $V$ has a one dimensional\/ $F^e$--stable subspace. In
    particular, $V$ is a simple $k[F^e]$--module if and only if\/ $V=k$,
    \ie $V$ is simple as a\/ $k[F^e]$--module if and only if\/ $V$ is
    simple as a\/ $k$--vectorspace.
\end{proposition}
\begin{proof}
    If there is $v \in V$ such that $F^e(v)=0$, then the subspace
    generated by $v$ is stable under $F^e$. Thus we can assume that $F^e$
    acts injectively. Let $v \in V$ and replace $V$ by the $k$--span of
    $(v, F^e(v), F^{2e}(v), \dots )$. Thus we can assume that for some $r
    \geq 0$ the tuple $(v,F^e(v), \ldots,F^{re}(v))$ is a basis of $V$.
    For brevity we denote $p^e$ by $q$.

    We have to find $w$ such that $F^e(w)=\lambda w$. We can assume that
    $\lambda=1$ since if we replace $w$ by $\lambda^{-1/q}w$ ($k$
    algebraically closed), then $F^e(w)=w$. If we write $w$ as well as
    $F^{e(r+1)}(v)$ in terms of the above basis
    \begin{eqnarray*}
        w &=& a_0v + a_1F^e(v) + \ldots + a_{r}F^{er}(v) \\
        F^{e(r+1)}(v) &=& b_0v + b_1F^e(v) + \ldots + b_{r}F^{er}(v)
    \end{eqnarray*}
    then we get from comparing coefficients of the required identity
    $F(w)=w$ the following equations:
    \begin{eqnarray*}
         a_0 &=& a_r^qb_0\\
         a_1 &=& a_0^q + a_r^qb_1\\
         a_2 &=& a_1^q + a_r^qb_2\\
        &\vdots& \\
         a_r &=& a^q_{r-1} + a_r^qb_r
    \end{eqnarray*}
    We see that all the $a_i$'s for $i<r$ are determined by $a_r$.
    Successively substituting from top to bottom and setting $a_r=t$ we
    get an algebraic equation in $t$:
    \[
        t=t^{q^{r+1}}b_0^{q^r}+t^{q^r}b_1^{q^{r-1}}+ \ldots + t^qb_r
    \]
    As $k$ is algebraically closed we can find a nonzero solution $a$ of
    this algebraic equation. Setting $a_r=a$ will determine a solution for
    $F(w)=w$ by calculating the $a_i$ according to the equations above.
    This $w$ generates a one dimensional $F^e$--stable subspace of $V$.
\end{proof}
In fact, this proof shows that if $k$ is not algebraically closed one may
have to go to a finite algebraic extension $K$ to find a one dimensional
$F$--stable subspace of $K \tensor_k V$. If we use this proof as the start
of an induction we get the following corollary.
\begin{corollary}\label{cor.VhasFfixedBasis}
    Let\/ $k$ be an algebraically closed field and\/ $V$ a finite dimensional\/
    $k$--vectorspace. If the Frobenius\/ $F^e$ acts injectively on\/ $V$, then\/
    $V$ has a basis consisting of\/ $F^e$--fixed elements of\/ $V$.
\end{corollary}
\begin{proof}
In the above proof, continuing with the $k[F^e]$--module quotient $W$, we
can assume by induction on $r$ that $W$ has a basis $w_1 \ldots w_{r-1}$
such that $F_W^e(w_i)=w_i$. Taking nonzero elements $w'_i \in V$ which get
mapped to $w_i$ by the quotient map $V \to W$ we see that $F^e_V(w'_i) =
w'_i + a_iw$ for some $a \in k$. For $w''_i \defeq w'_i - a_i^{1/p^e}w$ we
have $F_V^e(w''_i) = w''_i$. Thus $(w, w''_1, \ldots, w''_{r-1})$ is a
basis of $F^e$--fixed elements.
\end{proof}

\subsection{Reduction to the vectorspace case}
With all the tools developed in this third chapter, it is quite simple to
prove Main Theorem 1 by reduction to Proposition \ref{prop.main1VScase} .

\mainthmonecontent   %This includes exactly what was typeset where the Main Theorem first appeared

\begin{proof}
For some $e>0$ the simple unit $R[F^\infty]$--module $M$ is a simple unit
$R[F^e]$--module. Theorem \ref{thm.unitRFsemisimpleDR} shows that $M$ is
semisimple as a $D_R$--module. By Proposition \ref{prop.IsotypicareRFsub}
each $D_R$--isotypic component of $M$ is a unit $R[F^{er}]$--submodule for
some $r>0$. As $M$ is simple as an $R[F^\infty]$--module, it is in
particular simple as an $R[F^{er}]$--module. Therefore $M$ must be
$N$--isotypic for some simple $D_R$--module $N$; \ie $M \cong N \tensor_k
V$ as $D_R$--modules for some finite dimensional $k$--vectorspace $V$.

We replace $e$ by $er$ and therefore again assume that $M$ is in fact a
simple unit $R[F^e]$--module which is isotypic as a $D_R$--module. We
denote the Frobenius action by $F^e_M$. We have seen in Proposition
\ref{prop.DRfulluRF} that the $F^e_M$ stable $D_R$--submodules of $M$ are
exactly the unit $R[F^e]$--submodules of $M$.  By Proposition
\ref{prop.MreducetoV} there is a Frobenius action $F^e_V$ on $V$ such that
the $F^e_M$--stable $D_R$--submodules of $M$ correspond to the $F^e_V$
stable submodules of $V$. Thus $M$ is a simple unit $R[F^e]$--module if
and only if $(V,F^e_V)$ is a simple $k[F^e]$--module. By Proposition
\ref{prop.main1VScase}, if $V$ is a simple $k[F^e]$--module, then $V$ is
one dimensional. Therefore, $M \cong N$ and thus $M$ is simple as a
$D_R$--module.
\end{proof}
\begin{remark}\label{rem.MainOneGeneral}
Note that by Remark \ref{rem.EndDNgeneral} just after the proof of Lemma
\ref{lem.EndNisAlgebraic}, it is enough to assume that $R$ is a regular
$k$ algebra such that the cardinality of the algebraically closed field
$k$ is strictly bigger than the cardinality of a $k$--basis of $R$. Thus
Theorem 1 is true with this more general hypothesis.
\end{remark}

Even though this shows the strong connection between the categories of
$D_R$ modules and finitely generated unit $R[F^e]$--modules it is not
enough to conclude that the latter is a full subcategory of the former. In
fact, a given simple $D_R$--submodule $N$ of a finitely generated unit
$R[F^e]$--module $M$ may not be fixed by \emph{any} power of the Frobenius
action as the following examples show. What is the case is that we can
find an isomorphic copy of $N$ inside of $M$ which is in fact a unit
$R[F^{er}]$--submodule. We have the following corollary:
\begin{corollary}\label{cor.MainOne}
    Let\/ $R$ be regular and essentially of finite type over an uncountable
    algebraically closed field\/ $k$. Let\/ $M$ be a finitely generated unit\/
    $R[F^\infty]$--module and\/ $N$ a\/ $D_R$--simple submodule of\/ $M$. For
    some\/ $e \geq 0$, the unit\/ $R[F^e]$--submodule\/ $R[F^e]N$ decomposes
    into a finite direct sum of simple\/ $R[F^e]$-submodules of\/ $M$; as\/
    $D_R$--modules, each summand is isomorphic to\/ $N$.
\end{corollary}
\begin{proof}
As in the proof of Theorem \ref{thm.unitRFsemisimpleDR}, $R[F^e]N = N
\oplus F^e(N) \oplus \ldots \oplus F^{e(r-1)}(N)$ is a unit
$R[F^e]$--submodule of $M$ for some $r$. Thus $R[F^e]N$ is semisimple as a
$D_R$--module. By Proposition \ref{prop.IsotypicareRFsub}, the $N$
isotypic component of $R[F^e]N$ is a unit $R[F^{er}]$--module for some
$r$. Increasing $e$ to $er$ makes $M'
\defeq R[F^e]N$ a $D_R$--semisimple and $N$--isotypic finitely generated
unit $R[F^e]$--submodule of $M$. As in the last proof we write $M' \cong N
\tensor_k V$ and find a Frobenius action $F^e_V$ on $V$ such that the
$F^e_V$--stable subspaces of $V$ correspond to the unit
$R[F^e]$--submodules of $M'$. By Corollary \ref{cor.VhasFfixedBasis}, $V$
decomposes into a direct sum of one dimensional $F^e_V$--stable subspaces
and thus $M'$ decomposes into a direct sum of $D_R$--simple unit
$R[F^e]$--submodules. Since $M'$ was $N$--isotypic these are all
isomorphic to $N$ as $D_R$--modules.
\end{proof}

\section{Examples and Counterexamples}\label{sec.ExAndCounterEx}
Here we collect some examples of $R[F^e]$--modules to illuminate the
theory developed in this chapter. In fact, we really have just one
recurring example of a free $R$--module $M$ of rank $2$ where the action
of the Frobenius $F^e$ with respect to some basis $\Aa=(e_1,e_2)$ is
represented by the matrix
\[
    A= \begin{pmatrix} 0 & 1 \\ 1 & x \end{pmatrix}
\]
for some element $x \in R$. What will change from example to example is
$R$ and $x$. The choice of basis $\Aa$ also induces a natural
$D_R$--structure on $M \cong R\oplus R$ by acting componentwise on the
summands. Unless otherwise specified, if we speak of the $D_R$--module
$M$, this is the structure we have in mind.

From Section \ref{sec.FreeRF}, we recall that the matrix $A_r$
representing the $r$TtH power of this Frobenius action $F^e$ with respect
to this basis is given by
\[
    A_R= AA^{[q]}\cdots A^{[q^{r-1}]}.
\]
Equivalently, $A_r$ can be described inductively by the equation $A_r =
A_{r-1}A^{[q^{r-1}]}$ which translates into an inductive formula for the
coefficients of $A_r$. One has
\begin{equation}
    A_r=\begin{pmatrix} a^q_{r-2} & a^q_{r-1} \\
                       a_{r-1}   & a_r
        \end{pmatrix}
\end{equation}
where $a_r=a_{r-2}+a_{r-1}x^{q^{r-1}}$ with $a_{-1}=0$ and $a_0=1$. So,
for example, this formula computes
\[
    A_1 = \begin{pmatrix} {0} &{1} \\ {1} &{x} \end{pmatrix}\ \text{ and }\
    A_2 = \begin{pmatrix} {1} &{x^q} \\ {x} &{x^{q+1}+1}
\end{pmatrix}
\]
which can be easily verified by hand. We prove these assertions by
induction
\begin{equation*}\begin{split}
    A_r = A_{r-1}A^{[q^{r-1}]} &=
      \begin{pmatrix}
            a^q_{r-2} & a^q_{i-3}+a^q_{r-2}x^{q^{r-1}}  \\
            a_{i-1}   & a_{a-2}+a_{i-1}x^{q^{r-1}}
      \end{pmatrix} \\  &= \begin{pmatrix}
                            a^q_{r-2} & (a_{r-3}+a_{r-2}x^{q^{r-2}})^q \\
                            a_{r-1}   & a_r
                     \end{pmatrix} = \begin{pmatrix} a^q_{2-1} & a^q_{r-1} \\
                                                    a_{r-1}   & a_r
                                    \end{pmatrix}
\end{split}\end{equation*}
where the base case is by the initial condition of the recursion for $a_r$.
Furthermore we note that thinking of $a_r$ as polynomials in $x$, the
degree in $x$ of $a_r$ is $\deg a_r = 1 + q + q^2 +\ldots+q^{r-1}$. Again
an induction argument shows this nicely:
\begin{equation*}\begin{split}
    \deg a_r &= \max\{\deg a_{r-2}, \deg a_{r-1}+q^{r-1}\} \\
             &= \max\{1+q+\ldots q^{r-3},1+q+\ldots+q^{r-2}+q^{r-1}\} \\
             &= 1+q+\ldots+q^{r-2}+q^{r-1}
\end{split}\end{equation*}
and for the start of the induction we just recall that $a_{-1}=0$ and
$a_0=1$. Now we turn to more concrete situations.
\begin{example}[A simple \protect{$R[F]$--module that is not $D_R$}--simple]
Let $R=\FF_3$ and $x=1$. Since the Frobenius $F$ is the identity on
${\FF_3}$, Frobenius actions are just linear maps. The linear map
represented by
\[
   A = \begin{pmatrix} 0 & 1 \\ 1 & 1 \end{pmatrix}
\]
is not diagonalizable since its characteristic polynomial
$P_A(t)=t(t-1)+1$ is irreducible over $\FF_3$. Thus $M$ is a simple
${\FF_3}[F]$--module. Since ${\FF_3}$ is perfect $D_{\FF_3}={\FF_3}$ and
thus $M$ is not simple as a $D_{\FF_3}$--module as it is a free
${\FF_3}$--module of rank 2. Note that, since $\End_{\FF_3}(M)$ is finite,
some power of $F$ will be the identity on $M$. An easy calculation shows
that $F^4=-\id_M$ and therefore $M$ is not simple as a
${\FF_3}[F^4]$--module.
\end{example}

\begin{example}[\protect{A $D_R$--submodule that is not an $R[F^\infty]$}--submodule]\label{ex.DRnorRFinfty}
Now let $R$ be a ring containing an infinite perfect field $k$ with $x \in
k$ transcendental over the prime field $\FF_p$. With this $x$ the matrix
$A$, in fact, represents a $D_R$--linear map of $M$ (this is because
differential operators $D_R$ are linear over any perfect subring of $R$).
As $a_r$ is a nonzero polynomial in $x$ with coefficients in $\FF_p$, and
since $x$ is transcendental over the prime field, $a_r$ is a
nonzero element of $k$. This implies that, for example, $Re_1$ is not
stable under any power of $F^e$ as this would be equivalent to the matrix
$A_r$ having a zero entry in the bottom left corner. But this entry is
$a_{r-1}$ for which we just argued is nonzero. Thus $Re_1$ is a (simple)
$D_R$--submodule that is not a unit $R[F^\infty]$--submodule of $M$.
\end{example}

\subsection{A simple $R[F^\infty]$--module that is not $D_R$--simple}
Now we come to the main example of a simple $R[F^\infty]$--module that is
not $D_R$-simple. This will show that the assumption of algebraic closure
in Main Theorem 1 cannot be dropped. Furthermore, it provides a
counterexample to Lyubeznik's Remark 5.6a in \cite{Lyub} where he claims
that simple $D_R$--submodules of simple unit $R[F^e]$--modules are in fact
$R[F^\infty]$--submodules.

Again, we recycle the examples above. With $A$ and $M$ as before, let
$R=k(x)^{1/p^\infty}$ where $x$ is a new variable and $k$ is perfect. We
denote by $F^e_\Aa$ the Frobenius action arising from the choice of basis
$\Aa$. The Frobenius action $F^e$ on $M$ we define as being represented by
$A$, \ie given by application of
\[
    F^e = \begin{pmatrix} 0 & 1 \\ 1 & x \end{pmatrix} F^e_R.
\]
with respect to the basis $\Aa$. Since $R$ is perfect, $R=D_R$ and
consequently this (and every) $F^e$--structure is compatible with the
$D_R$--structure.

We first want to show that with this Frobenius action, $M$ is a simple
$R[F^e]$--module. Assume the contrary and let $v=(\alpha, \beta)^t$ be
such that $F^e(v)=\lambda v$. If $\beta=0$, Example \ref{ex.DRnorRFinfty}
gives an immediate contradiction. Thus $\beta \neq 0$ and we change $v$ to
$\frac{1}{\beta}v$. Therefore we can assume that $v=(\alpha, 1)^t$ for
some $\alpha \in R$. From the assumption $F^e(v)=\lambda v$ we get
\[
    \begin{pmatrix}\lambda\alpha \\ \lambda \end{pmatrix} = \lambda v
    = F^e(v) = A\begin{pmatrix}\alpha^q \\ 1 \end{pmatrix} =
    \begin{pmatrix}1 \\ \alpha^q + x \end{pmatrix}.
\]
Comparing entries of the vectors at both ends of this equality and plugging
the second into the first we get a monic algebraic equation for $\alpha$
with coefficients in $\FF_p[x] \subseteq k[x]$:
\begin{equation}\label{eqn.zalpha}
    \alpha^{q+1} +x\alpha - 1 = 0
\end{equation}
As $k[x]^{1/p^\infty}$ is integrally closed in its field of fractions $R$,
we conclude that $\alpha \in k[x]^{1/p^\infty}$. Now let $t$ be smallest
such that $\beta(x) \defeq \alpha^{p^t} \in k[x]$; \ie $\alpha \in
k[x^{1/p^t}]$ and $\beta(x)$ is not a $p$-TtH power unless $e=0$. Now we
take the $p^e$TtH power of the Equation \eqnref{eqn.zalpha} and get
\begin{equation}\label{eqn.zbeta}
    \beta^{q+1} +x^{p^t}\beta - 1 = 0.
\end{equation}
This is an algebraic equation for $x$ over $k$. As $x$ is transcendental
over $k$ the left side must be the zero polynomial. To finally reach a
contradiction we discriminate two cases:
\begin{descrip}{$t>0\ :$}
\item[$t>0$\ :] Differentiate \eqnref{eqn.zbeta} with respect to $x$
    and get:
    \[
        \beta^q \diff \beta + x^{p^t} \diff \beta = 0
    \]
    By minimality of $e$, $\beta$ is not a $p$-TtH power and therefore
    $\diff \beta$ is nonzero. Thus we conclude that $\beta^q=-x^{p^t}$.
    Substituting this back into \eqnref{eqn.zbeta} one gets the
    contradiction $1=0$.
\item[$t=0$\ :] Let $n$ be the degree of $\beta$. Comparing highest
    degrees in \eqnref{eqn.zbeta} we get $(q+1)n = 1+n$ and thus $pn=1$
    which is a contradiction.
\end{descrip}
Thus we conclude that $(M,F^e)$ is a simple $R[F^e]$--module. To show that
it is simple as an $R[F^\infty]$--module we show that $M$ is a simple
$R[F^{er}]$--module for all $r$. The argument is similar to the case $r=1$
but with some extra twists. We proceed in 3 Steps:
\begin{descrip}{\textsc{Step 3}}
\item[\textsc{Step 1}]
    With respect to the basis $\Aa$ the action $F^{er}$ is represented by
    the matrix $A_r$. We change the basis appropriately to $\Bb=(f_1,f_2)$
    such that the representing matrix $B_r$ of $F^{er}$ with respect to
    the basis $\Bb$ is ``nice''; by this we mean that
    \[
     B_r = \begin{pmatrix} 0 & s_r \\ 1 & t_r \end{pmatrix}
    \]
    for some $s_r, t_r \in \FF_p[x]$.
\item[\textsc{Step 2}]
    As in the case $r=1$ we assume that there is $v \in M$ such that
    $F^{er}(v) = \lambda v$. This yields a monic algebraic equation which gives an
    algebraic equation for $x$.
\item[\textsc{Step 3}]
    We finish by discriminating the same cases as for $r=1$. One is
    treated by a degree argument, the other by differentiation.
\end{descrip}

Let us begin with \textsc{Step 1}: The basis that will lead to the matrix
$B_r$ of the desired shape is $f_1 = e_1$ and $f_2 = a^q_{r-2}e_1 +
a_{r-1}e_2$. Thus the matrix responsible for the base change from $\Aa$ to
$\Bb$ is
\[
    C_r \defeq \begin{pmatrix} 1 & a^q_{r-2} \\ 0 & a_{r-1} \end{pmatrix}.
\]
With the formalism developed in Section \ref{sec.FreeRF} we see that with
respect to the new basis $\Bb$ the Frobenius action $F^{er}$ is
represented by the matrix $B_r = C_r^{-1}A_r C_r^{[q^{r}]}$ (cf.\
Proposition \ref{prop.bcForRF}). To determine $s_r$ and $t_r$ we
explicitly calculate $B_r$:
\begin{equation*}
\begin{split}
B_r &= C_r^{-1}A_r C_r^{[q^r]} = C_r^{-1}
\begin{pmatrix}
    {a_{r-2}^q}&{a_{r-1}^q}\\
    {a_{r-1}}&{a_r}
\end{pmatrix}
\begin{pmatrix}
    {1}&{a_{r-2}^{q^r+q}}\\
    {0}&{a_{r-1}^{q^r}}
\end{pmatrix} \\
&=\frac{1}{a_{r-1}}
\begin{pmatrix}
    {a_{r-1}}&{-a_{r-2}^q} \\
    {0}& {1}
\end{pmatrix}
\begin{pmatrix}
    {{a_{r-2}^q}} & {{a_{r-2}^q a_{r-2}^{q^r+q} + a_{r-1}^q a_{r-1}^{q^r}}} \\
    {{a_{r-1}}} & {{a_{r-1}a_{r-2}^{q^r+q} + a_r a_{r-1}^{q^r}}}
\end{pmatrix} \\
&= \frac{1}{a_{r-1}}   % This might not be necessary...
\begin{pmatrix}
    {a_{r-1}a_{r-2}^q-a_{r-2}^q a_{r-1}}&
    \scriptstyle{{a_{r-1}a_{r-2}^{q^r+2q}+a_{r-1}^{q^r+q+1}-a_{r-1}a_{r-2}^{q^r+2q}-a_ia_{r-1}^{q^r}a_{r-2}^q}} \\
    {a_{r-1}} & {a_{r-1}a_{r-2}^{q^r+q} + a_r a_{r-1}^{q^r}}
\end{pmatrix}\\
&=
\begin{pmatrix}
    {0} &{-a_{r-1}^{q^r-1}\det A_r}\\
    {1} &{a_{r-2}^{q^r+q} + a_ra_{r-1}^{q^r-1}}
\end{pmatrix} \\
&=
\begin{pmatrix}
    {0}&{(-1)^{r-1}a_{r-1}^{q^r-1}} \\
    {1}&{a_{r-2}^{q^r+q} + a_ra_{r-1}^{q^r-1}}
\end{pmatrix}
\end{split}
\end{equation*}
Besides index juggling skills one only needs the equation $\det A_r =
(-1)^r$ which follows from the  recursive definition of $A_r$ and the fact
that $\det A=-1$. We can read off the desired expressions for $s_r$ and
$t_r$.
\[
    s_r=(-1)^{r-1}a_{r-1}^{q^r-1} \text{\quad and \quad}
t_r=a_{r-2}^{q^r+q} + a_r a_{r-1}^{q^r-1},
\]
Note that both are in $\FF_p[x]$ since $a_i \in \FF_p[x]$. We now
work over this new basis $\Bb$ and start with \textsc{Step 2}:
Assume we have $v=(\alpha,1)^t$ such that ${F'}^r(v)=\lambda v$
(as in the case $r=1$ one reduces from a general
$v=(\alpha,\beta)^t$ to this case with the help of Example
\ref{ex.DRnorRFinfty}). Considering $\lambda v = B_rF_R^{er}v$,
this yields 2 equations:
\begin{eqnarray*}
    s_r & = & \lambda \alpha \\
    \alpha^{q^r}+t_r & = & \lambda
\end{eqnarray*}
Substituting the latter in the former we get a monic algebraic equation
for $\alpha$ with coefficients in $k[x]$:
\begin{equation*}
    \alpha^{q^r+1} + t_r \alpha - s_r = 0
\end{equation*}
As $k[x]^{1/p^{\infty}}$ is integrally closed in $R$ we conclude that
$\alpha \in k[x]^{1/p^{\infty}}$. We choose $t$ minimal such that $\alpha
\in k[x^{1/p^t}]$. Then $\beta(x)=\alpha^{p^t}$ is in $k[x]$ and not a
$p$-TtH power unless $t=0$. Taking the $p^t$-TtH power of the last
equation we get
\begin{equation}\label{eqn.betar}
    \beta^{q^r+1} + t_r^{p^t} \beta - s_r^{p^t} = 0
\end{equation}
which is an algebraic relation for $x$ with coefficients in $k$. Thus it
is constant zero by transcendence of $x$. For \textsc{Step 3} we
distinguish again the following 2 cases:
\begin{descrip}{$t > 0\ :$}
\item[$t > 0$ :] Differentiating \eqnref{eqn.betar} with respect to $x$ we get:
\[
    \beta^{q^r}\diff \beta + t_r^{p^t} \diff \beta = 0
\]
As we chose $\beta$ not to be a $p$TtH power, its derivative is nonzero.
Thus we can divide the above by $\diff \beta$ and get
$\beta^{q^r}=-t_r^{p^t}$. Substituting this back into \eqnref{eqn.betar}
we get $s_r^{p^t}$=0. But this is a contradiction as $s_r=\pm
a_{r-1}^{q^r-1} \neq 0$.
\item[$t=0$ :]For this we have to determine the degrees of the terms in Equation
\eqnref{eqn.betar}. As observed earlier $\deg(a_r)=1+q+\ldots+q^{r-1}$.
Therefore,
\begin{eqnarray*}
  \hbox{\phantom{$t=0$}}\deg{s_r}&=&(q^r-1)\deg(a_{r-1})=-1-q-\ldots-q^{r-2}+q^r+\ldots+q^{2r-2} \\
  \deg{t_r} &\leq&
\max\{\deg(a_{r-2}^{q^r+q}),\deg({a_ia_{r-1}^{q^r-1}})\}=q^{r-1}+\ldots+q^{2r-2}.
\end{eqnarray*}
In fact equality prevails in the last inequality since the two entries in
the $\max$ are different (the second is always bigger). To be precise:
\begin{eqnarray*}
\deg(a_{r-2}^{q^r+q}) &=&  (q^r+q)(1+\ldots+q^{r-3}) \\
                      &=& q+\ldots+q^{r-2}+q^r+\ldots+q^{2r-3} \\
\hbox{\phantom{$t=0$}}\deg({a_ia_{r-1}^{q^r-1}}) &=& 1+q+\ldots+q^{r-1}-1-\ldots-q^{r-2}+q^r+\ldots+q^{2r-2} \\
                           &=& q^{r-1}+\ldots+q^{2r-2}
\end{eqnarray*}
Since $q^{r-1} > 1+q+\ldots+q^{r-2}$ we see that the second line is in
fact strictly bigger than the first. Thus the degree of $s_r$ is strictly
smaller than the degree of $t_r$, and therefore the first two terms of
\eqnref{eqn.betar} must have the same degree. If we denote the degree of
$\beta$ by $n$ we get
\[
    (q^r+1)n = \deg(t_r)+n
\]
and after dividing by $q^{r-1}$ this simplifies to
\[
    1 = qn-q-q^2-\ldots-q^{r-1}.
\]
The right side is divisible by $q$ but the left side certainly isn't. This
is a contradiction.
\end{descrip}
This finishes the proof that $M$ is a simple $R[F^{er}]$--module for all
$r>0$. Thus $M$ is a simple $R[F^\infty]$--module but $M$ is not simple as
a $D_R$--module since every one dimensional $R$--subspace is a nontrivial
$D_R$--submodule.

\subsection{Examples over the polynomial ring}
So far the examples were over a field. Starting with these examples it is
not hard to obtain equivalent examples over higher dimensional rings.
For this let $(V,F^e)$ be the simple unit $K[F^e]$--module of the last
example ($K = \FF_p(x)^{1/p^\infty}$). Let $R$ be a regular $K$--algebra,
essentially of finite type over $K$ (\eg $R=K \tensor_{\FF_p}
\FF_p[\xn]=K[\xn]$). Let $M = N \tensor_K V$ where $(N, F^e_N)$ is a
finitely generated unit $R[F^e]$--module which is simple as a
$D_R$--module ($R$ itself will work). Then $M$ carries a natural
$R[F^e]$--structure defined by $F^e_M(n \tensor v) = F^e_N(n) \tensor
F^e(v)$. This is exactly the situation of Proposition
\ref{prop.MreducetoV} applied in the reverse direction as we used it in
the proof of Main Theorem 1. Here we are given the Frobenius action $F^e$
on $V$ and construct an appropriate action $F^e_M = F^e_N \tensor F^e_V$
on $M$. Proposition \ref{prop.MreducetoV}, or merely its proof, then shows
that the unit $K[F^e]$--submodules of $V$ are in one-to-one correspondence
with the unit $R[F^e]$--submodules of $M$. Thus, the simplicity of $V$ as
a unit $K[F^\infty]$--module implies that $M$ is a simple unit
$R[F^\infty]$--module. Clearly, since $V$ is not $D_K$--simple, $M$ is not
$D_R$--simple.

However, we must be careful. In Proposition \ref{prop.MreducetoV} we
assumed Claim \ref{cl.EndDsimp}; \ie we assumed that $\End_{D_R}(N)$
is algebraic over $k$ and that $k$ is algebraically closed. What was
really used is that $\End_{D_R}(N)=k$ (see proof of Lemma
\ref{lem.EndDEndk}). So the above argument is valid whenever
$\End_{D_R}(N)=K$. At least in one case this is true: For $N = R$ we have
$\End_{D_R}(R) = \bigcap R^{p^e} =K$ in the case that $R$ is essentially
of finite type over the perfect field $K$. This shows that $R \tensor V$
is an example of a simple $R[F^\infty]$--module which is not simple as a
$D_R$--module.

%%%%%%%%%%%%%%%%%%%%%%%%%%%%%%%%%%%%%%%%%%%%%%%%%%%%%%%%%%%%%%
%%                                                          %%
%%   This is file: chapter4.tex                             %%
%%   It contains the Third Chapter of my                    %%
%%   dissertation: dissertation.tex                         %%
%%                                                          %%
%%%%%%%%%%%%%%%%%%%%%%%%%%%%%%%%%%%%%%%%%%%%%%%%%%%%%%%%%%%%%%

\chapter{Functors on $R[F]$--modules}\label{chap.Functors}

In this chapter we investigate how certain functors on the category of
$R$--modules behave when Frobenius actions are present. The goal we have
in mind is to extend the \idx{Matlis duality} functor $D=\Hom(\usc,E_R)$
in a way that incorporates Frobenius structures. The extension, $\Dd$, of
Matlis duality we describe, is the main tool to establish the connection
between the tight closure of zero in $H^d_m(R/I)$ and the unique simple
$D_R$--submodule of $H^c_I(R)$. Thus $\Dd$ is an integral part in the
construction of $\Ll(A,R)$ which is obtained in the next chapter.

To obtain this extension of the Matlis duality functor we take somewhat of
a detour. We begin with establishing a right adjoint to the forgetful
functor from unit $R[F^e]$--modules to $R[F^e]$--modules. Such an adjoint
is nothing but a natural way to make an $R[F^e]$--module into a
\emph{unit} $R[F^e]$--module. This right adjoint is found as Hartshorne
and Speiser's \idx{leveling functor} $G$ \cite{HaSp}. We review their
construction and show that it is, in fact, an adjoint, \ie we show that
for an $R[F^e]$--module $M$ and unit $R[F^e]$--module $N$ we have
$\Hom_{R[F^e]}(N,M) \cong \Hom_{R[F^e]}(N,G(M))$. Furthermore, as an
application of Frobenius descent, we recall a construction of Emerton and
Kisin \cite{Em.Kis2}, establishing a left adjoint of the forgetful functor
from unit $D_R[F^e]$--module to $D_R[F^e]$--modules.

Then we note some abstract constructions about extending functors on
$R$--mod to functors on $R[F^e]$--mod. In general it is fairly
straightforward to extend a covariant functor \label{x.CovFunct} $C: \Rmod
\to \Rmod[A]$ for commutative rings $R$ and $A$ to a functor from
$\Rmod[{R[F^e]}]$ to $\Rmod[{A[F^e]}]$. For this one only needs to have a
natural transformation of functors $\Psi: \F[e]_A \circ C \to C \circ
\F[e]_R$. Then, if $\theta^e$ is an $R[F^e]$--structure on $M$, the
composition
\[
    \F[e]_A(C(M)) \to[\Psi_M] C(\F[e]_R(M)) \to[C(\theta^e)] C(M)
\]
defines an $A[F^e]$--structure on $C(M)$. We already employed this when we
showed that for a map of rings $R \to[\pi] A$ the functor $\pi^* = A
\tensor_R \usc$ extends to a functor from $R[F^e]$--modules to
$A[F^e]$--modules (cf.\ beginning of Section \ref{sec.PropRFmod}). Also,
the local cohomology functors commute with the Frobenius for regular $R$
and thus they are naturally extended to functors on $R[F^e]$--modules.

More demanding is the case of contravariant functors. If $K$ denotes a
\idx{contravariant functor} from $R$--mod to $A$--mod it is, \emph{a
priori}, not clear how one obtains from an $R[F^e]$--structure $\theta^e$
on $M$ an $A[F^e]$--structure on $K(M)$. Even with a compatibility between
$\F[e]$ and $K$, the map $K(\theta^e)$ somehow points in the wrong
direction. This is not a problem if $(M,\theta^e)$ is unit, but if one is
interested in not necessarily unit modules (and we are), more drastic
measures are necessary. Roughly, one passes to the directed system
generated by $K(\theta^e)$ and defines $\Kk(M)$ as its direct limit. Then,
given a natural transformation $\Psi:\F[e]_A \circ K \to K \circ \F[e]_R$,
one shows that, in fact, $\Kk(M)$ carries a canonical $A[F^e]$--structure.

The application for this general construction we have in mind is the
\idx{Matlis dual} functor $D=\Hom(\usc,E_R)$ for a local ring $R$ with
$E_R$ the injective hull of the residue field. We show that our formalism
applies to this situation by the existence of the natural transformation
$\F[e]_R\Hom(\usc,E_R) \to \Hom(\F[e]_R(\usc),E_R)$. We obtain an
extension $\Dd$ of the Matlis dual functor. The functor $\Dd$ specializes
to the functor $\Hh_{R,A}$ of Lyubeznik \cite[Chapter 4]{Lyub} for the
class of $R[F]$--modules which are cofinite $R$--modules supported on
$\Spec A$. On this class of $R[F^e]$--modules the functor $\Dd$ is
particularly useful. This follows from the fact that, if $M$ is such an
$R[F^e]$--module, then $\Dd(M)$ is a finitely generated unit
$R[F^e]$--module, provided that $R$ is a regular, complete and local ring
(cf.\ Proposition \ref{prop.DdcofiniteisFG}). Furthermore, the unique
minimal root of $\Dd(M)$ is obtained as the Matlis dual $D(\Ffred{M})$ of
the $F$--full and $F$--reduced subquotient $\Ffred{M}$ of $M$ (cf.\
Corollary \ref{cor.Dduniqueroot}). This implies that, roughly,  $\Dd$
gives a one-to-one correspondence between the $R[F^e]$--module quotients
of $M$ (up to $F$--full and $F$--reduced parts) on the one side, and the
unit $R[F^e]$--submodule of $\Dd(M)$. The most important example is that
$\Dd(H^d_m(A)) = H^c_I(R)$ which enables us to relate the $R[F^e]$--module
structure of $H^d_m(A)$ to the unit $R[F^e]$--module structure (and thus
the $D_R$-module structure) of $H^c_I(R)$.

\section{Adjoints of forgetful
functors}\label{sec.AdjofForget}\index{adjoints} The existence of a
\idx{right adjoint} to the forgetful functor from unit $R[F^e]$--modules
to $R[F^e]$--modules follows from the fact that this forgetful functor
commutes with direct limits. Then, after some set-theoretic issues are
taken care of, an adjoint always exists (cf.\ \cite[Chapter 5]{MacLane}).
But, in fact, we can concretely describe this right adjoint, first
introduced in \cite{HaSp} as follows: Let $(M,\theta^e)$ be an
$R[F^e]$--module. Define $G(M)$ as the inverse limit generated by the
structural map $\theta^e$; \ie
\[
    G(M) \defeq \invlim( \ldots \to \F[3e]M \to[\protect{\F[2e]\theta^e}] \F[2e]M \to[\protect{\F[e]\theta^e}] \F[e]M \to[\theta^e]
    M)
\]
and $G(M)$ comes equipped with natural maps $\pi_e: G(M) \to \F[e]M$. To
define the natural $R[F^e]$--module structure on $G(M)$ consider the maps
$\F[e]\pi_r$ from $\F[e]G(M) \to \F[e(r+1)]M$. These maps are compatible
with the maps of the limit system defining $G(M)$ and thus, by the
universal property of inverse limits, lift uniquely to a map
\[
    \F[e]G(M) \to G(M)
\]
defining the structure of an $R[F^e]$--module on $G(M)$.

\begin{proposition}\label{prop.Greg}
    Let\/ $R$ be regular and\/ $F$--finite and let\/ $M$ be an\/ $R[F^e]$--module. Then\/ $G(M)$ is
    a unit\/ $R[F^e]$--module.
\end{proposition}

\begin{proof}
    Since all we need to check is that a given map is an isomorphism we
    can localize and assume that $R$ is regular, local and therefore, that $R^e$ is
    a finitely generated free right $R$--module. Thus, tensoring
    with $R^e$ from the left (\ie applying $\F[e]$) commutes
    with the formation of inverse limits (elementary to check) and we get
    \[
        \F[e]G(M) \cong \invlim (\ldots \to \F[3e]M \to \F[2e]M \to \F[e]M) \cong
        G(M).
    \]
    This shows that the structural morphism of $G(M)$ is an
    isomorphism and thus $G(M)$ is unit.
\end{proof}
Even if $R$ is not regular, if $M$ is a unit $R[F^e]$--module, then $M
\cong G(M)$ as all the maps in the limit system defining $G(M)$ are then
also isomorphisms. For regular $F$--finite rings $R$, the natural map
$\pi_0: G(M) \to M$ is an isomorphism if and only if $M$ is a unit
$R[F]$--module. We get the following proposition.
\begin{proposition}\label{prop.RightAdjoint}
    Let\/ $R$ be regular and\/ $F$--finite. Then\/ $G$ is the
    right adjoint for the forgetful functor from\/ $R[F^e]$--modules
    to\/ $R[F^e]$--modules.
\end{proposition}
\begin{proof}
This essentially follows from the construction. After the last proposition
it remains to show that there is a functorial isomorphism
\[
    \Hom_{R[F^e]}(N,M) \cong \Hom_{R[F^e]}(N,G(M))
\]
for $M$ an $R[F^e]$--module and $N$ a unit $R[F^e]$--module. Given $\phi:
N \to{} M$, a map of $R[F^e]$--modules, it induces maps
$\phi_r=\F[er]\phi:\F[er]N \to \F[er]M$ for all $r$. These maps are
compatible with the $R[F^e]$--module structures and we therefore get a map
of inverse systems
\[
\xymatrix{
    {\ldots} \ar[r] &{\F[er]N} \ar[r] \ar_{\phi_r}[d] &{\F[e(r-1)]N} \ar[r] \ar_{\phi_{r-1}}[d] &{\ldots} \\
    {\ldots} \ar[r] &{\F[er]M} \ar[r] &{\F[e(r-1)]M} \ar[r] &{\ldots}
}
\]
which will give a map of the limits $\bar{\phi}: G(N) \to G(M)$. As $N$ is
unit, the natural map $G(N) \to N$ is an isomorphism so we can consider
$\bar{\phi}$ as a map from $N$ to $G(M)$.

Conversely, given a map $\psi:N \to G(M)$, we just compose with the
natural map $G(M) \to M$ to get the desired map $N \to M$.

The two processes just explained are easily verified to be functorial in
$N$ and $M$, and inverse to each other.
\end{proof}

\begin{examples}[\protect{cf.\ \cite[Proposition 2.1]{HaSp}}]
Let $R$ be regular and $F$--finite. For an ideal $I$ of $R$ the
quotient $A=R/I$ is an $R[F]$--module but generally not a unit
$R[F]$--module (cf.\ Example \ref{ex.FonIdeal}). $G(R/I)$ is the
limit of the sequence of surjections
\[
    \ldots \to \frac{R}{I^{[p^3]}} \to \frac{R}{I^{[p^2]}} \to \frac{R}{I^{[p]}}
    \to \frac{R}{I}
\]
which is just the completion of $R$ along the ideal $I$; \ie
$G(R/I) = \widehat{R}^I$ and the structural morphism $\theta: R^1
\tensor \widehat{R}^I \to \widehat{R}^I$ sending $r \tensor r'$ to
$r{r'}^p$ is an isomorphism.

If $R$ is local with maximal ideal $m$ we observed in Examples
\ref{ex.BasicRFmods} that the local cohomology modules $H^i_m(R/I)$
obtain, by functoriality, an $R[F^e]$--module structure from the
$R[F^e]$--module structure on $R/I$. It is the map $R^e \tensor H^c_I(R/I)
\cong H^i_m(R/I^{[p^e]}) \to H^i_m(R/I)$ that is induced by the natural
projection $R/I^{[p^e]} \to R/I$, which is the $R[F^e]$--module structure
on $R/I$. To determine the value of $G$ on these modules we calculate
\[
    G(H^i_m(R/I)) = \invlim( \ldots \to H^i_m(R/I^{[p^{2e}]}) \to
    H^i_m(R/I^{[p^{e}]}) \to H^i_m(R/I))
\]
which is just the local cohomology module $H^i_{m}(\widehat{R}^I)$ of the
formal completion of $R$ along $I$ (cf.\ \cite[Proposition
2.2]{Ogus.LocCohm}) with its canonical unit $R[F^e]$--structure induced
from the unit $R[F^e]$--structure on $\widehat{R}^I$ as just discussed,
\ie $G(H^i_m(R/I))=H^i_m(G(R/I))=H^i(\widehat{R}^I)$.
\end{examples}

\subsection{A left adjoint via Frobenius descent} As a nice application of
Frobenius descent we show a \idx{left adjoint} of the forgetful functor
from unit $D_R[F^e]$--modules to $D_R[F^e]$--modules as introduced in
\cite{Em.Kis2}. Let $(M,\theta^e)$ be an $R[F^e]$--module. Applying $T^e$
to the map $\theta^e$, we get a map $\beta^e : M \cong T^e(\F[e](M))
\to[T^e\theta^e] T^e(M)$ which, by further applying $T^e$, generates a
direct limit system. We define $U(M)$ to be its limit
\[
    U(M) \defeq \dirlim ( M \to T^eM \to T^{2e}M \to T^{3e}M \to \ldots).
\]
As $\F[e]$ commutes with direct limits, $\F[e](U(M))$ is the limit of the
system $\F[e]M \to M \to T^eM \to T^{2e}M \to \ldots$ which is isomorphic
to the one defining $U(M)$. Therefore $U(M)$ is naturally a unit
$R[F^e]$--module. By definition of $U(M)$ as a direct limit we have a
natural map $M \to U(M)$ which is an isomorphism if and only if $M$ is a
unit $R[F^e]$--modules. Functoriality of $U$ follows from the
functoriality of $T^e$ and direct limits.
\begin{proposition}[see also \cite{Em.Kis2}]
    Let\/ $R$ be regular and\/ $F$--finite. Then\/ $U$ is a left adjoint to the
    forgetful functor from unit\/ $R[F^e]$--modules to\/ $D_R[F^e]$--modules.
\end{proposition}
\begin{proof}
We have to show that there is a functorial isomorphism
\[
    \Hom_{R[F^e]}(U(M),N) \cong \Hom_{D_R[F^e]}(M,N)
\]
for unit $R[F^e]$--modules $N$ and $D_R[F^e]$--modules $M$. Given a map
$\phi: U(M) \to N$ we just precompose with the natural map $M \to U(M)$ to
obtain an element of the right hand side. Conversely, any map $\psi: M \to
N$ gives rise to a map $U(M) \to U(N) \cong N$ as $N$ is unit. To check
that these assignments are inverse to each other is straightforward.
\end{proof}

Emerton and Kisin show in \cite{Em.Kis2} the following interesting
application of this left adjoint; the proof of this is easily found to be
the same as the proof of our Proposition \ref{prop.DRfulluRF} and in fact
was inspired by their observation.
\begin{corollary}
    Let\/ $R$ be regular and\/ $F$--finite. Let\/ $M$ be a\/ $D_R[F^e]$--module such
    that the structural map\/ $\F[e]M \to[\theta^e] M$ is an injection and\/ $U(M)$ is a finitely
    generated unit\/ $R[F^e]$--module. Then\/ $M$ was already unit and\/ $M = U(M)$.
\end{corollary}
\begin{proof}
The injectivity of $\theta^e$ implies that $T^e(\theta^e)$ is also
injective. Therefore the natural map $M \to U(M) = \dirlim(M \to T^e(M)
\to T^{2e}(M) \to \ldots)$ is an injection. This makes $M$ a
$D_R[F^e]$--submodule of the finitely generated unit $D_R[F^e]$--module
$U(M)$. By Proposition \ref{prop.DRfulluRF} $M$ must be a unit
$D_R[F^e]$--module.
\end{proof}

\section{Extending Functors from $R$--mod to $R[F^e]$--mod}
The next task will be to develop a framework for deciding under what
circumstances it is possible to extend a functor from $R$--mod to $A$--mod
for rings $R$ and $A$ to a functor from $R[F^e]$--mod to $A[F^e]$--mod.
The case of a covariant functor we discussed in the introduction to this
chapter on page \pageref{x.CovFunct}. The contravariant case needs to be
treated with more detail. This will be our first task, but soon we
specialize to the Matlis dual functor which interests us most.

\subsection{Contravariant Functors on $R[F^e]$--modules}
Let $K$ be a contravariant functor from $R$--mod to $A$--mod. We have the
following proposition.
\begin{proposition}\label{prop.ContraExtends}
    Let\/ $K$ be a contravariant functor from\/ $R$--mod to\/ $A$--mod
    and let\/ $\Psi: \F[e] \circ K \to K \circ \F[e]$ be a natural
    transformation of functors. Then\/ $K$ naturally extends to a
    contravariant functor\/ $\cal{K}$ from\/ $R[F^e]$--mod to\/ $A[F^e]$--mod.
\end{proposition}
\begin{proof}
Given the natural transformation of functors $\Psi: \F[e] \circ K \to K
\circ \F[e]$ we can extend $K$ to a functor from $A[F^e]$--mod to
$R[F^e]$--mod. This extension $\Kk$ of $K$ depends on the natural
transformation $\Psi$. For an $R[F^e]$--module $M$ the construction of
$\Kk(M)$ goes as follows. We consider the inverse system generated by the
structural morphism $\theta^e: \F[e]M \to M$
\[
    \cdots \to \F[3e]M \to[\protect{\F[2e]}(\theta^e)] \F[2e](M) \to[\protect{\F[e]}(\theta^e)]
    \F[e]M \to[\theta^e] M
\]
to which we apply $K$ to obtain a direct system
\[
    K(M) \to[K(\theta^e)] K(\F[e](M)) \to[K(\protect{\F[e]}(\theta^e))] K(\F[2e](M))
    \to[] \cdots
\]
and $\cal{K}(M)$ is defined as its limit. The natural transformation
$\Psi$ allows us to endow $\cal{K}(M)$ with an $R[F^e]$--module structure.
As $\F[e]$ commutes with direct limits the map $\F[e]\cal{K}(M) \to
\cal{K}(M)$ is given by a map of directed systems
\begin{equation}\label{eqn.KkDef}
\begin{split}
\xymatrix{
    &{\F[e]K(M)} \ar[r] \ar[d]^{\Psi_M} & {\F[e]K\F[e](M)} \ar[r] \ar[d]^{\Psi_{\F[e]M}} & {\F[e]K\F[2e](M)} \ar[r] \ar[d]^{\Psi_{\F[2e]M}} & {\cdots} \\
{K(M)} \ar[r] &{K\F[e](M)} \ar[r] & {K\F[2e](M)} \ar[r] & {K\F[3e](M)} \ar[r] & {\cdots} \\
}
\end{split}
\end{equation}
The direct limit of the first row is $\F[e](\cal{K}(M))$ and the limit of
the second row is just $\cal{K}(M)$. By the naturality of $\Psi$ the
diagram commutes and we have obtained the structural morphism for
$\cal{K}(M)$. The fact that $\cal{K}$ is functorial follows from the
functoriality of $K$, $F^*$, $\dirlim$ and the naturality of $\Psi$.
\end{proof}

An important case is when the natural transformation $\Psi$ is an
isomorphism. If this is the case for all Frobenius powers of a given
$R[F^e]$--module $M$, then all the maps in Diagram \eqnref{eqn.KkDef} are
isomorphisms. Thus $\Kk(M)$ is, in fact, a unit $A[F^e]$--module. The fact
that $\Psi$ is an isomorphism also allows for an alternative way to
calculate $\Kk(M)$:
\begin{proposition}\label{prop.ContraExtendsWell}
    Let\/ $K$ be a contravariant functor together with the natural
    transformation\/
    $\Psi$ as before. If\/ $(M,\theta^e)$ is an\/ $R[F^e]$--module and\/ $\Psi_{\F[er]M}$
    is an isomorphism for all\/ $r \geq 0$, then\/ $\Kk(M)$ is the unit\/ $A[F^e]$--module
    generated by the map
    \[
        \beta^e: K(M) \to[K(\theta^e)] K(\F[e]M) \To[\cong]{\ \Psi_M^{-1}\ } \F[e]K(M).
    \]
    Thus, if\/ $\Psi$ is an isomorphism of functors, then\/ $\Kk$ is a functor
    from\/ $R[F^e]$--modules to unit\/ $A[F^e]$--modules.
\end{proposition}
\begin{proof}
One has to show that the direct limit system defining $\Kk(M)$ is
isomorphic to the one arising from the map $\beta^e: K(M) \to \F[e]K(M)$
by taking Frobenius powers (cf.\ the definition of generator on page
\pageref{x.DefGenerator}). This follows from the naturality of $\Psi$ and
from the assumption that $\Psi_{\F[er]M}: \F[e]K(\F[er]M) \to
K\F[e](\F[er]M)$ is an isomorphism for all $r \geq 0$ in a straightforward
fashion. We leave the details to the reader.
\end{proof}

Note that if $M$ is a unit $R[F^e]$--module, then $K(M) \cong \Kk(M)$
since the maps in the limit defining $\Kk(M)$ are all isomorphisms.
Therefore, if one only is interested in unit $R[F^e]$--modules $M$, then
$K(M)$ carries already a $A[F^e]$--module structure. In this case, $K(M)$
is also unit if and only if the fixed natural transformation $\Psi_M$ is
an isomorphism.

We want to be able to relate the functor $\Kk$ to the right adjoint of the
forgetful functor $G$ introduced in Section \ref{sec.AdjofForget}. This is
possible if we assume that $K$ commutes with limits: more precisely we
need an isomorphism of functors $K \circ \invlim = \dirlim \circ K$.
\begin{proposition}\label{prop.KkRelatesG}
    Let $R$ be regular and\/ $F$--finite and\/ $K$ contravariant with\/ $\Psi$
    as before. If we have an isomorphism of functors\/ $K \circ \invlim
    \cong \dirlim \circ K$, then\/ $\Kk \cong K \circ G = \Kk \circ G$.
\end{proposition}
\begin{proof}
The assumptions on $R$ assure that $G(M)$ is a unit $R[F^e]$--module by
Proposition \ref{prop.RightAdjoint}. By definition we have that
\[
    K(G(M)) = K( \invlim \F[er]M) = \dirlim K\F[er]M = \Kk(M).
\]
This finishes the argument.
\end{proof}

\subsection{\protect{$\Hom$} as a Functor on $R[F^e]$--modules}
We aim to apply the results from the preceding section to the
contravariant functor $\Hom(\usc,N)$ for some $R[F^e]$--module
$(N,\tau^e)$. For this observe that we have a natural transformation
$\Psi$ given for every $R$--module $M$ as the composition
\begin{equation}\label{eqn.NatTransForHom}
    \Psi: \F[e]\Hom(M,N) \to[\psi] \Hom(\F[e]M,\F[e]N) \to \Hom(\F[e]M,N)
\end{equation}
where the second map is just the one induced by $\tau^e$. The first map is
the natural transformation from $S \tensor \Hom_R(M,N) \to \Hom_S(S
\tensor M,S \tensor N)$ given by sending $s \tensor \phi$ to the map
$s\cdot (\id_S \tensor \phi)$ for a map of rings $R \to S$ (here $R=S$ and
the map is the Frobenius $F^e: R \to R$). In many good cases this
transformation is an isomorphism; we recall the following proposition.
\begin{proposition}\label{prop.HomAndTensor}
    Let\/ $f: R \to S$ be a flat map of noetherian rings. Then the natural map
    \[
        \psi: S \tensor \Hom_R(M,N) \to \Hom_S(S \tensor M,S \tensor N)
    \]
    is an isomorphism if\/ $S$ is module finite over\/ $R$ or if\/ $M$ is finitely
    presented.
\end{proposition}
\begin{proof}
The second part of this proposition is well known and follows easily by
applying $\Hom$ and $S \tensor \usc$ in either order to a presentation of
$M$ (see \cite[Lemma 3.2.4]{Trav.Phd}, for example). Then the flatness of
$S$ and the five lemma give the result. The first assertion is not quite
as simple and we will give a detailed argument.

Using geometric notation to denote restriction and extension of scalars by
$f_*$ and $f^*$ respectively we have to show that $f^*\Hom(M,N) \cong
\Hom(f^*M,f^*N)$  via the map $\psi$ with $\psi(s \tensor \phi)=s \cdot
f^*(\psi)$. Showing that this is an isomorphism can be done locally and
thus we assume that $S=\oplus_{i=1}^n e_iR$ is a finitely generated free
$R$--module. We define a new map $\psi'$ by going the bottom way through
the following diagram.
\[
\xymatrix{
    {f_*f^*\Hom_R(M,N)} \ar[r]^{\psi'} \ar[d]^{\cong} &{f_* \Hom_S(f^*M,f^*N)} \\
    {\oplus e_i\Hom_R(M,N)} \ar[d]^{\cong} &  \\
    {\Hom_R(M,\oplus e_iN)} \ar[r]^{\cong} & {\Hom_R(M,f_*f^*N)} \ar[uu]^{\cong}_{\text{adj}}
}
\]
Besides the adjointness of $f^*$ and $f_*$, which is responsible for the
right map, all other isomorphisms either come from the direct sum
decomposition of $S$ as an $R$ module or from $\Hom$ commuting with finite
direct sums in the second argument. Chasing through this definition, an
$R$--module generator $e_i \tensor \phi $ of the left hand side gets
mapped as follows:
\[
    e_i \tensor \phi \ \mapsto  \ e_i\phi \ \mapsto \ (e_i\cdot\usc ) \circ \phi \
    \mapsto \ (e_i \tensor\usc)\circ \phi \ \mapsto \ e_i \cdot f^*(\phi)
\]
Thus, by $R$--linearity of $\psi'$, we have $\psi'(s \tensor \phi) =
s\cdot f^*(\phi) = \psi(s \tensor \phi)$ for all $s \in S$ and $\phi \in
\Hom_R(M,N)$. Thus $\psi = \psi'$ which is an isomorphism by construction.
\end{proof}
If $R$ is regular and $F$--finite the conditions of the last proposition
are satisfied for the Frobenius map $F: R \to R$ and we get as a
corollary:
\begin{corollary}\label{cor.NatTransIsom}
    Let\/ $R$ be regular and\/ $F$--finite. Then the natural map
    \[
        \psi: \F[e]\Hom(M,N) \to \Hom(\F[e]M,\F[e]N)
    \]
    sending\/ $r \tensor \phi$ to $r\F[e](\phi)$ is an isomorphism for
    all\/
    $R$--modules\/ $M$ and\/ $N$.
\end{corollary}
This puts us in the situation were we can extend the functor
$\Hom(\usc,N)$ on $R$--modules to a functor on $R[F^e]$--modules provided
that $N$ is an $R[F^e]$--module.
\begin{proposition}\label{prop.HomExtendsToRF}
    Let $R$ be noetherian and\/ $N$ an\/ $R[F^e]$--module. Then\/
    $\Hom(\usc,N)$ naturally extends to a functor\/ $\Hhom(\usc,N)$ from\/
    $R[F^e]$--modules to\/ $R[F^e]$--modules.

    If $R$ is regular and $F$--finite and $N$ is unit, then
    $\Hhom(M,N)\cong\Hom(G(M),N)$ is in fact a unit $R[F^e]$--module for every
    $R[F^e]$--module $M$.
\end{proposition}
\begin{proof}
For the first part we use the natural transformation $\Psi$ described in
\eqnref{eqn.NatTransForHom} to apply Proposition \ref{prop.ContraExtends}.
This shows that we can extend $\Hom(\usc,N)$ naturally to a functor on
$R[F^e]$--modules.

Under the assumptions of the second part we see that the transformation
$\Psi$ is an isomorphism by Corollary \ref{cor.NatTransIsom} and the fact
that $N$ is assumed to be unit. Furthermore, it is well known that
$\Hom(\usc,N)$ interchanges direct and inverse limits and therefore
everything follows from Proposition \ref{prop.ContraExtendsWell} and
Proposition \ref{prop.KkRelatesG}.
\end{proof}

\section{Lyubeznik's $\Hh_{R,A}$ and Dualization}\label{sec.LyubandDd}
We are ready for the main task of this chapter: extending the Matlis
duality functor \index{$D$}$D=\Hom(\usc,E_R)$ to a functor from
$R[F^e]$--modules to $R[F^e]$--modules. To apply the general results of
the last section on extending $\Hom$ functors we need a natural
$R[F^e]$--module structure on $E_R$, the injective hull of the residue
field of $R$. If $R$ is regular, then we saw that all local cohomology
modules carry a natural unit $R[F^e]$--module structure. For a regular
local $(R,m)$ the injective hull $E_R$ of the residue field of $R$ is
isomorphic to the top local cohomology module $H^n_m(R)$ of $R$ with
support in $m$. The choice of such an isomorphism endows $E_R$ with a unit
$R[F^e]$--module structure $\theta^e: \F[e]E_R \to E_R$ which we will now
fix. By the results of the preceding section (in particular Proposition
\ref{prop.HomExtendsToRF}), the unit $R[F^e]$--module structure on $E_R$
allows us to extend the Matlis dual functor $D = \Hom(\usc,E_R)$ to a
functor \index{$\Dd$}$\Dd$ from $R[F^e]$--modules to $R[F^e]$--modules.
Concretely, for an $R[F^e]$--module $(M,\theta^e)$ we have
\[
    \Dd(M)= \dirlim(D(M) \to[D(\theta^e)] D(\F[e]M) \to[D(\protect{\F[e]}(\theta))]
    D(\F[2e]M) \to \ldots\ ).
\]
If $R$ is also $F$--finite, then the second part of Proposition
\ref{prop.HomExtendsToRF} shows that $\Dd(M)$ is a unit $R[F^e]$--module.
In order to not be restricted to the case that $R$ is $F$--finite we show
the following lemma which should be viewed as an extension of Proposition
\ref{prop.HomAndTensor} for $\Hom(\usc,E_R)$.
\begin{lemma}\label{lem.FcommutesD}
    Let\/ $R$ be regular and local. The natural map\/ $\psi_M:
    \F[e]\Hom(M,E_R) \to \Hom(\F[e]M,\F[e]E_R)$ is an isomorphism for all
    cofinite\/ $R$--modules\/ $M$.

    Since\/ $E_R$ is a unit\/ $R[F^e]$--module this implies that we have a
    natural isomorphism\/ $\F[e](D(M)) \to D(\F[e](M))$ for cofinite\/
    $R[F^e]$--modules\/ $M$.
\end{lemma}
\begin{proof}
Let $M$ be a cofinite $R$--module. We take a resolution $0 \to M \to E_1
\to E_2$ where $E_i$ is a finite direct sum of copies of $E_R$. By
applying $\F[e]$ and $\Hom$ in either order to this resolution we reduce
the task of showing that $\psi_M$ is an isomorphism to the case $M = E_i$.
Since $\Hom$ and $\F[e]$ both commute with finite direct sums we further
reduce to the case $M=E_R$. Thus, it remains to show that $\psi: \F[e]
\Hom(E_R,E_R) \to \Hom(\F[e]E_R,\F[e]E_R)$ is an isomorphism. This is
clear after one observes that both sides are canonically isomorphic to
$\widehat{R}$.
\end{proof}
In \cite[Section 4]{Lyub}, this lemma is the starting point for defining a
functor $\Hh_{R,A}$ for $A = R/I$ and $R$ complete, regular and local.
Lyubeznik's $\Hh_{R,A}$ is a functor from cofinite $A[F]$--modules to unit
$R[F]$--modules and for such an $A[F]$--module $(M,\theta)$ it is defined
as the unit $R[F]$--module generated by the map
\[
    \beta: D(M) \to[D(\theta)] D(\F(M)) \to[\Psi^{-1}] \F(D(M)).
\]
By our Proposition \ref{prop.ContraExtendsWell} this functor is therefore
just the restriction of $\Dd$ to the cofinite $R[F^e]$--modules supported
on $\Spec A$.
\begin{proposition}\label{prop.MalisDualonRF}
    For a regular, local ring\/ $R$, the Matlis dual functor\/ $D=\Hom(\usc,E_R)$
    naturally extends to a functor\/ $\Dd$ from\/ $R[F^e]$--modules to\/
    $R[F^e]$--modules.

    If\/ $(M,\theta^e)$ is cofinite or finitely generated as an\/ $R$--module
    or if\/ $R$ is\/ $F$--finite, then\/ $\Dd(M)$ is the unit\/
    $R[F^e]$--module generated by the map
    \[
        \beta^e: D(M) \to[D(\theta^e)] D(\F[e](M)) \to[\Psi] \F[e](D(M)).
    \]
    Furthermore,\/ $\Dd$ is an exact functor for this class of modules
    (respectively, for this class of rings).
\end{proposition}
\begin{proof}
The first part is just an application of the abstract machinery developed
in the previous section, in particular Proposition
\ref{prop.HomExtendsToRF}.

In the second part the assumptions on $M$ or $R$ ensure that by Lemma
\ref{lem.FcommutesD} or Proposition \ref{prop.HomAndTensor} the natural
transformation $\F[e]\Hom(\usc,E_R) \to \Hom(\F[e]\usc,\F[e]E_R)$ is an
isomorphism. Together with a fixed unit $R[F^e]$--structure on $E_R$ this
shows that the natural transformation $\Psi: \F[e] \circ D \to D \circ
\F[e]$ is an isomorphism. As for a finitely generated (respectively
cofinite) $R$--module $M$ all modules $\F[er]M$ are also finitely
generated (respectively cofinite) an application of Proposition
\ref{prop.ContraExtendsWell} and Proposition \ref{prop.KkRelatesG} yields
all that is claimed. The exactness follows from the exactness of $D$,
$\F[e]$ and $\dirlim$ with the description of $\Dd$ as the unit
$R[F^e]$--module generated by $\beta^e$.
\end{proof}
\begin{remark}
    For the first part of the above proposition it is enough to assume $R$
    Gorenstein since all one needs is that $E_R$ carries an $R[F^e]$--module
    structure. If $R$ is only Gorenstein there is an
    almost canonical choice of such structure via an isomorphism
    of $E_R$ with the top local cohomology module $H^d_m(R)$.
    As we saw in Example \ref{ex.BasicRFmods} $H^d_m(R)$ carries a natural
    unit $R[F^e]$--structure even if $R$ is not regular.
\end{remark}
\begin{remark}\label{rem.DforRFinfty}
We have to comment on the compatibility of $\Dd$ with different powers of
the Frobenius. In the definition of $\Dd$ we used a certain unspecified
power $e$ of the Frobenius and thus \emph{a priori} defined Functors
$\Dd^e$ for all $e \geq 0$. For various $e$ these Functors $\Dd^e$ are
compatible with the inclusion of categories $\Rmod[\protect{R[F^{re}]}]
\subseteq \Rmod[\protect{R[F^e]}]$. To see this, let $(M,\theta^e)$ be an
$R[F^e]$--module. Then, viewed as an $R[F^{er}]$--module $(M,\theta^{er})$
we calculate $\Dd^{er}(M)$ as the direct limit of
\[
    D(M) \to D(\F[er](M)) \to D(\F[2er](M)) \to \ldots.
\]
This is just a subsystem of the directed system defining $\Dd^e(M)$. Thus
the limits coincide. To see that $\Dd^{e}(M)$ and $\Dd^{er}(M)$ in fact
carry the same $R[F^{er}]$--structure we only remark that the natural
transformations $\Psi^e :\F[e] \circ D \to D \circ \F[e]$ for all $e$ are
induced from the natural transformation $\F \circ D \to D \circ \F$, \ie
they originate from the $e=1$ case. With this in mind it follows that the
collection of functors $\{\Dd^e\}$ naturally define a functor $\Dd^\infty$
on $R[F^\infty]$--modules. Given these compatibilities, we abuse notation
and denote all these functors just by $\Dd$.
\end{remark}

\subsection{Properties of $\Dd$}\label{sec.BasicsOfDd}
We want to transfer certain aspects of Matlis duality to incorporate
Frobenius actions. Especially the property that for a complete local ring
$R$ the Matlis dual of a cofinite $R$--module is finitely generated. As
the functor $\Dd$ is one of our key tools in the next chapter, we describe
here some of its properties in detail. Much of what follows generalizes
the properties of $\Hh_{R,S}$ noted in \cite[Section 4]{Lyub}. As Matlis
Duality works best for complete rings we will stick to this case in this
section and therefore assume that $R$ be complete until further notice.

\begin{proposition}\label{prop.DdcofiniteisFG}
    Let\/ $R$ be complete and regular. On the subcategory of\/ $R[F^e]$--modules
    which are cofinite as $R$--modules $\Dd$ is exact and its values
    are finitely generated unit\/ $R[F^e]$--modules.
\end{proposition}
\begin{proof}
With the notation of Proposition \ref{prop.MalisDualonRF} a generator of
$\Dd(M)$ is given by $\beta^e : D(M) \to \F[e]D(M)$. Since $R$ is complete
$D(M)$ is a finitely generated $R$--module. Consequently $\Dd(M)$ is a
finitely generated unit $R[F^e]$--module.
\end{proof}
Another important property of the Matlis dual is that for finite or
cofinite $R$--modules $M$ one has $D(D(M))\cong M$. We cannot expect this
to hold for $\Dd$ since in any case $\Dd(\Dd(M))$ is unit even if $M$
wasn't by Proposition \ref{prop.MalisDualonRF}. An obvious guess for
$\Dd(\Dd(M))$ is that it is equal to $G(M)$. To see this let $M$ be
finitely generated or cofinite as an $R$--module. Then $\Dd(M)$ is a unit
$R[F^e]$--module and therefore $\Dd(\Dd(M)) \cong D(\Dd(M))$. Now use that
$D$ transforms direct to inverse limits and get
\begin{equation*}
\begin{split}
    D(\Dd(M)) &= D(\dirlim(D(M) \to D(\F[e]M) \to D(\F[2e]M) \to \ldots \
                 ) \\
              &\cong \invlim(\ldots \to D(D(\F[2e]M)) \to D(D(\F[e]M)) \to
                 D(D(M)) ) \\
              &\cong G(M)
\end{split}
\end{equation*}
where we used that for finitely generated (resp. cofinite) $M$ also
$\F[er]M$ is finitely generated (resp. cofinite) as an $R$--module and
thus $D(D(\F[er]M))\cong \F[er]M$. We get the following proposition.
\begin{proposition}
Let\/ $R$ be regular and complete. Let\/ $M$ be an\/ $R[F^e]$--module that
is finitely generated or cofinite as an\/ $R$--module. Then\/ $\Dd(\Dd(M))
\cong G(M)$.
\end{proposition}

Before proceeding we introduce some notation loosely motivated by the
notation in \cite{HaSp}. Let $(M,\theta^e,F^e)$ be an $R[F^e]$--module. An
element $m \in M$ is called $F$--nilpotent if $F^{re}(m)=0$ for some $r$.
Then $M$ is called \emph{$F$--nilpotent}\index{F-nilpotent@$F$--nilpotent}
if $F^{er}(M)=0$ for some $r \geq 0$. It is possible that every element of
$M$ is $F$--nilpotent but $M$ itself is not, since $F$--nilpotency for $M$
requires that all $m \in M$ are killed by the \emph{same} power of $F^e$.
In particular the sub $F[R^e]$--module consisting of all $F$--nilpotent
elements $\Fnil{M}$ need not be nilpotent in general. If $\theta^e$ is
surjective, then $M$ is called \emph{F--full}\index{F-full@$F$--full}.
Note that $F$--fullness does not mean $F^e$ is surjective but merely that
the submodule $F^e(M)=\theta^e(\F[e]M)$ is all of $M$. Finally we say that
$M$ is \emph{$F$--reduced}\index{F-reduced@$F$--reduced} if $F^e$ acts
injectively.

\begin{remark}
If $R$ is $F$--pure (\ie the Frobenius on $R$ is a pure map of rings),
then $F$--reducedness of $M$ is implied by $\theta^e$ being injective.
This follows since, by definition of $F$--purity, the map $F^e_R \tensor
\id_M$ is injective. Therefore, $F^e_M = \theta^e_M\circ(F^e_R \tensor
\id_M)$ is injective if $\theta^e_M$ is injective. This is false in
general. As an example take the top local cohomology module $H^d_m(R)$ of
any non $F$--injective\index{F-injective@$F$--injective} ring, \ie a ring
for which the Frobenius does not act injectively on the top local
cohomology module. $H^d_m(R)$ is a unit $R[F^e]$--module and thus
$\theta^e$ is injective. Such rings exist and a concrete example is
$\frac{k[x,y,z]}{x^4+y^4+z^4}$.
\end{remark}

The above notions are the same if we view $M$ as an $R[F^{er}]$--module
for some $r \geq 0$. Therefore they also apply to $R[F^\infty]$--modules
and thus we allow $e=\infty$ in what follows.

We are lead to some functorial constructions for $R[F^e]$--modules. The
$R[F^e]$--submodule consisting of all $F$--nilpotent elements of $M$ we
denote by \index{$\Fnil{M}$}$\Fnil{M}=\{\,m \in M\,|\,F^{er}(m)=0 \text{
for some } r\,\}$. The quotient $M/\Fnil{M}$ is the biggest $F$--reduced
quotient and denoted by \index{$\Fred{M}$}$\Fred{M}$. The
$R[F^e]$--submodule \index{$\Fful{M}$}$\Fful{M}=\bigcap F^{er}(M)$ is the
largest $F$--full submodule. If $M$ is a cofinite $R$--module, then the
decreasing chain of $R[F^e]$--submodules $F^{er}(M)$ stabilizes and we
have $\Fful{M}=F^{er}(M)$ for some $r > 0$. We note some properties in the
following lemma.
\begin{lemma}\begin{enumerate}
\item
    The operation of taking\/ $F$--nilpotent parts is left exact. If $N \subseteq M$, then
    $\Fred{N} \subseteq \Fred{M}$.
\item
    Submodules of\/ $F$--reduced\/ $R[F^e]$--modules are also\/ $F$--reduced. Quotients
    of\/ $F$--full\/ $R[F^e]$--modules are\/ $F$--full. The property of\/ $F$--nilpotency
    passes to quotients and submodules.
\item
    $M/F^{er}(M)$ is\/ $F$--nilpotent for all\/ $r$.
\item
    The operations\/ $\Fful{(\usc)}$
    and\/ $\Fred{(\usc)}$ mutually commute which makes the\/ $F$--full and
    $F$--reduced subquotient\/ $\Ffred{M}=\Fred{(\Fful{M})}=
    \Fful{(\Fred{M})}$ of an\/ $R[F^e]$--module\/ $M$ well defined.
\end{enumerate}
\end{lemma}
\begin{proof}
For (a), note that if $N \subseteq M$ then $\Fnil{N} = \Fnil{M} \cap N$.
This, together with the fact that $\Fnil{M}$ is a submodule of $M$ implies
that $\Fnil{(\usc)}$ is left exact. The same formula $\Fnil{M} = \Fnil{M}
\cap N$ also implies that $\Fred{N} = N/\Fnil{N}$ is a submodule of
$\Fred{M}$. Note that $\Fred{(\usc)}$ is not left exact, in general.

(b) and (c) are clear, but we point out that a quotient of an $F$--reduced
module might not be $F$--reduced, and similarly, submodules of $F$--full
$R[F^e]$--modules might not be $F$--full.

To show part (d) we observe that $\Fful{M} \onto \Fful{(\Fred{M})}$ is a
surjection with kernel $\Fful{M} \cap \Fnil{M}$. This follows since an
element $m \in \Fful{M}$ is mapped to zero in $\Fful{(\Fred{M})}$ if and
only if the image of $m$ in $\Fred{M}$ is zero (the full parts $\Fful{M}$
are submodules of $M$!). This is the case if and only if $m \in \Fnil{M}$.
On the other hand $\Fful{M} \cap \Fnil{M} = \Fnil{(\Fful{M})}$ and by
definition $\Fred{(\Fful{M})}$ is the cokernel of the inclusion
$\Fnil{(\Fful{M})} \into \Fful{M}$. Thus we conclude that, in fact,
$\Fful{(\Fred{M})}\cong\Fred{(\Fful{M})}$.
\end{proof}

The relevance of these notions for the study of $\Dd$ is demonstrated by
the following proposition.
\begin{proposition}\label{prop.DisZeroOnNilpotent}
    Let\/ $R$ be a complete, regular local ring. Let\/ $M$\ be an\/ $R[F^e]$--module
    which is a cofinite\/ $R$--module. Then $M$ is\/ $F$--nilpotent if and only
    if\/ $\Dd(M)=0$.
\end{proposition}
\begin{proof}
Note that $(M, \theta^e)$ is $F$--nilpotent if and only if for some $r
\geq 0$ we have $F^{er}(M)=0$. By Remark \ref{rem.DforRFinfty} we can
replace $e$ by $er$ and assume that $F^e(M)=0$, which is equivalent to
$\theta^e$ being the zero map. This implies that the generator $\beta^e =
\Psi_{M} \circ D(\theta^e)$ is also the zero map, and thus the module
$\Dd(M)$ generated by $\beta^e$ is also zero.

Conversely, $\Dd(M)$ is the unit $R[F^e]$--module generated by $\beta^e:
D(M) \to \F[e](D(M))$ by Proposition \ref{prop.MalisDualonRF}. Thus,
$\Dd(M)$ is zero if and only if the image of $D(M)$ in $\Dd(M)$ is zero.
This, in turn, is equivalent to $D(M) = \bigcup \ker \beta^{er}$ where
$\beta^{er}$ is defined inductively by
$\beta^{e(r+1)}=\beta^{er}\circ\F[er](\beta^e)$. Since $D(M)$ is finitely
generated this increasing union stabilizes. Therefore, $\beta^{er}=0$ for
some $r>0$. Up to the natural transformation $\Psi: \F[er] \circ D \cong D
\circ \F[er]$ the map $\beta^{er}$ is just $D(\theta^{er})$ and we
conclude that also $\theta^{er}=0$. This implies that $F^{er}(M)=0$, \ie
$M$ is $F$--nilpotent.
\end{proof}
Hartshorne and Speiser show in \cite[Proposition 1.11]{HaSp} that if $M$
is an $R[F^e]$--module which is cofinite as an $R$--module, and $M =
\Fnil{M}$, then, in fact, $M$ is $F$--nilpotent, \ie a uniform power of
$F$ kills all elements of $M$. Thus, \emph{a posteriori}, in the above
proposition one could replace the condition that $M$ be $F$--nilpotent
with the weaker condition that $M = \Fnil{M}$. It would be interesting to
know if this modified version is true for not necessarily cofinite
$R[F^e]$--modules. We ask: Is it true for all $R[F^e]$--modules that
$\Dd(M)=0$ if and only if $M = \Fnil{M}$?

%In our treatment of the properties of $\Dd$ on cofinite $R[F^e]$--modules
%we will now use this fact from \cite{HaSp} that $\Fnil{M}$ is
%$F$--nilpotent for a cofinite $R[F^e]$--module $(M,\theta)$. This differs
%from the treatment in \cite{Lyub} where he shows the properties of
%$\Hh_{R,S}$ (the equivalent of our $\Dd$ in his notation) while proving
%the statement of \cite[Proposition 1.11]{HaSp} on the fly (it really is
%part (b) and (c) of \cite[Lemma 4.3]{Lyub} or then explicitly Proposition
%4.4).

We saw that $F$--nilpotency of $(M,\theta^e)$ forces the map $\beta^{er}$
to be zero for some $r>0$, and thus one concludes that $\Dd(M)=0$. Now we
investigate when $\beta^e$ is injective. If $M$ is cofinite, then $D(M)$
is finitely generated and the injectivity of $\beta^e$ then means nothing
but that $D(M)$ is a root of $\Dd(M)$. Up to the natural transformation
$\Psi$ the map $\beta^e$ is just $D(\theta^e)$. Thus $\beta^e$ is
injective if and only if $\theta^e$ is surjective which, by definition, is
the case if and only if $M$ is $F$--full. Thus we get:
\begin{proposition}\label{prop.Full=Root}
    Let\/ $R$ be regular, complete and local. Let\/ $(M,\theta^e)$ be an\/
    $R[F^e]$--module which is a cofinite\/ $R$--module. Then\/ $\beta^e: D(M) \to
    \F[e]D(M)$ is a root of\/ $\Dd(M)$ if and only if\/ $M$ is\/ $F$--full.
\end{proposition}
Now, assuming that $M$ is cofinite and $F$--full we want to connect the
roots of $\Dd(M)$ with the $R[F^e]$--quotients of $M$. For this we recall
the following result.
\begin{lemma}[\protect{\cite[Lemma 4.3]{Lyub}}]\label{lem.DandRoots}
    Let\/ $R$ be regular, complete and local. Let\/ $(M,\theta^e)$ be a
    cofinite and\/
    $F$--full\/ $R[F^e]$--module.\/ $D(N) \subseteq D(M)$ is a root of a
    unit\/
    $R[F^e]$--submodule of\/ $\Dd(M)$ if and only if\/ $N$ is a
    $R[F^e]$--module quotient of\/ $M$. $D(N)$ is a root of\/ $\Dd(N)$ itself
    if and only if\/ $\ker(M \to N)$ is\/ $F$--nilpotent.
\end{lemma}
\begin{proof}
$D(N)$ is a root of some unit $R[F^e]$--submodule if and only if $N
\subseteq F^e_{\Dd(N)}(N)$. By definition of the $R[F^e]$--structure on
$\Dd(M)$, this is equivalent to $\beta^e(D(N))\subseteq\F[e]D(N)$. By the
naturality of the transformation $\Psi: D \circ \F[e] \to \F[e] \circ D$
this is equivalent to $D(\theta^e)(D(N)) \subseteq D(\F[e]N)$. Dualising
this last inclusion we get a surjection $\theta^e: \F[e]N \to N$ (here we
use that $R$ is complete) showing that $N$ is an $R[F^e]$--module quotient
of $M$.

By exactness of $\Dd$ we have $\Dd(M) = \Dd(N)$ if and only if $\Dd(\ker(M
\to N)=0$. By Proposition \ref{prop.DisZeroOnNilpotent} this is the case
if and only if $\ker(M \to N)$ is $F$--nilpotent.
\end{proof}
\begin{proposition}\label{prop.Dd=ifFfred=}
    Let\/ $R$ be regular, complete and local. Let\/ $M$ and\/
    $N$ be two\/ $R[F^e]$--modules which are cofinite as\/ $R$--modules.
    Then\/
    $\Ffred{M} \cong \Ffred{N}$ if and only if\/ $\Dd(M) \cong \Dd(N)$.
\end{proposition}
\begin{proof}
First, we observe that $\Dd(M)\cong\Dd(\Fful{M})$: Since $M$ is cofinite
the decreasing chain of submodules $F^{er}(M)$ stabilizes at its limit in
finitely many steps; \ie for some $r > 0$ we have $F^{er}(M) = \Fful{M}$.
Thus the quotient $M/\Fful{M}=M/F^{er}M$ is killed by $F^{er}$ and thus is
$F$--nilpotent. Therefore, $\Dd(M/\Fful{M})=0$ and by exactness of $\Dd$
applied to the sequence
\[
    0 \to \Fful{M} \to M \to M/\Fful{M} \to 0
\]
it follows that $\Dd(M) \cong \Dd(\Fful{M})$.

Secondly, one shows similarly that $\Dd(\Fred{M})\cong\Dd(M)$. With
Proposition 1.11 of \cite{HaSp} (or \cite[Lemma 4.3]{Lyub}) we find
$\Fnil{M}$ to be nilpotent and therefore $\Dd(\Fnil{M})=0$ by Proposition
\ref{prop.DisZeroOnNilpotent}. Applying $\Dd$ to the exact sequence $0\to
\Fnil{M} \to M \to \Fred{M} \to 0$ gives the desired isomorphism
$\Dd(\Fred{M}) \cong \Dd(M)$. Putting these two observation together
completes one direction of the equivalence.

For the converse we can assume that $M$ and $N$ are $F$--reduced and
$F$--full (\ie $F$ injective and $\theta^e$ is surjective). Thus $D(M)$ is
a root of $\Dd(M)$ by Proposition \ref{prop.Full=Root}. Since $R$ is
complete $\Dd(M)$ has a unique minimal root $M_0 \subseteq D(M)$ by
Proposition \ref{prop.MinimalRoot}. By completeness of $R$ and Matlis
duality we can write $M_0=D(M')$ for a quotient $M'$ of $M$ where the fact
that $M_0$ is a root of $\Dd(M)$ (\ie $M_0 \subseteq F^*M_0$) ensures that
$M'$ is an $R[F^e]$--module quotient of $M$. Since $\Dd(M) = \Dd(M')$
Lemma \ref{lem.DandRoots} shows that the kernel of $M \onto M'$ must be
$F$--nilpotent. Since $M$ is assumed $F$--reduced this kernel must be
zero. Therefore $M_0=D(M')=D(M)$. Thus $D(M)$ is the unique minimal root
of $\Dd(M)$. The same argument holds for $D(N)$ which shows that
$D(N)\cong D(M)$ and therefore $N \cong M$ as required.
\end{proof}
As a corollary of the proof of this we get:
\begin{corollary}\label{cor.Dduniqueroot}
    Let\/ $R$ be complete, regular and local. Let\/ $M$ be a
    cofinite\/ $R[F^e]$--module. Then\/ $D(\Ffred{M})$ is the unique minimal root
    of\/ $\Dd(M)$.
\end{corollary}
With the same ideas as in this proof we are able to show our main
characterization of $\Dd$:
\begin{theorem}\label{thm.DdofSimpisSimp}
    Let\/ $R$ be complete, regular and local. Let\/ $(M, \theta^e)$ be a
    cofinite\/
    $R[F^e]$--module. Every unit\/ $R[F^e]$--submodule of\/ $\Dd(M)$ arises
    as\/
    $\Dd(N)$ for some\/ $R[F^e]$--module quotient of\/ $\Ffred{M}$.

    In fact,\/ $\Dd$ is an isomorphism between the lattice of\/
    $R[F^e]$--module quotients of\/ $M$ (up to\/ $\Ffred{(\usc)}$) and the
    unit\/
    $R[F^e]$--submodules of\/ $M$. Consequently\/ $\Dd(M)$ is
    a simple unit\/ $R[F^e]$--module if and only if\/ $\Ffred{M}$ is a
    simple\/
    $R[F^e]$--module.
\end{theorem}
\begin{proof}
We can assume that $M$ is $F$--full and $F$--reduced by the last
proposition. Then $D(M) \to[\beta^e] \F[e]D(M)$ is the unique minimal root
of $\Dd(M)$. Let $\Nn$ be a unit $R[F^e]$--submodule of $\Dd(M)$, then
$N_0
\defeq \Nn \cap D(M)$ is a root of $\Nn$. In particular $\beta^e$
restricts to an injection $N_0 \to \F[e](N_0)$ which, in turn, is
equivalent to $\theta^e$ inducing a surjection $\F[e]D(N_0) \to D(N_0)$.
But this says nothing but that $D(N_0)$ is an $R[F^e]$--module quotient of
$M$, obviously $\Dd(D(N_0)) = \Nn$. In this argument we used, without
mention, the natural isomorphism of functors $D \circ \F[e] \cong \F[e]
\circ D$ and that up to this identification $\beta^e = D(\theta^e)$.

The second part now follows trivially from the first part together with
the last Proposition saying that $\Dd(N) \cong \Dd(N')$ if and only if
$\Ffred{N} \cong \Ffred{N'}$ for $R[F^e]$--quotients $N$ and $N'$ of $M$.
\end{proof}
As a final remark we point out that if $M$ is a simple $F$--full
$R[F^e]$--module, then $\Dd(M)$ is nonzero and therefore a simple unit
$R[F^e]$--module. This follows since a simple $R[F^e]$--module is
$F$--full if and only if it is $F$--reduced and therefore by the last
Proposition $\Dd(M)$ is simple (and automatically nonzero). If $F^e$ had a
kernel it would be a nontrivial $R[F^e]$--submodule and thus if $M$ is
simple the kernel of $F^e$ must be all of $M$. Thus $F^e(M)=0$ which
contradicts the $F$--fullness since this exactly means that $F^e(M)=M$.

\subsection{The main example: $H^d_m(A)$}\label{sec.HdRelateHi}
We want to apply our extension of the Matlis dual to the top local
cohomology module of a quotient $(A,m)$ of a complete regular local ring
$(R,m)$. If $\pi: R \to A$ is the quotient map and $I$ is its kernel we
denote by $n$ the dimension of $R$. The local cohomology module $H^i_m(A)$
is a $A$--module and, by restriction, also an $R$--module. As it is a
cofinite $A$--module it is also cofinite as an $R$--module. As we
discussed earlier (cf.\ Example \ref{ex.BasicRFmods}), the natural unit
$A[F]$--structure on $A$ induces a unit $A[F]$--structure $\theta_A$ on
$H^i_m(A)$. Similarly, the natural $R[F]$--structure on $A=R/I$ (cf.\
Example \ref{ex.FonIdeal}) induces an $R[F]$--structure $\theta_R$ on the
$R$--module $H^i_m(R/I)$. This $R[F]$--structure is not unit since the
$R[F]$--structure on $A$ was not unit either. The following diagram shows
the connections between these structures:
\[
\xymatrix@C=4pc{
    {R^1 \tensor_R H^i_m(R/I)} \ar^{\theta_R}[rd]\ar_{\pi \tensor \id}[d] \ar^{\cong}[r]  &  {H^i_m(R/I^{[p]})} \ar[d] \\
    {A^1 \tensor_A H^i_m(A)} \ar_{\theta_A}[r] & {H^i_m(R/I)}
}
\]
The right vertical map is induced from the natural surjection $R/I^{[p]}
\to R/I$ which can be identified with the natural $R[F]$--module structure
on $A=R/I$. Therefore, by definition $\Dd(H^i_m(R/I))$ is the direct limit
of
\begin{equation}
    D(H^i_m(R/I)) \to D(H^i_m(R/I^{[p]})) \to D(H^i_m(R/I^{[p^2]})) \to
    \ldots
\end{equation}
where we use the natural isomorphism $R^e \tensor H^i_m(R/I) \cong
H^i_m(R/I^{[p^e]})$ which follows from the flatness of $R^e$ as a right
$R$--module. As in the case $r=1$ the maps are the ones induced from the
surjections $R/I^{[p^{r+1}]} \to R/I^{[p^{r}]}$. Using local duality
\cite[Theorem 3.5.8]{BrunsHerzog} for the complete, regular and local ring
$R$ this directed sequence is isomorphic to the following
\begin{equation}\label{eqn.DdfromExt}
    \Ext_R^{n-i}(R/I,R) \to \Ext_R^{n-i}(R/I^{[p]},R) \to \Ext_R^{n-i}(R/I^{[p^2]},R)
    \to \ldots
\end{equation}
where, again, the maps are the ones induced from the natural projections.
Since the Frobenius powers of an ideal are cofinite within the normal
powers, we get that the limit of this sequence is just $H^{n-i}_I(R)$.
This is because an alternative definition of $H^{n-i}_I(R)$ is as the
right derived functor of the functor $\Gamma_I(M)=\dirlim \Hom(R/I^{t},R)$
of sections with support in $\Spec R/I$ \footnote{see \cite[Theorem
3.5.6]{BrunsHerzog} for the equivalence with our definition of local
cohomology via \Czech complexes}. Then using that $\dirlim$ is an exact
functor one gets
\begin{equation*}
\begin{split}
    H^{n-i}_I(R) &= \mathbf{R}^{n-i}(\dirlim \Hom(R/I^{t},R)) = \dirlim \mathbf{R}^{n-i}\Hom(R/I^{t},R) \\
                 &= \dirlim \Ext_R^{n-i}(R/I^t,R)
\end{split}
\end{equation*}
by definition of $\Ext$ as the right derived functor of $\Hom$. Thus in
fact $\Dd(H^i(R/I)) \cong H^{n-i}_I(R)$. It remains to show that under
this identification the respective unit $R[F]$--module structures are
preserved, \ie that $\Dd(H^i_m(R/I))$ and $H^{n-i}_I(R)$ are in fact
isomorphic as $R[F^e]$--modules. But this is an immediate consequence of
\cite[Propositions 1.8 and 1.11]{Lyub}. Summarizing we get.
\begin{proposition}\label{prop.DdofHmisHi}
    Let\/ $(R,m)$ be regular, local, complete and\/ $F$--finite. Let\/ $A=R/I$
    for some ideal\/ $I$ of\/ $R$. Then
    \[
        \Dd(H^i_m(R/I)) \cong H^{n-i}_I(R)
    \]
    as unit\/ $R[F]$--modules.
\end{proposition}
Using Proposition \ref{prop.MalisDualonRF} we know that a generator of
$\Dd(H^i_m(R/I))$ is given by the natural map
\[
    \beta: D(H^i_m(R/I)) \to D(R^1 \tensor H^i_m(R/I)) \to[\Psi^{-1}] R^1
    \tensor D(H^i_m(R/I)).
\]
Applying local duality in the same fashion as above, $\beta$ gets
identified with the map
\begin{equation}\label{eqn.DdfromExt2}
    \beta: \Ext_R^{n-i}(R/I,R) \to \Ext_R^{n-i}(R/I^{[p]},R) \to[\Psi^{-1}] R^1
    \tensor \Ext_R^{n-i}(R/I,R)
\end{equation}
where the first part is induced from the natural surjection $R/I^{[p]} \to
R/I$, and $\Psi$ is the natural isomorphism coming from the natural
transformation $\Psi: R^1 \tensor \Hom(\usc,R) \cong \Hom(R^1 \tensor
\usc, R)$ (cf.\ Equation \eqnref{eqn.NatTransForHom} on page
\pageref{eqn.NatTransForHom}). It is straightforward that this natural
transformation for $\Hom$ induces a natural transformation on its right
derived functors $\Ext$.

If we specialize to the case that $i = d = \dim R/I$, then $H^d_m(R/I)$ is
the top local cohomology module of $A=R/I$, and $\Dd(H^d_m(R/I))=H^c_I(R)$
where $c$ denotes the codimension of $\Spec A$ in $\Spec R$. Since
$H^d_m(R/I)$ is a unit $A[F]$--module it is a $F$--full $R[F]$--module.
All this means is that the map
\[
    \theta: R^1 \tensor_R H^d_m(R/I) \to[\pi \tensor \id] A^1 \tensor_A H^d_m(R/I)
    \to[\cong] H^d_m(R/I)
\]
is surjective, which is clear since the first one is induced from the
projection $R \onto A$,  and the second is an isomorphism. Therefore,
$\beta: \Ext_R^{c}(R/I,R) \to R^1 \tensor \Ext_R^{c}(R/I,R)$ is a root
morphism for $H^i_I(R)$. By definition, $Ext^c_R(A,R)=\omega_A$ is the
canonical module for $A$. With this identification the generator $\beta$
of above can be written as
\begin{equation}\label{eqn.HcIviaomega}
    \beta: \omega_A \to R^1 \tensor \omega_A.
\end{equation}

In the next chapter we are dealing with rings that are not necessarily
complete and our technique is to complete and use the constructions
introduced here. To receive information about the original (uncompleted)
objects we have to understand which of them arise from objects via
completion. The following is a typical statement of that kind. Let $(R,m)$
be regular local and $A=R/I$ we denote by $\wh{R}$ the completion of $R$
along $m$, by $\wh{I}=I\wh{R}$ the expansion of $I$ to $\wh{R}$ and by
$\wh{A}=\wh{R}/\wh{I}=\wh{R}\tensor R/I$ the quotient, which is the same
as the completion of $A$ along its maximal ideal. When we refer to Matlis
duality $D$ we now mean the functor $\Hom_{\wh{R}}(\usc, E_{\wh{R}})$
which for cofinite $R$--modules is the same as $\Hom_R(\usc,E_R)$.
\begin{proposition}\label{prop.DdasLimitofExt}
    With the notation just indicated\/ $\Dd(H^i_m(A)) \cong \wh{R} \tensor_R
    H^{n-i}_I(R)$ and it arises as the direct limit of
    \[
        \wh{R} \tensor \Ext^{n-i}(R/I,R) \to \wh{R} \tensor R^e \tensor
        \Ext^{n-i}(R/I,R) \to \ldots
    \]
    where the map is the identity on\/ $\wh{R}$ tensored with the one
    induced form the surjection\/ $R/I^{[p^e]} \to R/I$ on\/ $\Ext^{n-i}$.
\end{proposition}
\begin{proof}
We first observe that $H^{n-i}_{\wh{I}}(\wh{R}) \cong \wh{R} \tensor
H^{n-i}_I(R)$, by flatness of completion. This can be checked to be an
isomorphism of unit $\wh{R}[F^e]$--modules where the structure on the
right is induced from the natural unit $R[F^e]$--structure on
$H^{n-i}_I(R)$ (cf.\ page \pageref{sec.PropRFmod}). The natural
$R[F^e]$--structure on $H^{n-i}_I(R)$ is obtained by writing it as the
unit $R[F^e]$--module generated by
\[
    \beta^e: \Ext^{n-i}(R/I, R) \to \Ext^{n-i}(R/I^{[p^e]},R) \to[\Psi] R^e \tensor \Ext^{n-i}(R/I,R)
\]
where the first part is induced by the projection $R/I^{p^e} \to R/I$ and
$\Psi$ is the natural transformation of Equation
\eqnref{eqn.NatTransForHom}. Completing $\beta^e$ obviously yields the
natural unit $R[F^e]$--module generator
\[
    \wh{\beta}^e: \Ext^{n-i}(\wh{R} \tensor R/I, \wh{R}) \to {\wh{R}^e \tensor \Ext^{n-i}(\wh{R} \tensor R/I,\wh{R})}
\]
of $H^{n-i}_{\wh{I}}(\wh{R})$. We used the flatness of completion to
obtain the natural isomorphism $\wh{R} \tensor \Ext^{n-i}(\wh{R} \tensor
R/I,R) \cong \Ext^{n-i}(\wh{R} \tensor R/I, \wh{R})$ from Proposition
\ref{prop.HomAndTensor} together with the natural isomorphism $\wh{R}^e
\cong \wh{R} \tensor R^e$ of $\wh{R}$--$R$--bimodules.
\end{proof}
Again, specializing to the case $i = d = \dim R/I$, we see that $H^c_I(R)$
is the unit $R[F^e]$--module generated by the map $\omega_A \to R^e
\tensor \omega_A$ induced from $\beta^e$ under the identification of
$\Ext^c(R/I,R)$ with $\omega_A$. Now $\omega_{\wh{A}} \cong \wh{R} \tensor
\omega_A$ by \cite[Theorem 3.3.5]{BrunsHerzog}. Thus $\Dd(H^d_m(R)) =
\wh{R} \tensor H^c_I(R)$ is the $R[F^e]$--module generated by the
completion of this map; \ie
\begin{equation}\label{eqn.LlviaOmega}
\begin{split}
    \Dd(H^d_m(R)) &= \dirlim( \wh{R} \tensor \omega_A \to \wh{R} \tensor
    R^e \tensor \omega_A \to \wh{R} \tensor R^{2e} \tensor \omega_A \to
    \ldots)  \\
                  &= \wh{R} \tensor (\dirlim (\omega_A \to R^1 \tensor
                  \omega_A \to R^2 \tensor \omega_A \to \ldots\ )).
\end{split}
\end{equation}
This description of $\Dd(H^d_m(A))$ shows that the canonical module
$\omega_A$ is a root of $\Dd(H^d_m(A))$.

\subsection{Graded case}\label{sec.GradedDd}
The use of Matlis duality in the previous section required the
completeness of the ring $R$. This allows us to conclude that the newly
defined dual $\Dd(M)$ of an $R[F^e]$--module $M$ which is a cofinite
$R$--module is a \emph{finitely generated} unit $R[F^e]$--module. Instead
of working in the category of modules over a complete ring we can also
work in the category of modules over a graded ring, since there, too, we
have a graded Matlis Duality theory available which is an equivalence
between finitely generated and cofinite $R$--modules. The setup is
literally the same as what we discussed so far, replacing the word
``complete'' by ``graded'' and using the appropriate graded equivalents of
the injective hull $E_R$, of $\Hom$ and of $D= \Hom(\usc,E)$. A good
reference for Matlis duality in the graded context is \cite[Section
3.6]{BrunsHerzog} where the theory is also developed in analogy with
Matlis duality in the complete case. Since the reason for treating the
graded case for us is its usefulness for concrete examples we treat it
with efficiency and not generality in mind. Thus, by a graded ring $(R,m)$
we always mean a finitely generated $\NN$--graded $k$--algebra where the
degree zero part $k$ is a perfect field of characteristic $p$. We denote
the graded maximal ideal of $R$ consisting of the elements of positive
degree by $m$.
\begin{proposition}
    Let\/ $(R,m)$ be a finitely generated, regular, graded\/ $k$--algebra. Then the graded Matlis dual
    functor\/
    ${}^*D={}^*\Hom_R(\usc,{}^*E_R)$ is naturally extended to an exact functor\/ $\Dd$
    on\/
    $R[F^e]$--modules. Furthermore if\/ $(M,\theta)$ is cofinite as an\/
    $R$--module, then
    \[
        \Dd(M) = \dirlim( {}^*D(M) \to R^e \tensor {}^*D(M) \to R^{2e} \tensor
        {}^*D(M) \to \ldots\ )
    \]
    is a finitely generated unit\/ $R[F^e]$--module.
\end{proposition}
\begin{proof}
 The
${}^*$ in the above statement is a reminder that we are in fact working
with the graded homomorphisms and consider the natural grading on the
injective hull of $R/m$ (cf.\ notation in \cite[Section
3.6]{BrunsHerzog}). Our formalism for extending functors on $R$--modules
to functors on $R[F^e]$--modules also applies to this graded Matlis dual
functor. Furthermore, one can show the graded analog of Proposition
\ref{prop.HomAndTensor}, \ie there is a natural isomorphism of functors
\[
    \F[e] \circ {}^*\Hom_R(\usc,{}^*E_R) \to[\cong]
    \Hom_R(\F[e](\usc),\F[e]({}^*E_R))
\]
Here we understand the grading on $\F[e](M)=R^e \tensor M$ such that $r
\tensor m$ has degree $\deg(r) + p^e\deg(m)$. With this grading, the
natural $R[F^e]$--module structure on $E_R$ is graded of degree zero.
Proposition \ref{prop.ContraExtendsWell} now shows that we can extend
${}^*D$ to a functor ${}^*\Dd$ from $R[F^e]$--modules to unit
$R[F^e]$--modules. ${}^*\Dd(M)$ is the unit $R[F^e]$--module generated by
\[
    {}^*D(M) \to \F[e]({}^*D(M)).
\]
If $M$ is cofinite, then ${}^*D(M)$ is a finitely generated $R$--module.
Therefore, ${}^*\Dd(M)$ is a finitely generated unit $R[F^e]$--module.
\end{proof}
The remaining properties of ${}^*\Dd$ are derived equivalently to the
complete case. This is fairly straightforward, but to deal with this in
complete detail one has to develop a sensible theory of graded
$R[F]$--modules, including notions of graded roots. This, also, can be
done quite naturally. As already indicated in the last proof, for a graded
$R$--module $M$, one gives $\F[e]M$ the grading such that $\deg(r' \tensor
m) = \deg(r') + p^e \deg(m)$.
\begin{definition}
Let $R$ be a graded $k$--algebra. A graded $R[F^e]$--module is a graded
$R$--module $M$, together with a graded map
\[
\theta^e: \F[e]M \to M
\]
of degree zero. We call $(M, \theta^e)$ a  graded unit $R[F^e]$--module if
$\theta^e$ is an isomorphism.
\end{definition}\index{graded $R[F^e]$--module}
Of course, the notion of a graded root, graded generator and so forth is
defined. The theory develops with the same results as in the non--graded
(complete) case. In particular, we have the existence of graded roots for
graded unit $R[F^e]$--modules which are $R[F^e]$--generated by finitely
many homogeneous elements. Even better, since Chevalleys Theorem, in the
form of Lemma \ref{lem.Chev}, holds also in the graded context, we get
that finitely generated graded unit $R[F^e]$--modules have a \emph{unique}
minimal graded root (cf.\ Proposition \ref{prop.MinimalRoot})\footnote{The
proofs of the graded version of Chevalleys Lemma and Proposition
\ref{prop.MinimalRoot} are the same as in the complete case, found in
\cite[Lemma 3.3]{Lyub} or on page \pageref{prop.MinimalRoot} of this
text}. With a developed theory of graded $R[F^e]$--modules as indicated it
is not hard to show the properties of ${}^*\Dd$, analogous to the ones of
$\Dd$ developed before. To ease notation we now drop the star ${}^*$ and
denote ${}^*\Dd$ just by $\Dd$. We summarize the results and only refer to
the proofs of the complete analogues, which mostly transfer word for word
to the graded case.

\begin{proposition}\label{prop.DdGraded}
    Let\/ $(R,m)$ be a regular, graded\/ $k$--algebra and let\/ $M$ be a graded\/
    $R[F^e]$--module that is cofinite as an\/ $R$--module. Then
    \begin{enumerate}
    \item $\Dd(M)=0$ if and only if\/ $M$ is\/ $F$--nilpotent.\/ $\Dd(M) \cong
        \Dd(N)$ if and only if\/ $\Ffred{M} = \Ffred{N}$.
    \item If\/ $M$ is\/ $F$--full, then\/ $D(M)$ is a root of\/ $\Dd(M)$. If\/ $M$ is
        also\/ $F$--reduced, then\/ $D(M)$ is the unique minimal root.
    \item Every graded unit\/ $R[F^e]$--submodule of\/ $\Dd(M)$ arises as\/ $\Dd(N)$
        for some graded\/ $R[F^e]$--submodule of\/ $M$.
    \item $\Dd$ is an isomorphism between the lattice of graded\/
        $R[F^e]$--modules quotients of\/ $M$ (up to $\Ffred{(\usc)}$) and the lattice of
        unit\/
        $R[F^e]$--submodules of\/ $\Dd(M)$.
    \end{enumerate}
\end{proposition}
\begin{proof}
At the beginning, the proof of (a) is the same as the proof of Proposition
\ref{prop.DisZeroOnNilpotent} replacing ``complete'' by ``graded'' and
working with the graded notion of $R[F]$--module, $D$, $\Dd$ and so forth.
For the second part of (a) one has to be more careful. The key ingredient
is the unique minimal root in the graded category, whose existence we
discussed above. Given this, the statement follows analogous to the proof
of Lemma \ref{lem.DandRoots} and Proposition \ref{prop.Dd=ifFfred=}.

For the proof of (b) the same note of caution is at order. As soon as we
have the existence of the unique minimal root in the graded category the
proof is the same as in the complete case (cf. Proposition
\ref{prop.Full=Root} and Corollary \ref{cor.Dduniqueroot}).

Part (c) and (d) are proven exactly like Theorem \ref{thm.DdofSimpisSimp}.
\end{proof}

%%%%%%%%%%%%%%%%%%%%%%%%%%%%%%%%%%%%%%%%%%%%%%%%%%%%%%%%%%%%%%
%%                                                          %%
%%   This is file: chapter5.tex                             %%
%%   It contains the Fourth Chapter of my                   %%
%%   dissertation: dissertation.tex                         %%
%%                                                          %%
%%%%%%%%%%%%%%%%%%%%%%%%%%%%%%%%%%%%%%%%%%%%%%%%%%%%%%%%%%%%%%

\chapter{The intersection homology $D_R$--module in finite
characteristic}\label{chap.IntHomD_R--mod}

In this chapter the second central result of the dissertation is proved,
and applications of the construction are given. This main result is the
existence of an analog of the Kashiwara--Brylinski intersection cohomology
$D_R$--module in finite characteristic.

\newcommand{\mainthmtwocontent}{
\begin{maintheorem2}
    Let\/ $(R,m)$ be regular, local and\/ $F$--finite. Let\/ $I$ be an ideal
    with\/ $\height I = c$ such that\/ $A=R/I$ is analytically irreducible.
    Then the local cohomology module\/ $H^c_I(R)$ has a unique
    simple\/ $D_R$--submodule\/ $\Ll(A,R)$.
\end{maintheorem2}}
\mainthmtwocontent

To sketch a proof, we first assume that $R$ is complete. Recall that
\idx{analytically irreducible} just means that the completion of $A$ is a
domain. One key ingredient in the proof of Main Theorem 2 is a result from
the theory of tight closure which says that the tight closure of zero,
$0^*_{H^d_m(A)}$, is the unique maximal proper $R[F^e]$--submodule of
$H^d_m(R/I)$ for all $e > 0$. With this at hand, one uses the extended
Matlis dual functor $\Dd$ of the last chapter to obtain a unique simple
$R[F^\infty]$--submodule $L$ of $H^c_I(R)$. That is,
$L=\Dd(H^d_m(R/I)/0^*_{H^d_m(R/I)})$.

Now, the property of having a unique simple unit
$R[F^{\infty}]$--submodule descends from faithfully flat extensions. Thus
we showed that $H^c_I(R)$ has a unique simple unit
$R[F^\infty]$--submodule $L$, regardless if $R$ was complete. For the same
reason, namely that the property of having a unique simple
$D_R$--submodule descends from faithfully flat extensions, we are allowed
to extend the field $\FF_p$ to an algebraically closed and sufficiently
huge field $K$, \ie we replace $R$ by $K \tensor_{\FF_p} R$. Then we can
apply Main Theorem 1 to conclude that $L$ is $D_R$--simple and Main
Theorem 2 follows.

Below, we first review some of the notation from the theory of tight
closure. On the way we recall a proof of the aforementioned result, that
$H^d_m(A)$ has a unique maximal proper $R[F^\infty]$--submodule, provided
the completion of $A$ is a domain. Also, we recall some results about test
ideals and test modules which play a role in understanding the concrete
structure of $\Ll(A,R)$.

Then we give the proof or Main Theorem 2 as outlined above. This proof is
mainly existential, \ie we do not pay attention to what can be said about
the construction of the unique simple $D_R$--submodule of $H^c_I(R)$. The
section following the proof is spent on sharpening the argument in order
to get a concrete description of this simple submodule. The problems one
encounters stem from the two reduction steps of completing, and making a
field extension. These problems are overcome by showing the respective
base change properties of $\Ll(R,A)$. A consequence of this is that the
unique simple $D_R$--submodule is at the same time the unique simple
$R[F^e]$--module. Furthermore $\Ll(R,A)$ arises as a direct limit of the
Frobenius powers of the parameter test module. This, in turn, shows that
the parameter test module commutes with completion.

The concrete description of $\Ll(A,R)$ gives also an explicit criterion
for when $H^c_I(R)$ is $D_R$--simple. Namely, the tight closure of zero,
$0^*_{H^d_m(A)}$, is $F$--nilpotent, if and only if $H^d_m(A)$ is
$D_R$--simple. This implies that if $A$ is $F$--rational, then $H^c_I(R)$
is $D_R$--simple.

We finish this chapter with some examples of graded complete intersections
where we explicitly calculate $\Ll(A,R)$. An application of our simplicity
criterion to rings $A$ for which the map from their $F$--rational
normalization is injective is also given. This shows the finite
characteristic analog of a theorem of S.P.Smith \cite{SmithSP.IntHom} on
the $D_R$--simplicity of $H^1_I(R)$ where $A=R/I$ is a plane cusp.

\section{Background in tight closure theory}
Tight closure\index{tight closure} is a powerful tool in commutative
algebra introduced by Mel \idx{Hochster} and Craig \idx{Huneke} about
fifteen years ago \cite{HH88}. It is a beautiful theory with many
applications to algebraic geometry. Among its greatest accomplishments are
a very simple proof of the Hochster--Roberts theorem on the
Cohen--Macaulayness of rings of invariants of reductive groups (see
\cite{HoRo} for original and \cite{HH89} for the tight closure proof), new
proofs and generalizations of the Brian\c{c}on--Skoda theorem and the
Syzygy theorem of Evans and Griffith \cite{HH90} and uniform Artin--Rees
theorems \cite{Hune.UniBounds}.

On a more geometric side, tight closure is intimately related to Kodaira
vanishing \cite{HuSmith.KodVan} and the Fujita conjectures
\cite{Smith.fuji}. There is a strong connection between the singularities
arising in the minimal model program, and singularities obtained from
tight closure theory \cite{Smith.sing}. One of the most significant is the
equivalence of the notions of rational singularity and $F$--rational type
which was established by Smith \cite{Smith.rat} and Hara \cite{Hara}. The
notion of $F$--rationality arises naturally from tight closure: the local
ring $(A,m)$ is called $F$--rational if all ideals $I$ generated by a full
system of parameters are tightly closed, \ie $I = I^*$. Some of the
techniques which are used to obtain this equivalence also play an
important role in establishing our second Main Theorem; but first, we have
to introduce the background in tight closure theory we require. For a more
detailed introduction to this beautiful subject we recommend
\cite{Smith.IntroTight,Hune.tight} and later the more technical original
papers \cite{HH90,HH89}.

Let $A$ be a noetherian ring. We denote by $A^\circ$ the subset of
elements of $r$ that are not contained in any minimal prime of $A$. Let $N
\subseteq M$ be a submodule of $M$. The tight closure $N^*_M$ of $N$
inside of $M$ is defined as follows:
\begin{definition}
Let $A$ be noetherian and $N \subseteq M$. The tight closure $N^*_M$ (or
just $N^*$ if $M$ is clear from the context) consists of all elements $m
\in M$, such that there exists a $c \in A^\circ$, such that for all $e \gg
0$
\[
    c \tensor m \in N^{[p^e]}.
\]
Here $N^{[p^e]}$ denotes the image of $\F[e]N$ in $\F[e]M$ and $c \tensor
m$ is an element of $\F[e]M$.
\end{definition}\index{$N^*$}
If $N=I$ is just an ideal of $A$, the definition is much more transparent.
In this case $r \in A$ is in $I^*$ if and only if there is $c \in A^\circ$
such that $cr^{p^e} \in I^{[p^e]}$ for all $e \gg 0$. A module is tightly
closed if $N^* = N$. We have that $N \subseteq N^*$ as one expects from a
decent closure operation. If $N$ is noetherian, then $N^* = (N^*)^*$.
There are two related closure operations which are important for us.
\begin{definition}
    Let $N \subseteq M$ be $A$--modules. The \emph{finitistic tight closure}\index{finitistic tight closure}
    of
    $N$ inside of $M$ consists of all elements $m \in (N \cap M_0)^*_{M_0}$
    for some finitely generated $M_0 \subseteq M$. It is denoted by
    $N^{*fg}_M$.

    The \emph{Frobenius closure}\index{Frobenius closure} $N^F_M$ consists of all elements $m \in
    M$ such that $1 \tensor m \in N^{[p^e]}$ for some $e \geq 0$.
\end{definition}
We immediately see that $N^{*fg} \subseteq N^*$ and that equality holds if
$M$ is finitely generated. Clearly, $N^F \subseteq N^*$. For the zero
submodule of the top local cohomology module of an excellent, local,
equidimensional ring $A$, the finitistic tight closure is equal to the
tight closure, \ie $0^{*fg}_{H^d_m(A)}=0^*_{H^d_m(A)}$ (see
\cite[Proposition 3.1.1]{SmithDiss}). In general, it is a hard question to
decide if the tight closure equals the finitistic tight closure, and it is
related to aspects of the localization problem in tight closure theory
(cf.\ \cite{LyubSmith.Comm}).

As our focus lies on modules with Frobenius actions we study the above
closure operations in this case more closely. The following is an
important proposition.
\begin{proposition}[cf.\ \protect{\cite[Proposition 4.2]{LyubSmith.Comm}}]
    Let\/ $A$ be noetherian and let\/ $M$ be an\/ $A[F^e]$--module. If\/ $N$ is
    a\/
    $A[F^e]$--submodule, then so are\/ $N^*_M$, $N^{*fg}_M$ and\/ $N^F_M$.
\end{proposition}
\begin{proof}
Let $x \in N^*$, \ie we have a $c \in A^\circ$ such that $c \tensor x \in
N^{[p^{e'}]}$ for all $e' \gg 0$. Denote by $F^e$ and $\theta^e$ the given
Frobenius action and operation on $M$. We have to show that $y
\defeq F^e(x)$ is also in $N^*$. For this first note that, by definition,
$y = \theta^e(1 \tensor x)$. Now consider the following commutative
diagram
\[
\xymatrix@C=4pc{
            {A^{e'+e}\tensor N} \ar^{\id_{A^{e'}}\tensor \theta^e}[r] \ar^{}[d] &{A^{e'} \tensor N} \ar[d] \\
            {A^{e'+e}\tensor M} \ar^{\id_{A^{e'}}\tensor \theta^e}[r] &{A^{e'}\tensor M}
}
\]
This diagram makes sense by the assumption that $N$ is an
$A[F^e]$--submodule, \ie $\theta^e$ naturally induces a Frobenius
structure (also denoted by $\theta^e$) on $N$, such that the diagram
commutes. The vertical maps are of course induced from the inclusion $N
\subseteq M$. Observe that $c \tensor y \in A^{e'} \tensor M$ is the image
of $c \tensor x \in A^{e'+e} \tensor M$, by definition of $y$. But $c
\tensor x$ is, for $e' \gg 0$, in the image of $A^{e'+e} \tensor
N=N^{[p^{e'+e}]}$. By the commutativity of the above diagram, this implies
that $c \tensor y$ is in the image of $A^{e'} \tensor N$ in $A^{e'}\tensor
M$, \ie $y \in N^*$. This shows that $N^*$ is in fact a
$A[F^e]$--submodule of $M$.

The case of $N^{*fg}_M$ can be proved similarly. For $N^F_M$ we just note
that $N^F_M$ consists exactly of the elements of $M$ which are mapped to
the $F$--nilpotent part of $M/N$.
\end{proof}
An alternative proof is found by observing that $(N^*)^{[p^e]} \subseteq
(N^{[p^e]})^*$. Then apply $\theta^e$ and use the easily verifiable fact
that $\theta^e(\usc^*) \subseteq \theta^e(\usc)^*$ to see that
\[
    F^e(N^*)=\theta^e((N^*)^{[p^e]}) \subseteq \theta^e((N^{[p^e]})^*)
    \subseteq (F^e(N))^* \subseteq N^*
\]
which finishes the argument. From this we get as an immediate corollary
that the tight closure of the zero $A[F^e]$--submodule is a Frobenius
stable submodule of any $A[F^e]$--module.
\begin{corollary}\label{cor.0*isRFsub}
    Let\/ $A$ be a ring and let\/ $(M,F^e)$ be an\/ $A[F^e]$--module.
    Then\/
    $0^{*fg}_M$, $0^*_M$ and\/ $0^F_M=\Fnil{M}$ are\/ $A[F^e]$--submodules of\/ $M$.
\end{corollary}\index{$0^*$}
Before proceeding further we have to study an important tool of tight
closure theory more closely.

\subsection{Test ideals and test modules}
The elements $c$ occurring in the definition of tight closure play a
special role. Those amongst them, that work for all tight closure tests
for all submodules of all finitely generated $A$--modules are called the
test elements of $A$.
\begin{definition}
    An element $c \in A^\circ$ is called a \emph{\idx{test element}} if for all
    submodules $N \subseteq M$, of every finitely generated $A$--module
    $M$,
    we have $cN^*_M \subseteq N$. A test
    element is called \emph{completely stable test element} if its image
    in the completion of every local ring of $A$ is a test element.
\end{definition}
It is shown in \cite[Proposition 8.33]{HH90}, that it is enough to range
over all ideals of $A$ in this definition, \ie $c$ is a test element if
and only if for all ideals $I$ and all $x \in I^*$ we have $cx^{p^e} \in
I^{[p^e]}$ for all $e \geq 0$. Thus, the test elements are those elements
$c$ occurring in the definition of tight closure which work for all tight
closure memberships of all submodules of all finitely generated
$A$--modules. A nontrivial key result is that in most cases, test elements
(and even completely stable test elements) exist:
\begin{proposition}\label{prop.TestEl}
    Let\/ $A$ be reduced and of finite type over an excellent local ring.
    Then\/ $A$ has completely stable test elements. Specifically, any element\/ $c
    \in A^\circ$ such that\/ $A_c$ regular has a power which is a completely
    stable test element.
\end{proposition}
The proof of this is quite technical and can be found in \cite[Chapter 6
]{HH89}. Results in lesser generality (for example, when $A$ is
$F$--finite) are obtained fairly easily: for a good account see
\cite{Smith.IntroTight,Hune.tight}.

The ideal $\tau_A$ generated by all test elements is called the test
ideal. As remarked, $\tau_A = \bigcap (I:_A I^*)$ where the intersection
ranges over all ideals $I$ of $A$. This naturally leads one to consider
variants of the test ideal by restricting the class of ideals this
intersection ranges over. The \emph{\idx{parameter test ideal}} of a local
ring $(A,m)$ is the ideal $\wt{\tau}_A = \bigcap (I:_A I^*)$ where the
intersection ranges over all ideals generated by a full system of
parameters. If $A$ is Cohen--Macaulay, it follows from the definition of
$H^d_m(A)$ as $\dirlim A/(\xd)^{[p^e]}$ that $\wt{\tau}_A =
\Ann_A(0^*_{H^d_m(A)})$ \cite[Proposition 4.1.4]{SmithDiss} where $\xd$ is
a system of parameters for the local ring $(A,m)$. If $A$ is only an
excellent domain, then $\wt{\tau}_A \subseteq \Ann_A(0^*_{H^d_m(A)})$.
Further generalizing, the \emph{\idx{parameter test module}} is defined as
$\tau_{\omega_A} = \Ann_{\omega_A} 0^*_{H^d_m(A)} = \omega_A \cap
\Ann_{\omega_{\wh{A}}}0^*_{H^d_m(\wh{A})}$ where the action of $\omega_A$
on $H^d_m(A)$ is the one coming from the Matlis duality pairing $H^d_m(A)
\times \omega_{\wh{A}} \to E_A$. Of course we require here that $A$ has a
canonical module.
\begin{lemma}\label{lem.tauNonzero}
    Let\/ $A$ be reduced, excellent, local and equidimensional with canonical
    module\/ $\omega_A$. If\/ $c$ is a parameter test element, then\/
    $c\omega_A \subseteq \tau_{\omega_A}$. In particular,\/ $\tau_{\omega_A}$ is
    nonzero.
\end{lemma}
\begin{proof}
Let $c$ be a parameter test element. In particular, $c$ annihilates the
finitistic tight closure of zero in $H^d_m(A)$. Therefore, for every $\phi
\in \omega_A$ and $\eta \in 0^*_{H^d_m(A)}=0^{*fg}_{H^d_m(A)}$ we have
$c\phi \cdot \eta = \phi \cdot (c\eta) = \phi \cdot 0 = 0$ where
``$\cdot$'' represents the Matlis duality pairing. This shows that
$c\omega_A \subseteq \tau_{\omega_A}$. The hypotheses on $A$ ensure by
\cite[Remark 2.2(e)]{HH.CanMod} that the canonical module is faithful, \ie
$c\omega_A \neq 0$. Therefore the last part of the lemma follows from the
existence of test elements (Proposition \ref{prop.TestEl}), since a test
element is also a parameter test element.
\end{proof}

One of the persistently difficult and only partially solved problems in
tight closure theory is the behavior of various aspects of this theory
under localization and completion. For the test ideals and test modules we
ask the following. Is the localization/completion of the (parameter) test
ideal/module of $A$ the (parameter) test ideal/module of the
localization/completion of $A$? For the test ideal some known cases where
it, in fact, commutes with localization and completion are $A$ an isolated
singularity, $A$ $\QQ$--Gorenstein or $A$ Gorenstein on its punctured
spectrum \cite{LyubSmith.Comm}. Since the parameter test ideal, the
parameter test module, and the test ideal are equal, in the case that $A$
is Gorenstein, all three commute with localization and completion in this
case. If we assume that $A$ is Cohen--Macaulay, and complete and local
then Smith shows in \cite{Smith.test} that the parameter test module and
the parameter test ideal commute with localization. This is unknown for
the test ideal.

For our purposes the parameter test module is most important, since, under
Matlis duality it exactly corresponds to the tight closure of zero in
$H^d_m(A)$. But not only do we use the parameter test module, we also
obtain results about it. As a consequence of our theory, we will show
(Corollary \ref{cor.TestCommCompl}) that the parameter test module
commutes with completion. The next lemma will be needed later: it is a
step towards a base change property for the parameter test module.
\begin{lemma}\label{lem.tauFieldExt}
    Let\/ $A=R/I$ be a local algebra over a perfect field\/ $k$ where\/ $R$ is regular.
    Let\/ $A_K = K \tensor_k A$ where\/ $K$ is an extension field of\/ $k$.
    Then\/ $K \tensor_k \tau_{\omega_A} \subseteq \tau_{\omega_{A_K}}$ and
    thus\/ $\tau_{\omega_A} \subseteq A \cap \tau_{\omega_{A_K}}$ is
    nonzero.
\end{lemma}
\begin{proof}
Since $A \to A_K$ is faithfully flat we see that $K \tensor_k H^d_m(R)
\cong H^d_m(A_K)$. For the same reason
\[
    \omega_{A_K} = \Ext^c_{R_K}(A_K, R_K) \cong K \tensor_k
    \Ext^c_R(A,R) = K \tensor_k \omega_A.
\]
Here we used that, since $k$ is perfect, $K$ is separable over $k$ and
thus $K \tensor_k R$ is also regular. Using \cite[Theorem
4.1]{AbEn.TestIdAndBaseChange}, it follows that $K \tensor_k
0^*_{H^d_m(A)} = 0^*_{H^d_m(A_K)}$ and therefore
\begin{equation*}
    \tau_{\omega_{A_K}} = \Ann_{\omega_{A_K}} 0^*_{H^d_m(A_K)}
                         = \Ann_{A \tensor_k \omega_{A}} K \tensor_k
                         0^*_{H^d_m(A)}.
\end{equation*}
Now one immediately sees that this contains $K \tensor_k \Ann_{\omega_A}
0^*_{H^d_m(A)} = K \tensor_k \tau_{\omega_A}$. The second statement
follows observing that by Lemma \ref{lem.tauNonzero}, $\tau_{\omega_A}$ is
nonzero and thus so is $K \cap \omega_{A_K}$.
\end{proof}

\subsection{$F$-rationality and local cohomology}
The tight closure of zero in the top local cohomology module $H^d_m(A)$ of
a local ring $(A,m)$ plays a role as the obstruction to
\index{F-rationality@$F$--rationality} $F$--rationality of $A$. Its
distinguishing property is that it is the maximal proper $A[F]$--submodule
of $H^d_m(A)$. Precisely the following is the case:
\begin{theorem}\label{thm.MaxRFofH}
    Let\/ $(A,m)$ be reduced, excellent and analytically irreducible.
    Then, the tight closure of zero,\/ $0^*_{H^d_m(A)}$, in\/ $H^d_m(A)$ is the unique maximal
    proper\/
    $A[F^e]$--submodule of\/ $H^d_m(A)$. The quotient\/ $H^d_m(A)/0^*_{H^d_m(A)}$ is a nonzero
    simple\/
    $F$--reduced and\/ $F$--full\/ $A[F^\infty]$--module.
\end{theorem}
Before we prove this result we have to recall the following fact from
tight closure theory. It says that for finitely generated $M$, in order to
check the tight closure membership, it is enough to require $c \tensor m
\in N_M^{[p^{e'}]}$ for infinitely many $e'$ (instead of all but finitely
many as in the definition of tight closure). This fact is needed to show
that $0^*_{H^d_m(A)}$ is the maximal proper $A[F^e]$--module for all $e$,
and not just the maximal proper $A[F]$--module. The case $e = 1$ of
Theorem \ref{thm.MaxRFofH} was shown by Smith as \cite[Theorem
3.1.4]{SmithDiss} and with the just mentioned fact, our proof is just a
straightforward adaption of the proof there.
\begin{proposition}\label{prop.TightEnough}
    Let\/ $A$ be reduced and excellent.
    For a submodule\/ $N$ of an\/ $A$--module\/ $M$ we have that\/ $m \in N^{*fg}$ if
    and only if there exists a\/ $c \in A^\circ$ such that\/ $c \tensor m \in (N \cap M_0)^{[p^e]}$
    for infinitely many\/ $e \geq 0$ and some finitely generated\/ $M_0 \subseteq M$.

    If\/ $N^* = N^{*fg}$, then\/ $m \in N^*$ if and only if there is\/ $c
    \in A^\circ$ such that\/ $c \tensor m \in N^{[p^e]}$ for infinitely
    many\/
    $e \geq 0$.
\end{proposition}
\begin{proof} Exercise 2.3 of \cite{Hune.tight} shows this
for ideals and the same proof (solution in the back of \cite{Hune.tight})
works for finitely generated $A$--modules. As membership in the finitistic
tight closure is really a membership question of tight closure in a
finitely generated setting this applies to our statement. The last part is
also proven analogously.
\end{proof}

\begin{proof}[Proof of Theorem \ref{thm.MaxRFofH}] First, Corollary
\ref{cor.0*isRFsub} shows that $0^*_{H^d_m(A)}$ is an
$A[F^e]$--sub\-mo\-du\-le for all $e \geq 0$. The assumption that $A$ is
an excellent domain ensures the existence of completely stable test
elements, by Proposition \ref{prop.TestEl}. With such an element $c \in A$
it follows immediately that the tight closure of zero in $H^d_m(A)$ is the
same, whether we think of $H^d_m(A)$ as a module over $A$ or over
$\widehat{A}$. As it is also clear that the $\widehat{A}[F^e]$--submodules
and the $A[F^e]$--submodules are the same we can assume that $A$ is a
complete domain. The same test element $c$ also shows that the annihilator
of $0^*_{H^d_m(A)}=0^{*fg}_{H^d_m(A)}$ in $A$ has positive height and
therefore, $0^*_{H^d_m(A)}$ cannot be all of $H^d_m(A)$.

We have to show that, if $m \in H^d_m(A)$ generates a proper
$A[F^e]$--submodule of $H^d_m(A)$, then $m \in 0^*_{H^d_m(A)}$. So, let
$M=A[F^e]m$ be the proper $A[F^e]$--submodule generated by $m$ and denote
by $C$ the cokernel of the inclusion $M \subseteq H^d_m(A)$. Applying
Matlis duality to the resulting short exact sequence we get a short exact
sequence of finitely generated $A$--modules.
\[
    0 \to D(C) \to \omega_{A} \to D(M) \to 0.
\]
Since $A$ is a domain, the canonical module $\omega_{A}$ is torsion free
of rank one \cite{HH.CanMod}. Tensoring this exact sequence with the field
of fractions $K$ of $A$, we conclude that either $K \tensor D(C)=0$ or $K
\tensor D(M)=0$. In the first case it follows that $D(C)=0$ as it is also
torsion free (as a submodule of the torsion free $\omega_{A}$). Thus $C =
0$ and $M = H^d_m(A)$, contradicting the assumption that $M$ is a proper
$A[F^e]$--submodule.

Thus $K \tensor D(M) = 0$ which, by the noetherian property of $D(M)$,
implies that some $c \in \widehat{A}$ annihilates $D(M)$ and therefore
also $M=D(D(M))$. But this means that $c$ simultaneously kills all
$F^{re}(m)=\theta^{er}(1 \tensor m)$. As $\theta^e$ is an isomorphism,
this is  equivalent to $c \tensor m = 0$ in $\F[er]H^d_m(A)$ for all $r
\geq 0$. Since by \cite[Proposition 3.1.1]{SmithDiss}, $0^{*fg}_{H^d_m(A)}
= 0^{*}_{H^d_m(A)}$ it follows from Proposition \ref{prop.TightEnough}
that $m \in 0^*_{H^d_m(A)}$.

Because $0^*_{H^d_m(A)}$ is the maximal proper $A[F^e]$--submodule for all
$e$, we know that $H^d_m(A)/0^*_{H^d_m(A)}$ is the unique nonzero simple
$A[F^\infty]$--module quotient. It remains to show that it is $F$--reduced
(a simple $A[F^e]$--module is $F$--full if and only if it is
$F$--reduced). For this note that the kernel of $F$ is a
$A[F^e]$--submodule and, by simplicity, it must either be zero
($F$--reduced) or all of $H^d_m(A)/0^*_{H^d_m(A)}$. In the second case,
this implies that $F(H^d_m(A)) \subseteq 0^*_{H^d_m(A)}$. Since $H^d_m(A)$
is a unit $A[F]$--module (enough that the structural map $\theta$ is
surjective) we have that $F(H^d_m(A)) = H^d_m(A)$. This contradicts the
fact that $0^*_{H^d_m(A)}$ is a proper submodule which we showed above.
Thus the quotient is $F$--reduced and $F$--full.
\end{proof}
The vanishing of this unique maximal $A[F^e]$--submodule of $H^d_m(A)$
characterizes $F$--rationality of $A$. By definition, $A$ is called
$F$--rational\index{F-rational@$F$--rational} if and only if every ideal
that is generated by a system of parameters is tightly closed. Then, an
argument \cite[Proposition 3.1.1]{SmithDiss} involving the colon capturing
property of tight closure and the description of $H^d_m(A)$ as the limit
$\dirlim A/(\xd)^{[p^e]}$ for some parameter ideal $\xd$, shows that
$z+(\xd) \in H^d_m(A)$ is in the tight closure of zero if and only if $z
\in (\xd)^*$. For a Cohen--Macaulay ring $R$ it follows that $A$ is
$F$--rational if and only if $0^*_{H^d_m(A)}=0$.

Smith showed in \cite{Smith.rat} that if $A$ is $F$--rational, then it is
pseudo--rational. Pseudo--rationality is a resolution free
characterization of \idx{rational singularities}, thus it applies also to
characteristic $p$ where it is unknown whether there are resolutions of
singularities in general. The converse of Smith's result was shown by Hara
\cite{Hara} and, independently, by Metha and Srinivas \cite{MehtaSr}; they
show that if $A$ is a generic characteristic $p$ model obtained by
reduction mod $p$ from a ring $\bar{A}$ in characteristic zero and if
$\bar{A}$ has rational singularities, then $A$ is $F$--rational.

\section{The unique $D_A$--simple submodule of $H^c_I(R)$}
To cope with the mentioned reductions to the completion and field
extension we first note that the property of having a unique simple
submodule descends down from faithfully flat extensions. This is best
observed by the following very general lemma.
\begin{lemma}\label{lem.descend}
    Let\/ $M$ have finite length in some abelian category. Then\/ $M$ has a
    \emph{unique} simple sub-object if and only if\/ $M$ does not contain
    internal direct sums.

    Let\/ $C$ be a faithfully exact functor between abelian categories. Let\/ $M$ be
    such that\/ $M$ and\/ $C(M)$ have finite length. If\/ $C(M)$ has a
    unique simple sub-object, then so does\/ $M$.
\end{lemma}
\begin{proof}
The finite length condition ensures that every submodule contains a
nonzero simple submodule. Therefore, having no internal sums is equivalent
to having no internal sums of simple submodules. Now the assertion is
clear.

For the second part, assume that $M$ does not have a unique simple
submodule. By the first part, this means that $M$ contains a direct sum
$N' \oplus N'' \subseteq M$. Applying the faithfully exact functor $C$ we
get $C(N') \oplus C(N'') \subseteq C(M)$. The directness follows form the
exactness of $C$, being faithful ensures that $C(N')$ and $C(N'')$ are
nonzero. Thus $C(M)$ has an internal direct sum.
\end{proof}
The abelian categories we have in mind are, of course, $D_R$--mod and
$R[F^e]$--mod. The functor is either extending a perfect field contained
in $R$ or completion along the maximal ideal.

Now, all techniques to prove the second main theorem are present and we
can start the proof without further delay.

\mainthmtwocontent

\begin{proof}
First, we show that $H^c_I(R)$ has a unique simple
$R[F^\infty]$--submodule. The existence of such a submodule descends down
from a faithfully flat extension. Thus we prove the result first for the
completion $\widehat{R}$ of $R$ along $m$. Thus we assume for now that $R$
and $A$ are complete. In this case $L(A,R) \defeq
\Dd(H^c_m(R/I)/0^*_{H^d_m(A)})$ is the unique simple $R[F^e]$--submodule
for all $e \geq 1$. First, since by Theorem \ref{thm.MaxRFofH},
$0^*_{H^d_m(A)}$ is the maximal nonzero $R[F^\infty]$--submodule of
$H^d_m(R)$, the quotient $H^d_m(R/I)/0^*_{H^d_m(A)}$ is the unique simple
$R[F^\infty]$--quotient of $H^d_m(R/I)$. Applying $\Dd$ to this surjection
we get that $L(A,R) = \Dd(H^d_m(R/I)/0^*_{H^d_m(A)})$ is the unique simple
unit $R[F^\infty]$--submodule of $\Dd(H^d_m(R/I)) = H^c_I(R)$ where the
$F$--fullness part in Theorem \ref{thm.MaxRFofH} ensures that $L(A,R)$ is
nonzero. To see that $L(A,R)$ is the \emph{unique} simple unit
$R[F^\infty]$--submodule let $\Nn$ be a nonzero unit
$R[F^\infty]$--submodule of $H^c_I(R)$. By Theorem
\ref{thm.DdofSimpisSimp} $\Nn$ arises as $\Dd(N)$ for some nonzero
$R[F^\infty]$--module quotient $N$ of $H^c_m(R/I)$. Therefore, the kernel
of this quotient map is a proper $R[F^\infty]$--submodule of $H^d_m(R/I)$.
Therefore, it is contained in $0^*_{H^d_m(R/I)}$, the maximal proper
$R[F^\infty]$--submodule of $H^d_m(R/I)$. This implies that $L(A,R)
\subseteq \Nn$. This shows that for $R$ complete $L(A,R) =
\Dd(H^d_m(R/I)/0^*_{H^d_m(A)})$ is the unique simple
$R[F^\infty]$--submodule of $H^c_I(R)$.

Now let $R$ and $A$ be as before. Theorem \ref{thm.fguRFhaveDCC} ensures
that $H^c_I(M)$ and $H^c_{\wh{I}}(\wh{R})=\wh{R} \tensor H^c_I(R)$ have
finite length as $R[F^e]$--modules ($\wh{R}[F^e]$--modules resp.). Since
$\wh{R}\tensor H^c_I(R)$ has a unique simple $\wh{R}[F^e]$--submodule, and
by Lemma \ref{lem.descend}, $H^c_I(R)$ has a unique simple
$R[F^e]$--module.

The next step will be to show that this implies the existence of a unique
simple $D_R$--submodule. In the case that the coefficient field $k$ is
algebraically closed and has strictly bigger cardinality than the
$k$--dimension of $R$ we can apply Main Theorem 1 (resp. Remark
\ref{rem.MainOneGeneral}), and conclude right away that the unique simple
$R[F^\infty]$--submodule is, in fact, $D_R$--simple. To show it is the
unique simple $D_R$-submodule, let $N$ be a simple $D_R$--submodule of
$H^c_I(R/I)$. Then, Corollary \ref{cor.MainOne} implies that for some $e
\geq 0$, the $R[F^e]$--submodule $R[F^e]N$ decomposes into a direct sum of
simple unit $R[F^e]$--submodules, each $D_R$--isomorphic to $N$. Since we
have a unique simple unit $R[F^e]$--submodule, $R[F^e]N$ consists of just
one summand, which means that $N = R[F^e]N = L(A,R)$.

In general, we extend the coefficient field to an algebraically closed
field $K$ of sufficient cardinality. This is a faithfully flat extension
and by the same argument as we used in the complete case above (using
Lemma \ref{lem.descend}) we conclude that $H^c_I(R)$ has a unique simple
$D_R$--submodule since $H^c_{I(K\tensor_k R)}(K \tensor_k R)$ has a unique
simple $D_{{K \tensor_k R}/K}=K\tensor_k D_{R/k}$--module. This proves the
existence of the unique simple $D_R$--submodule of $H^c_I(R)$ which we
will denote by $\Ll(A,R)$.
\end{proof}
Along the way in this proof we also showed that $H^c_I(R)$ has a unique
simple $R[F^\infty]$--module $L(A,R)$:
\begin{theorem}
    Let\/ $R$ be regular, local and\/ $F$--finite. Let\/ $A=R/I$ of codimension\/ $c$
    be such that its completion is a domain. Then\/ $H^c_I(R)$ has a unique
    simple\/ $R[F^e]$-module\/ $L(A,R)$.
\end{theorem}

The above proof gives, \emph{a priori}, only limited information about the
concrete structure of $\Ll(A,R)$. The reason for this is that we use two
reduction steps. First, we reduce to the case that $R$ is complete in
order to apply the functor $\Dd$ and derive that $H^c_I(R)$ has a unique
simple $R[F^e]$--module. Secondly, we tensor with a huge algebraically
closed field in order to apply Main Theorem 1. Even though the unique
simple $R[F^\infty]$--submodule of $H^c_I(R)$ for $R$ complete is very
concretely constructed as $\Dd(H^d_m(R/I)/0^*_{H^d_m(R/I)})$, the
concreteness might get lost through the two reductions indicated. The aim
of the next section is to describe what information we still have
available about $\Ll(A,R)$ if $R$ is not complete.

\subsection{Construction of $\Ll(A,R)$}
We are looking for a more concrete description of $\Ll(A,R)$. We show that
the unique simple $R[F^\infty]$--submodule of the proof of Main Theorem 2
is obtained as a direct limit just as in the complete case. Thus we first
focus on the $R[F^e]$--structure and deal with the $D_R$--structure
separately. Let $L(A,R)$ denote the unique simple $R[F^e]$--submodule of
$H^c_I(R)$ to differentiate it from $\Ll(A,R)$, the unique simple
$D_R$--submodule.
\begin{proposition}\label{prop.LlcommutesComlpetion}
    Let\/ $(R,m)$ be regular, local and\/ $F$--finite. Let\/ $A=R/I$ be a domain.
    If\/ $L(A,R)$ denotes the unique simple\/ $R[F^\infty]$--submodule of\/ $H^c_I(A)$,
    then\/ $L(\wh{A},\wh{R}) \cong \wh{R} \tensor L(A,R)$.
\end{proposition}
\begin{proof}
In the complete case $L(\wh{A},\wh{R}) = \Dd(H^d_m(R/I)/{0^*_{H^d_m(A)}})$
and it arises as the direct limit
\[
    \tau_{\omega_{\wh{A}}} \to R^1 \tensor \tau_{\omega_{\wh{A}}} \to R^2 \tensor
    \tau_{\omega_{\wh{A}}} \ldots
\]
Furthermore, this description makes $\tau_{\omega_{\wh{A}}}$ the unique
minimal root of $L(\wh{A},\wh{R})$ since by Theorem \ref{thm.MaxRFofH}
$H^d_m(R/I)/{0^*_{H^d_m(A)}}$ is $F$--full and $F$--reduced and then by
Corollary \ref{cor.Dduniqueroot} its Matlis dual $\tau_{\wh{A}}$ is the
unique minimal root. Lemma \ref{lem.SubsofCompl} below shows that $N
\defeq H^c_I(R) \cap \Ll(\wh{A},\wh{R})$ is a unit $R[F^e]$--submodule of
$H^c_I(R)$. Thus a root of it is found by intersecting with the root
$\omega_A$ of $H^c_I(R)$, \ie $\omega_A \cap N = \omega_A \cap
\tau_{\omega_{\wh{A}}}=\tau_{\omega_A}$ is a root of $N$. By Lemma
\ref{lem.tauNonzero} $\tau_{\omega_A}$ is nonzero and thus $N$ is a
\emph{nonzero} unit $R[F^e]$--submodule of $H^c_I(R)$. Thus $\wh{R}
\tensor N$ is a nonzero unit $R[F^e]$--submodule of $L(\wh{A},\wh{R})$,
and since the latter is simple, they have to be equal. This implies that
also $N$ is simple, and therefore $N= L(A,R)$. Note that the root of
$\wh{R} \tensor L(A,R)=L(\wh{A},\wh{R})$ is $\wh{R} \tensor_R
\tau_{\omega_A} = \wh{A} \tensor_A \tau_{\omega_A}$. Since the
\emph{unique minimal} root of $L(\wh{A},\wh{R})$ is
$\tau_{\omega_{\wh{A}}}$, they must be equal, i.e. $\tau_{\omega_{\wh{A}}}
= \wh{A} \tensor_A \tau_{\omega_A}$.
\end{proof}
We still need the following general lemma about the behavior of unit
$R[F^e]$--modules under completion.
\begin{lemma}\label{lem.SubsofCompl}
    Let\/ $R$ be regular, local and\/ $F$--finite. Let\/ $H$ be a finitely generated
    unit\/
    $R[F^e]$--module and\/ $L$ a unit\/ $\wh{R}[F^e]$--submodule of\/ $\wh{R}
    \tensor H$. Then\/ $L \cap H$ is a unit\/ $R[F^e]$--submodule of\/ $H$.
\end{lemma}
\begin{proof}
By Proposition \ref{prop.LlcommutesComlpetion} it follows that the unit
$\wh{R}[F^e]$--module $\wh{R} \tensor H$ is in fact a unit
$R[F^e]$--module, and $L$ is also a unit $R[F^e]$--submodule of $\wh{R}
\tensor H$. Thus the intersection $L \cap H$ is indeed a unit
$R[F^e]$--submodule of $H$.
\end{proof}
As a corollary of the proof of Proposition \ref{prop.LlcommutesComlpetion}
we obtain that the parameter test module commutes with completion for
quotients of regular rings. Furthermore we have the concrete description
of $L(A,R)$ as the limit
\[
    \tau_{\omega_A} \to R^1 \tensor \tau_{\omega_A} \to R^2 \tensor
    \tau_{\omega_A} \to \ldots
\]
just as we have in the complete case.
\begin{corollary}\label{cor.TestCommCompl}
    Let\/ $A = R/I$ be a normal domain; $R$ regular, local and\/ $F$--finite. Then the
    \idx{parameter test module} commutes with completion. \ie
    $\tau_{\omega_{\wh{A}}} = \wh{A} \tensor_A \tau_{\omega_A}$.
\end{corollary}
This nicely complements a result of Smith \cite{Smith.test} on the
commutation of the parameter test module with localization for a complete
Cohen--Macaulay domain $A$.

The next step is to show that the unique simple unit
$R[F^\infty]$--submodule $L(A,R)$ is simple as a $D_R$--module; \ie
$L(A,R) = \Ll(A,R)$. This is clear whenever the cardinality of an
algebraically closed field $k$ contained in $R$ exceeds the cardinality of
$R$ as a $k$--vectorspace. Then Main Theorem 1 immediately shows that
$L(A,R)$ is $D_R$--simple. To show this in general we show first that
$L(A,R)$ behaves well under extension of the field $k$.
\begin{lemma}\label{lem.LlforFieldEx}
    Let\/ $R$ regular, local and\/ $F$--finite. Let\/ $A=R/I$ analytically irreducible. Let\/ $k
    \subseteq R$ be a perfect field and\/ $K$ be a perfect extension field of
    $k$. We denote\/ $K \tensor_k R$ by $R_K$ and similarly for\/ $A_K$.
    Then\/
    $K \tensor_k L(A,R) = L(R_K, A_K)$.
\end{lemma}
\begin{proof}
The proof is similar to the one for the complete case above, Thus, we only
give an outline. We show that $N \defeq L(R_K,A_K) \cap H^c_I(R) \subseteq
K \tensor H^c_I(R) = H^c_{IR_K}(R_K)$ is a nonzero unit $R[F^e]$--module.
For this, first observe that $K \tensor H^c_I(R)$ is a unit
$R[F^e]$--module by Proposition \ref{prop.LlcommutesComlpetion}. This
shows that $N$, as the intersection of two unit $R[F^e]$--submodules, is a
unit $R[F^e]$--submodule. It is nonzero by looking at its root
$\tau_{\omega_{A_K}} \cap \omega_A$, which by Lemma \ref{lem.tauFieldExt}
contains the nonzero $\tau_{\omega_A}$.

Then $K \tensor_k N$ is a nonzero unit $R_K[F^e]$--module contained in
$L(R_K, A_K)$ and thus equal to it. By flatness of $K \tensor_k \usc $ we
conclude that $N$ is simple and thus $N = L(A,R)$. Therefore, $K \tensor_k
L(A,R) = L(R_K,A_K)$.
\end{proof}
Now, the $D_R$--simplicity of $L(A,R)$ follows easily. Take $K$
algebraically closed of sufficiently big cardinality. Then $L(R_K, A_K) =
R_K \tensor L(A,R)$ is $D_{R_K}$--simple and therefore $L(A,R)$ is
$D_R$--simple by faithfully flatness of $R \to R_K$. Thus $L(A,R) =
\Ll(A,R)$. Combining these observations we get:
\begin{theorem}\label{thm.Llconcrete}
    Let\/ $A = R/I$ as in Main Theorem 2. Then the unique simple\/
    $D_R$--submodule\/ $\Ll(A,R)$ is the unique simple\/ $R[F^\infty]$--module
    of\/ $H^c_I(R)$. It arises naturally as the direct limit
    \[
        \tau_{\omega_A} \to R^1 \tensor \tau_{\omega_A} \to R^2 \tensor
        \tau_{\omega_A} \to \ldots
    \]
    where the map is the restriction of the natural map\/ $\omega_A \to R^1
    \tensor \omega_A$ dual to the\/ $R[F^e]$--module structure map\/ $R^1 \tensor H^d_m(R/I)
    \to[\theta] H^d_m(R/I)$ via local duality.
\end{theorem}

\section{$F$--rationality and simplicity of $H^c_I(R)$}\label{sec.FratSimpl}
The concrete description of the last section immediately allows us to
derive a characterization of the $D_R$--simplicity (equivalently,
u$R[F]$--simplicity) of $H^c_I(R)$.
\begin{theorem}\label{thm.HcIisSimpleCrit}
    Let\/ $R$ be regular local and\/ $F$--finite. Let\/ $I$ be an ideal such
    that\/
    $A=R/I$ is analytically irreducible. Then\/ $H^c_I(R)$ is\/ $D_R$--simple if
    and only if the tight closure of zero in\/ $H^d_m(A)$ is\/ $F$--nilpotent.
\end{theorem}
\begin{proof}
$H^c_I(R)$ is $D_R$--simple if and only if it is equal to $\Ll(A,R)$. By
Proposition \ref{prop.LlcommutesComlpetion} together with the fact that
$\Ll(A,R)=L(A,R)$, we can reduce to the complete case and thus assume that
$R$ and $A$ are complete. Then $L(A,R)=\Dd(H^d_m(A)/0^*_{H^d_m(A)})$ is
all of $H^c_I(R)$ if and only if $\Dd(0^*_{H^d_m(A)})=0$, by exactness of
$\Dd$. This is the case if and only if $0^*_{H^d_m(A)}$ is $F$--nilpotent
by Proposition \ref{prop.DisZeroOnNilpotent}.
\end{proof}
A trivial case where the latter condition is satisfied is when $A$ if
$F$--rational. In this case the tight closure of zero in $H^d_m(R)$ is
equal to zero. As a corollary we get:
\begin{corollary}
    Let\/ $R$ be regular, local and\/ $F$--finite. Let $I$ be an ideal such
    that\/
    $A=R/I$ is analytically irreducible. If\/ $A$ is\/ $F$--rational, then\/ $H^c_I(R)$
    is\/
    $D_R$--simple. If\/ $A$ is\/ $F$--injective (\ie $F$ acts injectively
    on\/
    $H^d_m(A)$), then\/ $A$ is\/ $F$--rational if and only if\/ $H^c_I(R)$
    is\/
    $D_R$--simple.
\end{corollary}
This should be compared to the following characterization of
$F$--regularity in terms of $D_A$--simplicity due to Smith:
\begin{proposition}[\protect{\cite[2.2(4)]{SmithDmod}}]
    Let\/ $A$ be an\/ $F$--finite domain which is\/ $F$--split. Then\/ $A$ is
    strongly\/ $F$--regular\index{F-regular@$F$--regular} if and only if\/ $A$ is simple as a\/
    $D_A$--module.
\end{proposition}
Note that this proposition is a statement about the $D_A$--module
structure of $A$, \ie a statement about the differential operators on $A$
itself. This is different from our approach as we work with the
differential operators $D_R$ of the regular ambient ring $R$.
Nevertheless, the similarity of the results are striking.

Roughly, this proposition is obtained by observing that a certain variant
$\wt{\tau}$ of the test ideal is a $D_A$--submodule of $A$. The variant of
$\tau$ we have in mind is the \emph{asymptotic} test ideal
\[
    \wt{\tau} = \bigcap_{I \subseteq A} \bigcap_{e \gg 0} (I^{[p^e]}:(I^*)^{[p^e]})
\]
whose elements work in tight closure tests for $e \gg 0$. If this ideal is
the unit ideal, then the Frobenius closure is equal to the tight closure
in $A$. This is clear since if $1 \in \wt{\tau}$, then $x \in I^*$ implies
that $1\cdot x^{p^e} \in I^{[p^e]}$ for some $e \gg 0$ and therefore $x
\in I^F$, by definition. Thus if $A$ is $D_A$--simple, then $\wt{\tau} =
A$ since $\wt{\tau}$ is a $D_A$--submodule of $A$. Therefore $I^F=I^*$ for
all ideals $I$ of $A$. Picking up on this analogy we can show the
following.
\begin{proposition}\label{prop.HcIsimpleCrit}
    Let\/ $R$ be regular, local and\/ $F$--finite. Let\/ $I$ be an ideal such
    that\/
    $A=R/I$ is analytically irreducible. If for all parameter ideals of\/ $A$ we
    have\/
    $J^F=J^*$, then\/ $H^c_I(R)$ is\/ $D_R$--simple.

    If\/ $A$ is Cohen--Macaulay, then\/ $H^c_I(R)$ is\/ $D_R$--simple if and only
    if\/
    $J^*=J^F$ for all parameter ideals\/ $J$.
\end{proposition}
\begin{proof}
We show that if $J^*=J^F$ for all parameter ideals, then $0^*_{H^d_m(A)}$
is $F$--nilpotent, \ie $0^*_{H^d_m(A)}=0^F_{H^d_m(A)}$. Let $\eta \in
H^d_m(A)$ represented by $z + (\xd)$ for some parameter ideal $J = (\xd)$,
thinking of $H^d_m(A)$ as the limit $\dirlim A/J^{[p^e]}$. Then the colon
capturing property of tight closure shows that $z \in 0^*_{H^d_m(A)}$ if
and only if $z \in J^*$ (cf.\ \cite[Proposition 3.1.1 ]{SmithDiss}). By
our assumption $J^*=J^F$, this implies that $z^{p^e} \in J^{[p^e]}$ for
some $e > 0$. Consequently, $F^e(\eta)=z^{p^e}+J^{[p^e]}$ is zero and thus
every element of $0^*_{H^d_m(A)}$ is $F$--nilpotent.

Under the assumption that $A$ is Cohen--Macaulay the same argument can be
reversed using that the limit system defining $H^c_m(A)$ is injective.
\end{proof}

\section{Examples and Applications}
Applying the extension of Matlis duality $\Dd$ to the lattice of
$R[F^e]$--submodules of $H^d_m(R/I)$, we obtain a complete description of
the lattice of unit $R[F^e]$--sub\-mo\-du\-les of $H^c_I(R)$. In this
section we carry this out in the case of a few examples. We focus on the
case of a graded complete intersection with at worst isolated non
$F$--rational locus. In this case $0^*_{H^d_m(A)} \subseteq H^d_m(A)$ is
well understood. The following was conjectured by Smith and Huneke in
\cite{HuSmith.KodVan}, and proven by Hara \cite[6.1]{Hara.GeomTight}
\begin{proposition}\label{prop.TightisPosDeg}
    Let\/ $(A,m)$ be a graded\/ $k$--algebra with\/ $k$ perfect which is obtained by
    reduction from characteristic zero. Assume that\/ $A$ is\/ $F$--rational
    away from the graded maximal ideal\/ $m$, then for\/ $p \gg 0$ the tight
    closure of zero in\/ $H^d_m(A)$ consists exactly of the elements of
    non--negative degree.
\end{proposition}
For the technique of \idx{reduction to finite characteristic} see
\cite{Smith.sing} or Section \ref{sec.relcharzero} below. Since we are
working with a graded quotient of a graded regular ring assuming that $R$
arises from a ring of characteristic zero is not really a serious
restriction. So as long as we keep in mind to let the characteristic be
sufficiently large the above proposition will apply.

\subsection{Complete Intersections and Gorenstein rings}
If the local ring $A=R/I$ is Gorenstein then the canonical module
$\omega_A$ is isomorphic to $A$ itself. Therefore, the sequence
\[
    \Ext^c(R/I,R) \to \Ext^c(R/I^{[p^e]},R) \to \Ext^c(R/I^{[p^e]},R) \to
    \ldots
\]
calculating $H^c_I(R)$ is, as we noticed before in
\eqnref{eqn.HcIviaomega} on page \pageref{eqn.HcIviaomega}, via Matlis
duality isomorphic to the sequence
\begin{equation}\label{eqn.HcIforGor}
    R/I \to R/I^{[p^e]}  \to R/I^{[p^{2e}]} \to R/I^{[p^{3e}]} \to \ldots\
    .
\end{equation}
Under this identification, the map is given by multiplication by the
$A$--module generator of $\frac{(I^{[p^e]}:I)}{I^{[p^e]}}$. Recall that by
Section \ref{sec.InjHullEx} this is indeed a free $A$--module of rank one
and the canonical generator is obtained via a choice of an isomorphism
$H^d_m(A) \cong E_A$. Let us fix a representative $y \in R$ of the
$A$--module generator of $\frac{(I^{[p^e]}:I)}{I^{[p^e]}}$. Then the maps
in Sequence \eqnref{eqn.HcIforGor} are just multiplication by $y$,
$y^{p^e}$ and so forth. For example, if $A$ is a complete intersection in
$R$, say $I=(\xc)$ is an ideal generated by a regular sequence of elements
of $R$, then
\[
    \frac{(I^{[p^e]} : I)}{I^{[p^e]}} = \frac{(x_1\cdot \ldots \cdot
    x_c)^{p^e-1}}{I^{[p^e]}}.
\]
Therefore $y = x^{p^e-1} \defeq (x_1 \cdot \ldots \cdot x_c)^{p^e-1}$ is a
sensible choice for $y$. By Theorem \ref{thm.Llconcrete}, $\Ll(A,R)$
arises by restricting the directed system \eqnref{eqn.HcIforGor} to the
parameter test ideal. Therefore, if we denote the pullback of the
parameter test ideal $\wt{\tau}_A\subseteq A$ to $R$ by $\tau$, we get
\[
    \Ll(A,R) = \dirlim( \tau/I \to[y] \tau^{[p^e]}/I^{[p^e]} \to[y^{p^e}]
    \tau^{[p^{2e}]}/I^{[p^{2e}]}\to \ldots \ ).
\]

In the graded or complete case the same holds true for the annihilator in
$R$ of \emph{any} $R[F^e]$--submodule of $H^d_m(A)$. If $M  \subseteq
H^d_m(A)$ is such an $R[F^e]$--submodule, Matlis Duality shows that
$\Ann_A(M) \cong D(H^d_m(A)/M) \subseteq D(H^d_m(A)) \cong A$ (we assume
$A$ Gorenstein, otherwise replace $A$ by $\omega_A$). If we denote by
$\tau_M$ the pullback of $\Ann_A M$ to $R$, \ie $\tau_M = \Ann_R M$, then
Proposition \ref{prop.DdGraded} shows that the map $R/I \to[\cdot y]
R/I^{[p^e]}$ restricts to a map $\tau_M \to[\cdot y]
\tau_M^{[p^e]}/I^{[p^e]}$, and
\[
    \Dd(M) \cong \dirlim(\tau_M \to[\cdot y] \tau_M^{[p^e]}/I^{[p^e]} \to[\cdot y^{p^e}]
    \tau_M^{[p^{2e}]}/I^{[p^{2e}]}\to \ldots \ ).
\]
This leads to an interesting observation about colon ideals of
annihilators of $R[F^e]$--submodule of $H^c_I(R)$.
\begin{proposition}
    Let\/ $A=R/I$ be a Gorenstein quotient of the complete (or graded) regular local ring\/ $R$.
    For an\/ $R$--submodule\/ $M \subseteq H^d_m(A)$ we denote by\/ $\tau_M = \Ann_R(M)$
    its annihilator in\/ $R$. Then\/ $M$ is an\/ $R[F^e]$--submodule if and only
    if\/ $(I^{[p^e]} :_R I) \subseteq (\tau_M^{[p^e]} :_R \tau_M)$.
\end{proposition}
\begin{proof}
One implication we showed above. Namely, if $M$ is an $R[F^e]$--submodule
then the $A$--module generator $y$ of $(I^{[p^e]}:I)/I^{[p^e]}$ induces a
map $\tau_M/I \to \tau_M^{[p^e]}/I^{[p^e]}$. Thus $(I^{[p^e]}:I)\tau_M
\subseteq \tau_M^{[p^e]}$ as claimed.

Conversely, if we have an inclusion of colon ideals as above, we get that
multiplication by $y$ induces a map $\tau_M/I \to
\tau_M^{[p^e]}/I^{[p^e]}=R^e \tensor \tau_M/I$. This map is the root of
some unit $R[F^e]$--submodule $\Mm$ of $\Dd(H^d_m(A))$, and by Theorem
\ref{thm.DdofSimpisSimp} (more precisely its proof) it arises as $\Dd(M')$
for $M'= D(\Mm \cap A)$, an $R[F^e]$--module quotient of $H^d_m(A)$. Here
we use that $A \to[\cdot y] R^e \tensor A$ is a root of $H^d_m(R)$ and
therefore $\Mm \cap A = \tau_M/I$ and thus $M'=D(\Mm \cap A) = D(\tau_M/I)
= M$. Therefore $M=M'$ is an $R[F^e]$--submodule of $H^d_m(R)$.
\end{proof}
This last lemma is a special case of \cite[Proposition
5.2]{LyubSmith.Comm} where an equivalent fact is established for ideals of
$R$ whose annihilator in the injective hull is Frobenius stable, thus in
the Gorenstein case, they amount to the same statement.

\subsection{Some Concrete Examples}
With the insights developed in the last section we are able to describe
the structure of $H^c_I(R)$ as a unit $R[F^e]$--module (\ie the lattice of
unit $R[F^e]$--submodules) by investigating the structure of $H^d_m(A)$ as
an $R[F^e]$--module. We focus on the graded case and therefore make use of
our graded variant of $\Dd$ developed in Section \ref{sec.GradedDd}.

\begin{example}[Graded CI with zero $a$--invariant]
Let $A=k[\xn[y]]/I$ for $I=(\xc)$ where the $x_i$'s are a homogeneous
regular sequence and the sum of their degrees is $n$. We give a criterion
for when $H^c_I(R)$ is $D_R$--simple.
\begin{lemma}
    Let\/ $A=R/I$ with\/ $R=k[\xn[y]]$ a graded complete intersection
    with\/
    $I=(\xc)$ and\/ $\deg (x) = \deg(x_1 \cdot \ldots \cdot x_c) = n$ (\ie
    the\/ $a$--invariant of\/ $A$ is zero). Assume that\/ $A$ is\/ $F$--rational away from
    the graded maximal ideal of\/ $m$ and that the characteristic is sufficiently big.
    Then\/ $A$ is\/ $F$--injective if and only if\/ $H^c_I(R)$ is\/ $D_R$--simple, if and only
    if the coefficient of\/ $(y_1 \cdot \ldots \cdot y_n)^{p-1}$ in\/
    $x^{p-1}$ is nonzero.
\end{lemma}
\begin{proof}
The assumptions on $A$ allow us to apply Proposition
\ref{prop.TightisPosDeg}. Therefore, the tight closure of zero in
$H^d_m(A)$ is exactly its part of degree zero, \ie $0^*_{H^d_m(A)}$ is
precisely the one dimensional socle of $H^d_m(A)$. Thus, the parameter
test ideal $\tau_A$ is the graded maximal ideal of $A$. Its pullback to
$R=k[\xn[y]]$ is just the maximal ideal $m=(\xn[y])$ of $R$. Now,
$\Ll(A,R)=\Dd(H^d_m(A)/0^*_{H^d_m(A)}) = H^c_I(R)$ if and only if
$0^*_{H^d_m(A)}$ is $F$--nilpotent. This in turn is the case if and only
if $R^1 \tensor 0^*_{H^d_m(A)} \to 0^*_{H^d_m(A)}$ is the zero map since
$0^*_{H^d_m(A)}$ is one dimensional. This last map is the dual of the
generator $\beta: D(0^*_{H^d_m(A)}) \to R^1 \tensor D(0^*_{H^d_m(A)})$,
which by the identifications of the preceding paragraph is just the map
\[
    \beta: R/m \to[\cdot x^{p-1}] R/m^{[p]}.
\]
It is given by multiplication by the $(p-1)$st power of $x = x_1 \cdot
\ldots \cdot x_c$. This map is zero if and only if $x^{p-1} \in m^{p}$.
Since the degree of $x^{p-1}$ is $n(p-1)$, the only monomial of this
degree which is not in $m^{[p]}$ is $(y_1 \cdot \ldots \cdot y_n)^{p-1}$.
\end{proof}
This, of course, also follows from Fedder's criterion of $F$--purity
\cite{Fedder83}. In the case that $A$ is $F$--injective, the above
argument also shows that
\[
    \Dd(0^*_{H^d_m(A)}) \cong \dirlim( R/m \to[\cdot y^{p-1}] R/m^{[p]} \to[\cdot y^{p(p-1)}]
    R/m^{[p^2]} \to \ldots \ )
\]
where $y$ denotes the product of all the $y_i$ (up to a scalar multiple).
But this sequence is isomorphic (as an $R[F]$--module) to $H^n_m(R)$. Thus
the short exact sequence of $R[F]$--modules
\[
    0 \to 0^*_{H^d_m(A)} \to H^d_m(A) \to H^d_m(A)/0^*_{H^d_m(A)} \to 0
\]
gives after applying the exact functor $\Dd$ the short exact sequence of
unit $R[F^e]$--modules:
\[
    0 \to \Ll(A,R) \to H^c_I(R) \to H^n_m(R) \to 0
\]
In terms of the Koszul complex description of the local cohomology modules
the map on the right is just given by sending $r + I^{[p^e]}$ to
$r+m^{[p^e]}$. So $H^c_I(R)$ is a unit $R[F]$--module extension of the two
simple unit $R[F]$--module $\Ll(A,R)$ and $H^n_m(R)$. Thus $\Ll(A,R)$ is
concretely represented as the kernel of the map from $H^c_I(R) \to
H^n_m(R)$ as just described.
\end{example}
\begin{example}[Fermat Hypersurfaces]
Now let $A = \frac{k[\xn[y]]}{y_1^d+\cdots+y_{n-1}^d-y_n^d}$ and denote by
$R$ the polynomial ring $R = k[\xn[y]]$ and
$I=(y_1^d+\cdots+y_{n-1}^d-y_n^d)$. Then $A$ has an isolated singularity
at the graded maximal ideal and therefore the tight closure of zero in
$H^d_m(A)$ is just its part of nonnegative degree (cf.\ Proposition
\ref{prop.TightisPosDeg}). Since the action of the Frobenius on $H^d_m(A)$
multiplies degrees by $p$, the part of $H^d_m(A)$ of strictly positive
degree is $F$--nilpotent. Thus, for the purpose of $\Dd$, which vanishes
on $F$--nilpotent $R[F^e]$--modules we can ignore the strictly positive
part and concentrate on the part of degree zero. We have the following.
\begin{lemma}\label{lem.FermHypInjec}
    Let\/ $A=R/I$ be the \idx{Fermat hypersurface} described above. If\/ $p = 1 \mod{d}$, then
    the Frobenius acts injectively on the degree zero piece of\/ $H^d_m(R)$.
\end{lemma}
\begin{proof} If we compute $H^d_m(A)$ via the \Czech complex arising
from the system of parameters $(y_1, \ldots, y_{n-1})$, then a $k$ basis
of the degree zero part of $H^{n-1}_m(A)$ is given by $e_{s} =
\frac{y_n^{d-1}}{y_1\cdots y_{n-1}} \cdot s^{-1}$ where $s$ ranges through
all monomials of degree $d-n$ in $\xn[y]$. To determine the image of such
an element $e_{s}=\frac{y_n^{d-i_n}}{y_1^{i_1}\cdot\ldots\cdot
y_{n-1}^{i_{n-1}}}$ (then $s = y_1^{i_1-1}\cdot \ldots \cdot y_n^{i_n-1}$
and $i_j > 1$ and $\sum i_j = d$) under the Frobenius one just calculates
\[
\begin{split}  F(e_s) &= \frac{y_n^{p(d-i_n)}}{y_1^{pi_1}\cdot\ldots\cdot y_{n-1}^{pi_{n-1}}}
                       = y_n^{d-i_n}\frac{y_n^{(p-1)(d-i_n)}}{y_1^{pi_1}\cdot\ldots\cdot
                       y_{n-1}^{pi_{n-1}}} \\
                      &= y_n^{d-i_n}\frac{(y_1^d+\ldots+y_{n-1}^d)^{r(d-i_n)}}{y_1^{pi_1}\cdot\ldots\cdot
                       y_{n-1}^{pi_{n-1}}} \\
                      &= y_n^{d-i_n}\sum_{j_1+ \cdots+
                       j_{n-1}=r(d-i_n)}\frac{(r(d-i_n))!}{j_1!\cdot \ldots \cdot j_{n-1}!}
                       \frac{y_1^{dj_1}\cdot \ldots \cdot
                       y_{n-1}^{dj_{n-1}}}{y_1^{pi_1}\cdot\ldots\cdot y_{n-1}^{pi_{n-1}}}
\end{split}
\]
where we wrote $\frac{p-1}{d}=r$. For a term in this sum to be nonzero in
the local cohomology module $H^d_m(A)$, we must have that the exponent of
each $y_k$ in the denominator is strictly bigger than the one in the
numerator, \ie $dj_k < pi_k$ for all $k=1, \ldots, (n-1)$. Using $p-1=rd$
this amounts to $dj_k < dri_k + i_k$. This shows that the only surviving
term in this sum is the one corresponding to $j_k = ri_k$. Thus we get
\[
    F(e_s) =  \frac{(r(d-i_n))!}{(ri_1)!\cdot \ldots \cdot
    (ri_{n-1})!}
    \frac{y_n^{d-i_n}}{y_1^{i_1}\cdot\ldots\cdot y_{n-1}^{i_{n-1}}} = \frac{(r(d-i_n))!}{(ri_1)!\cdot \ldots \cdot
    (ri_{n-1})!} e_s
\]
and point out that the coefficient is nonzero since it is not divisible by
$p$ as the term inside the factorial expression of the numerator is
$r(d-i_n)=(p-1)-ri_n$ which is smaller than $p$. This shows that the
Frobenius acts injectively on the degree zero part of $H^d_m(A)$.
\end{proof}
It follows that for $p = 1 \mod{d}$ the Frobenius reduced part of
$0^*_{H^d_m(A)}$ is exactly its degree zero part. The $F$--reduced
quotient is therefore a $k[F]$--module and the matrix representing the
action of the Frobenius as in Section \ref{sec.FreeRF} with respect to the
basis $\{e_s\}$ is given by a diagonal matrix (since $F(e_s) =
\text{constant}\cdot e_s$ from the last proof). Thus, as a $k[F]$--module
(equivalently $R[F]$--module), $\Fred{(0^*_{H^d_m(A)})}$ decomposes into a
direct sum of simple $k[F]$--modules, namely $\Fred{(0^*_{H^d_m(A)})}
\cong \oplus ke_s$ is a direct sum of $k[F]$--modules.

By definition $\Ll(A,R)$ is the kernel of the map
$H^c_m(A)=\Dd(H^{n-1}_m(A)) \to
\Dd(0^*_{H^d_m(A)})=\Dd(\Fred{(0^*_{H^d_m(A)})})$. And the decomposition
of $\Fred{(0^*_{H^d_m(A)})}$ into the simple $R[F]$--submodules $ke_s$
implies that $\Ll(A,R)$ is, in fact, the intersection of the kernels of
the maps $\Dd(H^{n-1}_m) \to \Dd(ke_s)$.

As a final investigation in this example we describe these maps
concretely. From our setup, it follows that if
$s=y_1^{i_1-1}\cdot\ldots\cdot y_n^{i_n-1}$ is a monomial of degree $d-n$,
then the annihilator of $e_s\defeq\frac{y_n^{d-1}}{y_1\cdots y_{n-1}}
\cdot s^{-1}$ in $R$ is just the parameter ideal $\tau_s \defeq
(y_1^{i_1},\ldots,y_n^{i_n})$. Thus we can identify the Matlis dual of the
inclusion $Re_s \subseteq H^{n-1}_m(A)$ with the natural projection $R/I
\to R/\tau_s$. Furthermore, the above proof shows that under this
identification the map $D(Re_s) \to R^1 \tensor D(Re_s)$ is (up to
constant multiple) given by
\[
    \beta: R/\tau_{e_s} \to[\cdot (y_1^{i_1} \cdot\ldots\cdot
    y_n^{i_n})^{p-1}] R/\tau_{e_s}^{[p]}.
\]
The limit of the directed system generated by this map is
$\Dd(Re_s)=\Dd(ke_s)$. At the same time, this directed system computes
$H^n_{m}(R)$ as a limit of Koszul complexes arising from the system of
parameters $(y^{i_1}_1, \ldots, y^{i_n}_n)$. With this identification we
get a series of maps
\[
    H^{n-1}_I(R) \to[\phi_s] H^n_m(R) \cong E_R
\]
where $\phi_s$ sends $r+I^{[p^e]}$ to $r+\tau_s^{[p^e]}$ when thinking of
both local cohomology modules as arising from limits of Koszul homology as
described before. As we already pointed out, the intersection of the
kernels of these maps is $\Ll(A,R)$.
\end{example}

\begin{example}
Concretely, we consider $A=\frac{k[x,y,z]}{x^4+y^4-z^4}$. The degree zero
part of $H^2_m(A)$ has a $k$--basis represented by the elements
\[
\frac{z^3}{x^2y}\ ,\  \frac{z^3}{xy^2}\ ,\ \frac{z^2}{xy}
\]
and Lemma \ref{lem.FermHypInjec} shows that for $p=1 \mod{4}$ the
Frobenius acts injectively on each of them. A simple calculation shows
that in the case $p=3 \text{ mod } 4$ the Frobenius acts as zero. We show
this for the first of the three generators, for the other two one shows
this analogously.
\[
F(\frac{z^3}{x^2y}) = \frac{z^{3p-1}z}{x^{2p}y^p}=z
\frac{\sum_{4i+4j=3p-1} c_{ij} x^{4i}y^{4j}}{x^{2p}y^p}
\]
The degree of the numerator is $3p-1$ which implies that every monomial is
in the ideal generated by $(x^{2p},x^p)$ which implies that the above sum
is zero. Thus, if $p = 3 \mod{4}$ the Frobenius acts as zero on
$0^*_{H^d_m(A)}$, thus $\Dd(0^*_{H^d_m(A)})=0$ and therefore $\Ll(A,R) =
H^1_I(R)$ is a simple unit $R[F]$--module.
\end{example}

\begin{example}
The next harder case is $R=\frac{k[x,y,z]}{x^5+y^5-z^5}$. Similar
calculations can be made as in the previous example. One finds that for $p
\neq 1 \mod{5}$ the square of each of the generators of the degree zero
part vanishes and therefore $H^1_{(x^5+y^5-z^5)}(k[x,y,z])$ is a simple
unit $R[F^e]$--module. For $p = 1 \mod{5}$ there are $6$ maximal
non--trivial $R[F]$--submodules of $H^d_m(R)$, one for each generator of
the degree zero part of $H^1_I(R)$. Each of these $R[F]$--submodules is
the kernel of a map $H^1_I(R) \to H^3_m(R)=E_R$ and their intersection is
$\Ll(A,R)$.
\end{example}

\subsection{On a result of S.P.Smith}
Finally we show that if $A=R/I$ is a domain such that its normalization
$\bar{A}$ if $F$--rational and the natural projection $\Spec{\bar{A}} \to
\Spec{A}$ is injective, then $H^c_I(R)$ is $D_R$--simple. This can be
viewed as a generalization of the results of S.P. Smith
\cite{SmithSP.IntHom} in characteristic zero showing that if $f=0$ defines
a plane cusp, then $H^1_{(f)}(k[x,y])$ is $D_{k[x,y]}$--simple. This also
follows from the Riemann--Hilbert correspondence but \cite{SmithSP.IntHom}
gives a completely algebraic (elementary) proof of the result.

\begin{proposition}
    Let\/ $R$ be regular local and\/ $F$--finite. Let\/ $A=R/I$ be a domain
    with isolated singularity such that the normalization map\/
    $\Spec{\bar{A}} \to \Spec{A}$ is injective
    and\/ $\bar{A}$, the normalization itself, is\/ $F$--rational.
    Then\/ $H^c_I(A)$ is\/ $D_R$--simple.
\end{proposition}
\begin{proof}
We plan to apply Theorem \ref{thm.HcIisSimpleCrit}. For this one has to
show that, first, the completion of $A$ is a domain, and secondly, that
the tight closure equals the Frobenius closure for parameter ideals. The
assumption on the injectivity of the normalization map ensures that the
normalization $\bar{A}$ of $A$ is also a local ring. Furthermore, the
completion of $\bar{A}$ along its unique maximal ideal is just $\wh{A}
\tensor \bar{A}$, by module finiteness of $\bar{A}$ over $A$. By
faithfully flatness of completion we get that $\wh{A}$ is a subring of
$\wh{\bar{A}}$. Since $\bar{A}$ is a normal, excellent (implied by
$F$--finiteness) domain its completion is a domain. Therefore, $\wh{A}$ is
a domain.\note{refer to Matsumura}

Let $z \in J^*$ for a parameter ideal $J=(\xd[y])$ of $A$. Then, since
$\bar{A}$ is $F$--rational and the expansion $\bar{J}$ of $J$ to $\bar{A}$
is also a parameter ideal, one concludes that $z \in \bar{J}^* = \bar{J}$.
Let, for some $a_i \in \bar{A}$,
\[
    z = a_1y_1+ \ldots +a_dy_d
\]
be an equation witnessing this ideal membership. As we observe in Lemma
\ref{lem.IntClos} below, for some big enough $e$, all $a^{p^e}$ are in
$A$. Therefore $z^{p^e} = a_1^{p^e}y_1^{p^e}+ \ldots +a_d^{p^e}y_d^{p^e}$,
which shows that $z^{p^e} \in J^{[p^e]}$ since all $a^{p^e} \in A$. Thus
$J^* = J^F$ and Theorem \ref{thm.HcIisSimpleCrit} implies that $H^c_I(R)$
is $D_R$--simple.
\end{proof}
\begin{lemma}\label{lem.IntClos}
    Let\/ $A$ be as in the last proposition and let\/ $x \in \bar{A}$ be an
    element of its normalization. Then\/ $x^n \in A$ for some\/ $n \in \NN$.
\end{lemma}
\begin{proof}
The condition that the normalization map is injective means that for every
prime ideal $q$ of $A$ there is exactly one prime ideal $Q$ of $\bar{A}$
lying over $q$. This, in turn, implies that the radical of the expansion
$q\bar{A}$ of $q$ to $\bar{A}$ is just $Q$. As a formula:
$\sqrt{q\bar{A}}=Q$.

Keeping this in mind let $x\in \bar{A}$ a nonunit. Then $x$ lies inside
some prime ideal $Q$ of $\bar{A}$. Let us denote $q=Q\cap A$. Then
sufficiently high powers of $x$ are in $q\bar{A}$ since
$\sqrt{q\bar{A}}=Q$. We want to conclude that $x^n \in q$ and thus is in
$A$ for big enough $n$. For this we show that the conductor ideal $C=
(A:_{\bar{A}}\bar{A})$ contains sufficiently high powers of $x$. Since $A$
is assumed to have isolated singularities we observe that the conductor
$C=C\bar{A}$ is primary to the maximal ideal. Thus for every nonunit $x
\in A$ sufficiently high powers lie in the conductor $C$. Now, if $x^n \in
q\bar{A}$ and $x^m \in C$, then $x^{n+m}=x^nx^m$ is in $q$ itself.

Now let $u$ be a unit of $\bar{A}$. In the fraction field of $A$ we can
write $u=\frac{ux}{x}$ for $x$ a nonunit of $\bar{A}$. Since both, $x$ and
$ux$ are not units, some power is in $A$. Therefore the same power of $u$
will be in $A$ too. This finishes the proof.
\end{proof}

%%%%%%%%%%%%%%%%%%%%%%%%%%%%%%%%%%%%%%%%%%%%%%%%%%%%%%%%%%%%%%
%%                                                          %%
%%   This is file: chapter6.tex                             %%
%%   It contains the Fifth Chapter of my                    %%
%%   dissertation: dissertation.tex                         %%
%%                                                          %%
%%%%%%%%%%%%%%%%%%%%%%%%%%%%%%%%%%%%%%%%%%%%%%%%%%%%%%%%%%%%%%

\chapter{Problems}\label{chap.Problems}
The results in this dissertation give an (almost) complete answer in the
local case to the question of constructing the analog of the intersection
cohomology $D_R$--module in finite characteristic. We show that, under
reasonable assumptions, there is such an analog. Not only do we show its
existence, but perhaps more importantly, we give a concrete construction
of $\Ll(A,R)$, tying it to the tight closure theory on $A$. In this
respect our finite characteristic description of $\Ll(A,R)$ is superior to
the characteristic zero picture, since there only for a much smaller class
(complete intersections with isolated singularities \cite{Vil}) a concrete
description is known.

Of course, many questions about $\Ll(A,R)$ in finite characteristic still
need to be answered. Below we discuss the ones we consider most
interesting.
\section{Globalization}
To answer Question \ref{ques.LnFiniteChar} in the generality it is posed
in the introduction the most pressing issue is that of globalization; \ie
we are asking whether the following is true.\footnote{As Brian Conrad
points out, a smearing out technique a la EGA \cite{EGA4.3} might be
successfully applied here too. If $R$ denotes the ring of a neighborhood
of a point $x$ of $X$ the technique yields (after possibly shrinking that
neighborhood) a finitely generated unit $R[F^e]$--submodule $\Ll$ of
$H^c_I(R)$ such that $\Ll_x = \Ll(A_x,R_x)$. With this at hand, the
globalization problem comes down to showing that for other points $y$ in
this neighborhood of $x$ we have $\Ll_y=\Ll(A_y,R_y)$. A treatment of this
is in preparation.}
\begin{problem}\label{prob.Global}
    Let\/ $Y \subseteq X$ be a closed irreducible subscheme of the smooth
    scheme\/ $X$. For each point\/ $x \in X$ we showed the existence of a
    unique simple\/ $D_{\Oo_{x,X}}$--submodule\/ $L(\Oo_{x,Y},\Oo_{x,X})$
    of\/ $\Hh^c_{[Y]}(\Oo_X)_x$. Is there a unique simple $\Dd_X$--subsheaf\/ $\Ll(Y,X)$
    of\/ $\Hh^c_{[Y]}(\Oo_X)$ such that for all\/ $x \in X$,
    $\Ll(Y,X)_x=L(\Oo_{x,Y},\Oo_{x,X})$?
\end{problem}
The line of reasoning followed in this dissertation doesn't immediately
lend itself to these kind of globalization questions. Several crucial
steps of our construction are local in nature. First, we use the tight
closure of zero in $H^d_m(A)$ as the starting point in our construction of
$L(R,A)$. Secondly, local duality, in the form of its extension $\Dd$, is
used to relate the tight closure $0^*_{H^d_m(A)}$ to the simple
$R[F^\infty]$--module $L(A,R)$. Both of these constructions are not
available in a non local setting.

In an attempt to solve the above problem we first consider the affine
case. For $Y=\Spec A$ and $X=\Spec R$ we have to show that $H^c_I(R)$ has
a $D_R$--submodule $L(A,R)$ such that for all maximal ideals $m \in R$,
$L(A,R)_m = L(A_m,R_m)$. Having found such $L(A,R)$ we still have to
verify it is in fact the unique simple $D_R$--submodule of $H^c_I(R)$. As
it turns out, this is not so difficult.

We try to mimic the concrete construction of $L(A_m,R_m)$ as in Theorem
\ref{thm.Llconcrete}. There we show that $L(A_m,R_m) \subseteq H^c_I(R_m)$
arises as a sublimit of the direct limit computing $H^c_I(R_m)$ itself:
\[
    H^c_I(R_m)  =  \dirlim( \omega_{A_m} \to R^1 \tensor \omega_{A_m}
    \to R^2 \tensor \omega_{A_m} \to \ldots \ \ )
\]
where the maps are ultimately induced by the natural projection $R/I^{[p]}
\to R/I$. $L(A_m,R_m)$ is then obtained from the Frobenius powers of the
submodule $\tau_{\omega_{A_m}} \subseteq \omega_{A_m}$.

Attempting to find $L(A,R)$ similarly in the non--local situation we
define a global version of the parameter test module. For this we recall
the global version of the test module $\omega_A = \Ext^c(A,\omega_R)$. As
the canonical module of $R$ we take the highest exterior power of the
sheaf of K\"ahler differentials $\omega_R = \bigwedge^n \Omega_R$.
\begin{definition}
    Let $A$ be a ring with canonical module $\omega_A$. Then the
    \emph{\idx{parameter test module}} $\tau_{\omega_A}$ is defined as
    \[
        \tau_{\omega_A} \defeq \bigcap \omega_A \cap \tau_{\omega_{A_m}}
    \]
    where the intersection ranges over all maximal ideals (equivalently
    all prime ideals) of $A$.
\end{definition}
This seems to be a sensible definition as it agrees with the definition of
the parameter test module given earlier in the case that $A$ is local. The
key step for showing that the parameter test module is the correct object
(\ie it will give rise to the sought after $L(A,R)$ in the affine case) is
to show that it is well behaved under localization, \ie we have to solve
the following problem.
\begin{problem}\label{prob.TestCommutes}
    Let $A$ be a ring with canonical module $\omega_A$. Let $S \subseteq
    A$ be a multiplicatively closed subset of $A$. Is
    $\tau_{\omega_{S^{-1}A}}$ equal to  $S^{-1}\tensor_A \tau_{\omega_A}$?
\end{problem}
We show that the existence (and construction) of the unique simple
$R[F^\infty]$--submodule of $H^c_I(R)$ in the non--local case follows, if
we assume that the assertion of Problem \ref{prob.TestCommutes} is true
for $S=(R\setminus m)$, $m$ a maximal ideal of $R$.
\begin{proposition}\label{prop.Llglobal}
    Let\/ $R$ be regular and of finite type over a field. Let\/ $R/I=A$ be a
    normal domain. Assume that for all maximal ideals of\/ $R$ we have
    that\/
    $R_m \tensor \tau_{\omega_R} = \tau_{\omega_{R_m}}$. Then\/
    $\tau_{\omega_A} \subseteq F^e(\tau_{\omega_A})$
    as submodules of\/ $H^c_I(R)$. Furthermore, the unit\/ $R[F^e]$--module generated
    by\/
    $\tau_{\omega_A}$ is\/ $L(A,R)=\bigcup F^{er}(\tau_{\omega_R})$ is the unique
    simple\/
    $R[F^e]$--submodule of\/ $H^c_I(R)$ and\/ $R_m \tensor L(A,R) =
    L(A_m,R_m)$.
\end{proposition}
\begin{proof}
The inclusion
\begin{equation}\label{eqn.omegaisFann}
    \tau_{\omega_R} \subseteq F^e(\tau_{\omega_R})
\end{equation}
holds if and only if it holds after localizing at maximal ideals, \ie if
and only if for all maximal ideals $m$ of $R$ we have $R_m \tensor
\tau_{\omega_R} \subseteq R_m \tensor F^e(\tau_{\omega_R})$. Using the
assumption that the parameter test module commutes with localization at
maximal ideals\footnote{In \cite[Remark 3.2 (ii)]{Hara.GeomTight} Hara
claims that it is easy to see that $(\tau_{\omega_A})_m =
\tau_{\omega_{A_m}}$. To the best of our knowledge though, it is an open
question. If $A$ is Gorenstein, it is equivalent to the question whether
the test ideal commutes with localization, which, is still open in the
case that $A$ is not local.}, as well as the fact that the Frobenius
commutes with localization, this is just the inclusion
$\tau_{\omega_{R_m}} \subseteq F^e(\tau_{\omega_{R_m}})$ which was proven
in Theorem \ref{thm.Llconcrete}. The inclusion \eqnref{eqn.omegaisFann}
allows to define $L(A,R)$ to be the unit $R[F^e]$--submodule of $H^c_I(R)$
generated by $\tau_{\omega_R} \subseteq H^c_I(R)$. Clearly, localizing at
a maximal ideal we get
\[
    R_m \tensor L(A,R)=R_m \tensor \dirlim(F^{er}(\tau_{\omega_R}))
    = \dirlim F^{er}(\tau_{\omega_{R_m}})=L(A_m, R_m).
\]
It remains to show that $L(A,R)$ is the \emph{unique simple} submodule of
$H^c_I(R)$. For this, let $N$ be some simple unit $R[F^e]$--submodule of
$H^c_I(R)$ with root $N_0 = \omega_A \cap N$. For each maximal ideal, $R_m
\tensor N$ is a unit $R_m[F^e]$--submodule of $H^c_I(R_m)$ with root $R_m
\tensor N_0$. Since $N_0 \subseteq \omega_A$ is torsion free we conclude
that $R_m \tensor N_0$ is nonzero, and therefore $R_m \tensor N$ is
nonzero. Since $H^c_I(R_m)$ has a unique simple unit $R[F^e]$--submodule
$L(A_m,R_m)$, we get $L(A_m,R_m) \subseteq R_m \tensor N$. Therefore
$L(A,R) \subseteq N$. Since $N$ is simple we get $L(A,R) = N$ and thus
$L(A,R)$ is the unique simple $R[F^e]$--submodule of $H^c_I(R)$.
\end{proof}
This solves the affine case of the globalization problem. Similarly as in
the case of $L(A,R)$ commuting with completion it can be shown that if we
assume the existence of $L(A,R)$ and its commutation with localization on
maximal ideals, then the parameter test module also commutes with
localization at maximal ideals. This can be attributed to the fact that in
the local case the parameter test module is the unique minimal root of
$L(A,R)$. Therefore, in some sense, Problem \ref{prob.Global} is
equivalent to Problem \ref{prob.TestCommutes}.

From the affine case the general case follows easily. One has to show that
for a multiplicative subset $S \subseteq R$ consisting of the powers of a
single element $s \in R$, the localization $S^{-1}L(A,R)$ is the unique
simple $S^{-1}R[F^e]$--submodule of $H^c_{IS^{-1}R}(S^{-1}R)$. Since
$S^{-1}L(A,R)$ is a nonzero unit $S^{-1}R[F^]$--module (further localizing
to a maximal ideal is nonzero) we get $L(S^{-1}A,S^{-1}R) \subseteq
S^{-1}R[F^e]$ as the former is the unique simple unit $R[F^e]$--module.
After further localization at the maximal ideals of $S^{-1}R$ this
inclusion becomes equality, as both are equal to $L(A_m,R_m)$ by the last
proposition. Thus, if $Y \subseteq X$ is a normal subvariety we get a
global sheaf $L(Y,X) \subseteq H^c_{[Y]}(\Oo_X)$ by gluing the compatible
affine pieces $L(A,R)$ of an open affine covering of $X$.

Globalizing the concept of an $R[F^e]$--module we get an
$\Oo_X[F^e]$--module, i.e. a sheaf of $\Oo_X$--modules $M$, together with
a map $F^*M \to[\theta^e] M$ (cf.\ \cite{Em.Kis}). Such an
$\Oo_X[F^e]$--module $M$ is called \emph{unit} if $\theta^e:\F[e] \to M$
is an isomorphism. In fact, with the appropriate modifications, a large
part of the theory developed in Chapters \ref{chap.Fmod} and
\ref{chap.DRsub} carries through to the global case (cf. \cite{Em.Kis}).
With such a global theory available one can easily show that the so
constructed $L(R,A)$ is indeed the unique simple unit
$\Oo_X[F^e]$--submodule of $\Hh^c_{[Y]}(\Oo_X)$. Thus we showed the
following.
\begin{theorem}
    Let\/ $Y \subseteq X$ be a normal, irreducible and closed subvariety of
    codimension $c$ in the smooth variety\/ $X$. Under the assumption that the
    parameter test module commutes with localization at maximal ideals,\/
    $\Hh^c_{[Y]}(\Oo_X)$ has a unique simple unit\/ $\Oo_X[F^\infty]$--submodule\/ $L(Y,X)$.
\end{theorem}

The last issue to deal with now, is to see whether $L(Y,X)$ is
simple as a $\Dd_X$--module, as well. In the affine case, this can
be seen by extending the perfect field $k \subseteq R$ to an
algebraically closed field $K$ of huge cardinality. Similarly as
in the proof of Main Theorem 2 and Lemma \ref{lem.LlforFieldEx},
and using Main Theorem 1, one shows that $K \tensor_k L(A,R) =
L(K\tensor_k A, K \tensor_K R)$ is the unique simple $D_{K
\tensor_k} R$--submodule of $H^c_I(K \tensor_k R)$. Therefore,
$L(A,R)$ is the unique simple $D_R$--submodule of $H^c_I(R)$. Now,
if $Y\subseteq X$ is as before, let $N$ be a proper
$\Dd_X$--submodule $N \subseteq L(Y,X)$. Then for some open affine
subset $U$ of $X$ we have $N|_U$ is a proper $\Dd_U$--submodule
$L(Y,X)|_U$ of $M$. Thus, $N|_U$ is zero since $L(Y,X)|_U = L(Y
\cap U, U)$ is $\Dd_U$--simple. Therefore $N$ is supported on the
closed set $Z=X \setminus U$. Let $z \in Z$ be a point where the
stalk of $N$ is nonzero and $V$ an open affine neighborhood of
$z$. Then $N|_V$ is a nonzero $\Dd_V$--submodule of $L(Y,X)|_V =
L(Y \cap V, V)$. Since the latter is $\Dd_V$--simple, $N|_V =
L(Y,X)|_V$. By further restricting to $(U \cap V)$ we get a
contradiction. On the hand $N|_{U \cap V}=(N|_V)|_U$, and on the
other hand $N|_{U \cap V}=(N|_U)|_V$ is zero. Therefore $N=0$ and
$\Ll(Y,X)$ is $\Dd_X$--simple. We get the following corollary.

\begin{corollary}
    Let\/ $Y \subseteq X$ be a normal, irreducible and closed subvariety of
    codimension\/ $c$ in the smooth variety\/ $X$. Under the assumption that the
    parameter test module commutes with localization at maximal ideals,
    the unique simple\/ $\Oo_X[F^\infty]$--submodule\/ $L(Y,X)$ of\/
    $\Hh^c_{[Y]}(\Oo_X)$ is also its unique simple\/ $\Dd_X$--submodule\/ $\Ll(Y,X)$.
\end{corollary}

Thus we have reduced the problem of proving the existence and construction
of the unique simple $\Dd_X$--submodule of $\Hh^c_{[Y]}(\Oo_X)$ to the
problem of showing that the parameter test module commutes with
localization at maximal ideals; \ie we have to show that if $A$ is a
normal domain, essentially of finite type over the field $k$, then for
every maximal ideal $m$ of $A$ we have $A_m \tensor \tau_{\omega_A} =
\tau_{\omega_{R_m}}$.

\section{Relation to Characteristic zero}\label{sec.relcharzero}
As important as the globalization issue is the question about the relation
of the finite characteristic $L(A,R)$ with the characteristic zero
original. To make this precise we will have to introduce some notation
from the technique of ``reduction to finite characteristic''. We will
treat the affine case here; the general case is obtained by gluing.

Let $R$ be a finitely generated algebra over a field $k$ of characteristic
zero. We replace $k \to R$ by a flat map $T \to R_T$ where $T \subseteq k$
is a finitely generated algebra over $\ZZ$ such that $R_T \tensor_T k =
R$. This can always be achieved. If $R=\frac{k[\xn]}{(\xc[f])}$ is a
presentation of $R$, then we take $T$ to be the $\ZZ$--algebra generated
by the coefficients of the $f_i$'s and $R_T=\frac{T[\xn]}{(\xc[f])}$. To
ensure the flatness of $T \to R_T$, we localize $T$ at a single element.
Furthermore, we freely enlarge $T$ to make any finite amount of data
(finitely generated modules, ideals, maps\ldots) also defined, and free
over $T$.

This family $T \to R_T$ is our data to perform the reduction to finite
characteristic. For a closed point $s \in \Spec T$ denote the residue
field at this point by $\kappa(s)$. Then $R_{\kappa(s)} = \kappa(s)
\tensor_T R_T$ is a finitely generated algebra over the finite (and thus
perfect) field $\kappa(s)$ of finite characteristic $p$. Thus, the closed
fibers are the finite characteristic models of the family. The general
fiber $R_Q = Q \tensor_T R_T$ is the characteristic zero model and is
obtained from $R_T$ by tensoring with the field of fraction $Q$ of $T$.

As an example, let $A=R/I$ be a finitely generated $k$--algebra. We chose
some descent data (\ie the family $T \to R$) such that $I_T$ is defined
and free over $T$. Furthermore, we make sure that all the modules of the
\Czech complex calculating the local cohomology module $H^c_{I_T}(R_T)$
are also free. In particular, all the maps of the \Czech complex are split
(over $T$) and therefore tensoring over $T$ commutes with taking
cohomology. We get $H^c_{I_T}(R_T) \tensor_T \kappa(s)=
H^c_{I_{\kappa(s)}}(R_{\kappa(s)})$ and $H^c_{I_T}(R_T) \tensor_T Q =
H^c_{I_Q}(R_Q)$. This shows that, for almost all primes $s$ of $T$, the
local cohomology modules $H^c_{I_{\kappa(s)}}(R_{\kappa(s)})$ are obtained
by reduction mod $p$ from the characteristic zero model $H^c_I(R_T)$.

The first question to answer is whether the Brylinski--Kashiwara
$D_R$--module $\Ll(A,R)$ arises from a $D_{R_T}$--submodule $\Ll(A_T,R_T)$
of $H^c_{I_T}(R_T)$? \ie is there a unique simple $D_{R_T/T}$--submodule
$\Ll(A_T,R_T)$ of $H^c_{I_T}(R_T)$ such that $\Ll(A,R)=k \tensor_T
\Ll(A_T,R_T)$? If $N$ is a finitely generated $R$--module, then we can
just enlarge $T$ such that $N$ is defined by a free presentation
\[
    R^{\oplus n_2} \to[G] R^{\oplus n_1} \to N \to 0.
\]
The map $G$ is then given by a matrix with entries in $A_T$. But, quite
likely $\Ll(A,R)$ is \emph{not} finitely generated over $R$ and thus it is
not clear that $\Ll(A,R)$ is even defined over a essentially finite
algebraic extension of $T$. Therefore we have to use a different line of
reasoning. The idea is that, even though $\Ll(A,R)$ is not finitely
generated as an $R$--module, it is a coherent $D_R$--module (see
\cite{Borel.Dmod}). Therefore we have a presentation of $\Ll(A,R)$ by
finitely generated free $D_{R/k}$--modules, \ie we have
\[
    D_{R/k}^{\oplus n_2} \to[G] D_{R/k}^{\oplus n_1} \to \Ll(A,R) \to 0
\]
for some matrix $G$ with entries in $D_{R/k}$. Now we observe that, since
$k$ is flat over of $T$, the base change property of differential
operators \cite{EGA4} shows that $D_{R/k} = k \tensor_T D_{R_T/T}$. Thus,
after possibly enlarging $T$, we can assume that the entries in $G$ are
elements of $D_{R_T/T}$. Therefore, if we define the $D_R$--module
$L(A_T,R_T)$ as the cokernel of the map $D^{\oplus n_2}_{R_T/T} \to[G]
D^{\oplus n_1}_{R_T/T}$. Then, clearly $k \tensor_T L(A_T,R_T) =
\Ll(A,R)$, and $L(A_T,R_T)$ is the unique simple $D_{R_T/T}$--submodule of
$H^c_{I_T}(R_T)$. With this setup we pose the following question:
\begin{problem}\label{prob.red}
    Let\/ $R$ be a regular\/ $k$--algebra of finite type, and\/ $A=R/I$ be a normal
    domain. Let\/ $T \to R_T$ be descend data for reduction to finite
    characteristic such that\/ $\Ll(A,R) = k \tensor_T L(A_T,R_T)$. Is the
    reduction of\/ $L(A_T,R_T)$ to finite characteristic the unique simple
    $D_R$--submodule of\/$H^c_{I_{\kappa(s)}}(R_{\kappa(s)})$ for
    infinitely many primes\/ $s$ of\/ $T$? Phrased differently, is\/ $\kappa(s)
    \tensor_T L(A_T,R_T) = L(\kappa(s) \tensor A_T, \kappa(s) \tensor
    R_T)$?
\end{problem}

To investigate the behavior of simple $D_{R/k}$--modules under reduction
to finite characteristic a general study of $D_{R/k}$--modules under
reduction to finite characteristic is at order. In \cite{vdB.Smith} van
den Bergh and Smith development some techniques to reduce rings of
differential to finite characteristic. They are foremost concerned with
the reduction of the ring of differential operators of certain singular
rings, and encounter serious difficulties. Nevertheless, their setup
indicates that for smooth $k$--algebras the ring of differential operators
reduces nicely, as we already used above. Refining their techniques and
extending them to $D_{R/k}$--modules one should be able to develop a
satisfactory theory of reduction to characteristic $p$ for
$D_{R/k}$--modules over a smooth $k$--algebra, hopefully leading to a
positive answer of the above Problem \ref{prob.red}.

%Finally we want to point out that if the intersection homology
%$D_{R/k}$--module reduces to finite characteristic as just explained, then
%our results in Section \ref{sec.FratSimpl} imply (via reduction to finite
%characteristic), that if $A=R/I$ has only rational singularities, then
%$H^c_I(R)$ is $D_R$--simple.

\section{Further perspectives}
In this final section we collect some thoughts about the results obtained,
and techniques used in this dissertation.

\subsection{Endomorphisms of simple $D_R$--modules}
We begin with restating the question raised in Chapter
\ref{chap.DRsub} about the endomorphisms of simple $D_R$--modules in
finite characteristic.
\begin{problem}
Let $R$ be a regular ring, essentially of finite type over a field $k$.
Let $N$ be a simple $D_R$--module. Is $\End_k(N)$ algebraic over $k$.
\end{problem}
A positive answer to this question would at once remove the awkward
assumption of uncountability of the field $k$ in Main Theorem 1. Thus it
would lead to a more natural treatment of the remaining results in this
dissertation, since the repeated occurrence of the step ``extending the
field $k$ to be huge'' could be avoided.

The corresponding result in characteristic zero is a consequence or
Quillen's lemma \cite{Quill}. It is proven by using a almost commutative
variant of the generic freeness lemma. For this the fact that $D_R$ has an
commutative associated graded ring which is a finitely generated
$k$--algebra is crucial. Once more this is not the case in characteristic
$p > 0$. Therefore, Quillen's proof does not have an immediate analog in
finite characteristic. In order to proof the above, the Frobenius
techniques developed here (and used before in
\cite{SmithSP.diffop,SmithSP:DonLine,Haastert.DiffOp,Haastert.DirIm,Bog.DmodBorel})
are likely to be useful. But, at the moment, they don't seem to offer an
immediate solution. It seems that some novel technique must be developed.

%One of the problems in describing the endomorphism set of a simple
%$D_R$--module $N$ via the filtration $\D[e]_R$ of $D_R$ is the following:
%Even if we could control $\End_{\D[e]_R}(N)$ it is not clear how to draw
%conclusions about their intersection, which is $\End_{D_R}(N)$. One
%exceptional case where this method is successful is $N=R$. Then,
%$\End_{D_R}(R) = \bigcap\End{\D[e]_R}(R) = \bigcap R^{p^e} = k$, if $k$ is
%the perfect coefficient field of $R$.

\subsection{Berthelot's theory of arithmetic $D_R$--modules}
As we already discussed in the introduction it is very likely that
Berthelot's theory of arithmetic $D$--modules
\cite{Ber.OpDiff,Ber.FrobDesc,Ber.IntroDmodArith} offers an alternative
take on the existence of the unique simple $D_R$--submodule of $H^c_I(R)$.
As the existence of $\Ll(A,R)$ in characteristic zero follows fairly
formal \cite{Bry.Kash} using the duality functor on holonomic
$D_R$--modules, we can hope that, if such a formal setup is available also
in finite characteristic, then the existence of $\Ll(A,R)$ follows
similarly as in characteristic zero. Berthelot's Frobenius descent is
likely to provide such formalism. He establishes the required operations
of pullback, pushforward and duality of a $D_R$--module in finite
characteristic.

Furthermore, if the existence of the unique simple $D_R$--submodule of
$H^c_I(R)$ can be established along these lines it is automatically a
global object. Thus an investigation of $\Ll(A,R)$ in the light of
Berthelot's theory is likely to provide insight into the question about
globalization raised earlier.

\subsection{Emerton and Kisin's Riemann--Hilbert Correspondence}
Emerton and Kisin describe in \cite{Em.Kis,Em.Kis2} a certain
correspondence for unit $R[F^e]$--modules in finite characteristic, not
unlike the Riemann--Hilbert correspondence in characteristic zero.
Roughly, they show an anti-equivalence between the category of finitely
generated unit $\Oo_X[F^e]$--modules on one hand, and the category of
constructible $\FF_{p^e}$--sheaves on the \'etale site on the other hand.
Their correspondence is on the level of the respective derived categories.
We ask to find the complex of constructible \'etale sheaves which belongs
to $\Ll(Y,X)$ under this correspondence. In the characteristic zero theory
and under the Riemann--Hilbert correspondence this was the intersection
homology complex, an extremely important object.

Emerton and Kisin's correspondence is very new and a study of the objects
belonging to $\Ll(Y,X)$ is likely to improve the understanding of their
correspondence as well as it will construct a potentially significant
complex of constructible sheaves on the \'etale site. This will also
improve the understanding of the nontrivial $t$--structure on the derived
category of \'etale $\FF_{p^e}$--sheaves by describing a class of simple
such complexes. In the characteristic zero case this would correspond to
the fact that the perverse sheaves are the image of the category of
holonomic $D$--modules under the Riemann--Hilbert correspondence.

\subsection{Multiplier versus test ideals}
Our construction of $\Ll(A,R)$ connects it intimately to the parameter
test module. If $A$ is Gorenstein, the parameter test module is identified
with the test ideal. It was shown by Smith in \cite{Smith.multipl} that,
under reduction to characteristic zero the \idx{multiplier ideal} reduces
to the test ideal in the Gorenstein case. The multiplier ideal is an
important invariant in modern day algebraic geometry. Our connection
between the test ideal and the analog of the intersection homology
$D$--module in finite characteristic suggests, that there is also a
connection between the multiplier ideal and the intersection homology
$D$--module in characteristic zero. An investigation into this question is
likely to uncover subtle connections between multiplier ideals and
$D$--module techniques. Furthermore, since it is known that the multiplier
ideal behaves well under reduction to finite characteristic, it will also
improve the understanding of behavior of the intersection homology
$D$--module to finite characteristic. Thus it will be beneficial for the
solution of the problem raised earlier on the connection between
$\Ll(A,R)$ in characteristic zero and $\Ll(A,R)$ in finite characteristic.

%---------   Backmatter ----------------------------------------------------------------

\backmatter   % formatting for the back

\newcommand{\etalchar}[1]{$^{#1}$}
\providecommand{\bysame}{\leavevmode\hbox to3em{\hrulefill}\thinspace}

\printindex  %prints the Index
% \doublespacing

\end{document}